\documentclass{article}
\usepackage{float}
\usepackage{tikz, tikz-3dplot}
\usetikzlibrary{calc,arrows.meta}
\usepackage{graphicx} 
\usepackage{subcaption}
\usepackage{amsmath}
\usepackage{amsfonts}
\usepackage{amsthm}
\usepackage{xcolor}
\newtheorem{theorem}{Theorem}
\newtheorem{proposition}{Proposition}
\usepackage{cases}
\usepackage{comment}
\usepackage{hyperref}
\usepackage{tabularx}
\usepackage{pgfplots}
\pgfplotsset{compat=1.18}
\usepackage{authblk}
\usepackage[margin=1in]{geometry}
\usepackage[backend=biber, maxnames=99, minnames=99]{biblatex}
\usepackage[para]{footmisc}
\newtheorem{remark}{Remark}
\renewcommand{\d}{\text{d}}
\renewcommand{\dfrac}[2]{\frac{\d #1}{\d #2}}
\newcommand{\pdfrac}[2]{\frac{\partial #1}{\partial #2}}
\newcommand{\brackets}[1]{\left\langle #1 \right\rangle}
\newcommand{\e}{\varepsilon}
\newcommand{\R}{\mathbb{R}}

\renewcommand{\u}{u}
\newcommand{\rhobar}{\overline{\rho}}
\newcommand{\ul}[1]{\u_{#1}}
\newcommand{\fraction}[1]{\left(\frac{\, \rho_{#1} \, }{\rhobar}\right)}
\newcommand{\uzerobar}{\overline{u}_0}
\newcommand{\rhozerobar}{\overline{\rho}_0}
\newcommand{\ebar}{\overline{\e}}

\counterwithin{equation}{section}

\title{An Analysis of the Riemann Problem for a $2 \times 2$ System of Keyfitz-Kranzer Type Conservation Laws Using Shadow Waves and Dafermos Regularization}
\date{\today}
\author[1]{Josh Culver \thanks{culver23@purdue.edu}}
\author[2]{Aubrey Ayres \thanks{aayres@email.sc.edu}}
\author[3]{Evan Halloran\thanks{ehallor@iu.edu}}
\author[4]{Ryan Lin \thanks{rl00023@mix.wvu.edu}}
\author[5]{Emily Peng \thanks{emily.peng@yale.edu}}
\author[4]{Charis Tsikkou \thanks{tsikkou@math.wvu.edu}}
\affil[1]{Department of Mathematics, Purdue University, West Lafayette, IN 47907, USA}
\affil[2]{Department of Mathematics, University of South Carolina, Columbia, SC 29208, USA}
\affil[3]{Department of Mathematics, Indiana University, Bloomington, IN 47405, USA}
\affil[4]{School of Mathematical and Data Sciences, West Virginia University, Morgantown, WV 26506, USA}
\affil[5]{Department of Physics, Yale University, New Haven, CT 06511, USA}

\usepackage{float}
\begin{document}
\maketitle
\begin{abstract}\noindent
We consider a system of two conservation laws and provide a detailed description of both classical and non-classical self-similar Riemann solutions. In particular, we demonstrate the existence of overcompressive delta shocks as singular limits of the Dafermos regularization of the system. The system is chosen for its minimal yet representative structure, which captures the essential features of transport dynamics under density constraints. Our analysis is carried out using blow-up techniques within the framework of Geometric Singular Perturbation Theory (GSPT), allowing us to resolve the internal structure of these singular solutions. Despite its simplicity, the system serves as a versatile prototype for crowding-limited transport across a range of applications, including biological aggregation, ecological dispersal, granular compaction, and traffic congestion. Our findings are supported by numerical simulations using the Local Lax–Friedrichs scheme.
\end{abstract}

\vspace{5mm}


\vspace{5mm}


\section*{Notation}
We define $\left[\cdot\right] := \left[\cdot\right]_{\text{jump}} = \cdot_L - \cdot_R,$ $\cdot_- - \cdot_+$or $\cdot_0 - \cdot_1$.

\section{Introduction} \label{intro}
In this work, we study a nonlinear system of two conservation laws that model density-dependent transport without external forcing (for the case where a time-dependent external force is incorporated, see Culver et al. \cite{REU2025_1}). The system takes the form:
\begin{align}
        \rho_t + \left( \rho u \left(1 - \fraction{}^a \right) \right)_x =& \, \rho_t + \bigl(f_1(\rho,\u)\bigr)_x = 0, \label{eq:sys_1} \\
        (\rho u)_t + \left(\rho u^2 \left( 1 - \fraction{}^a \right) \right)_x =& \, (\rho \u)_t + \bigl(f_2(\rho,\u)\bigr)_x = 0, \label{eq:sys_2}
\end{align}
where the dependent variables are density $\rho$ and velocity $u.$ The constant $\bar{\rho}>0$ is a critical density (for example, a jamming or saturation threshold) and $a\in \mathbb{R} \setminus \{0\}$ governs the nonlinear degeneracy of the mobility term. The sign of $a$ has significant modeling and analytical consequences:
For $a>0$, mobility vanishes as $\rho \to \bar{\rho}$, modeling crowding or jamming effects due to high density;
for $a<0$, degeneracy occurs as $\rho \to 0$, modeling low-density inhibition such as clustering or dispersal suppression.

Although this system has a simple structure and does not include complex effects, it serves as a versatile prototype for transport processes affected by crowding or depletion effects, with applications in biological aggregation, ecological dispersal, traffic flow, and granular compaction. A key analytical feature of the model is the loss of strict hyperbolicity. The characteristic speeds of the system are 
$$\lambda_{a} = \, \u \left(1 - (a+1)\fraction{}^a\right), \quad \lambda_{0} = \, \u \left(1 - \fraction{}^a\right),$$
and, depending on the state, can coincide or vanish. In particular, degeneracy occurs along the critical density surface $\rho=\bar{\rho}$, where $\lambda_0 = 0$, along $\rho=\bar{\rho}\cdot(\frac{1}{a+1})^{1/a} $ when $a>-1$, where $\lambda_a=0$, partially degeneracy occurs at $\rho = 0$ when $a > 0$, where the system loses strict hyperbolicity and the two eigenvalues are equal, and fully degeneracy occurs at $u=0$,  where both eigenvalues vanish. This breakdown of strict hyperbolicity gives rise to novel behaviors. We construct and classify self-similar Riemann solutions, capturing both classical wave patterns (contact discontinuities, rarefactions, and shocks) and non-classical ones, such as overcompressive delta shocks that arise as bona fide limits of smooth solutions to the Dafermos regularization and possess internal structure that can be resolved using blow-up techniques from Geometric Singular Perturbation Theory. In addition to the singular solutions mentioned above, which are the main focus of this work, the Riemann problem also admits solutions that pass through vacuum states ($\rho = 0$) and across the critical density threshold $\rho = \bar{\rho}$, where mobility vanishes and one characteristic speed degenerates. A comprehensive treatment of all wave patterns is presented in Culver et al. \cite{REU2025_1}. Our analytical results are validated through numerical simulations using the local Lax–Friedrichs scheme, which confirms the emergence of singular profiles, vacuum layers, and degenerate wave transitions predicted by the theory.

The fundamental building block for establishing existence theorems for the Cauchy problem in systems of conservation laws (in one spatial dimension) is the solution to the Riemann problem. As such, constructing explicit Riemann solutions is essential from both analytical and numerical perspectives.

Singular solutions involving delta shocks, a more compressive generalization of classical shock waves, are characterized by the concentration of at least one state variable in a weighted Dirac delta distribution. Keyfitz and Kranzer first discovered delta shocks \cite{Ke_Kr_1, Ke_Kr_2,Ke_Kr_3} in the context of a strictly hyperbolic, genuinely nonlinear system derived from a one-dimensional model of isothermal gas dynamics. They identified a large region of the state space where the Riemann problem cannot be resolved using classical solutions composed solely of shocks and rarefactions. In this regime, the authors constructed approximate unbounded solutions that do not satisfy the conservation laws in the classical weak sense and observed that only the first component of the Rankine–Hugoniot condition is satisfied, resulting in a unique propagation speed $\eta$ for which any two given states can be connected.

Subsequently, Schecter \cite{Sc} provided a rigorous justification for delta shock formation by constructing self-similar viscous profiles using techniques from geometric singular perturbation theory. In particular, Schecter used a blow-up method to overcome the lack of normal hyperbolicity, based on the foundational work of Fenichel \cite{Fe} and Jones \cite{Jo}. For further developments and related approaches to singular solutions, see also \cite{Hsu, Ka_Mi, Ke, Ke_2, Ke_3, Ke_4, Le_Sl, Ma_Be} and the references therein.

Previous investigations of singular solutions have primarily focused on cases where only a single state variable concentrates into a Dirac delta distribution. However, physically relevant systems exist in which both state variables simultaneously develop singularities. For example, Mazzotti et al. \cite{Ma_1, Ma_2, Ma_3} numerically studied a model arising in two-component chromatography, which has significant industrial applications. This system exhibits singular solutions in which both components of the conserved variables concentrate and neither component satisfies the classical Rankine–Hugoniot condition. Tsikkou \cite{Ts} provided a detailed analysis of the same chromatography system, which undergoes a change of type (from hyperbolic to elliptic) in parts of the state space. By applying linear transformations to the conserved quantities, a simplified system was derived and studied, offering a coherent explanation for the unbounded solutions observed. In addition, Frew et al. \cite{REU2024} studied a system of two balance laws of Keyfitz–Kranzer type involving a generalized Chaplygin gas with time-dependent forcing. The authors provided a detailed explanation and classification of the resulting non-classical singular solutions. Notably, the explicit time dependence of the source term leads to a dynamic partitioning of the state space, resulting in temporal variations in the regions where combinations of classical and non-classical singular solutions are admissible. More generally, time-dependent and forced systems have increasingly become a focus of study for singular solutions, particularly in contexts where classical wave structures fail; see also Culver et al. \cite{REU2025_1}. These studies demonstrate that singular solutions are not limited to autonomous or conservative systems, but can also play a fundamental role in time-evolving forced systems with applications in gas dynamics, cosmology, and reactive flows.

These examples raise several natural questions. How can one predict the emergence of singular solutions? What is their physical interpretation? In what precise mathematical sense do they satisfy the original conservation laws? Can one formulate a more general definition of a solution that captures a broader class of singular behaviors, possibly extending beyond classical or weak solutions? And how are these phenomena related to structural features of the system, such as genuine nonlinearity or non-hyperbolicity?

The model considered in this paper contributes significantly to these questions. In addition to its physical relevance, it provides a framework for analyzing nonclassical and singular structures that arise in degenerate or non-strictly hyperbolic systems. From a mathematical point of view, this work seeks to expand the understanding of how singular solutions can serve as essential building blocks in the global resolution of Riemann and Cauchy problems with large initial data, potentially leading to new solution theories and generalized numerical schemes.

The system consisting of equations \eqref{eq:sys_1} and \eqref{eq:sys_2} is a special case of a generalized Keyfitz-Kranzer-type system 
\begin{equation}\label{EQMODEL}
    \begin{split}
    \begin{cases}
    \rho_t+\bigg(\rho \Phi\big(\rho, u\big)\bigg)_x=F\big(\rho, u\big),\\
    \big(\rho u\big)_t+\bigg(\rho u\Phi\big(\rho, u\big)\bigg)_x=G\big(\rho, u\big),
    \end{cases}
    \end{split}
\end{equation}
and has various applications depending on $\Phi,$ $F,$ and $G.$ For example, the pressureless Euler system and a simplified second-order traffic flow model structurally related to simplified versions of the Aw and Rascle \cite{Aw} and Zhang \cite{Zh_3} models correspond to $F=G=0, \ \Phi(\rho, u)=u,$ and $F=G=0, \ \ \Phi(\rho, u)=u-p(\rho),$ respectively. The literature, such as in Zhang \cite{Zh_1,Zh_2}, shows that the Riemann problem with pressure laws depending solely on density and with $F=0$ has been extensively studied. In this work, we consider equations \eqref{eq:sys_1} and \eqref{eq:sys_2} and describe solutions to the Riemann problem. The initial data are given by
\begin{align}
    (\rho, \u)(x,0) = \begin{cases}
        (\rho_L,u_L)& x < 0, \\
        (\rho_R,\u_R)& x > 0,
    \end{cases} \label{eq:initial_conditions}
\end{align}
which represents a jump discontinuity between two constant states in a strictly hyperbolic region.

The paper is organized as follows: In Section \ref{classical_analysis}, we provide a formal description of the classical Riemann solutions to the system of conservation laws. We derive the shock and contact discontinuity curves through a given left state using the Rankine–Hugoniot conditions, and construct the rarefaction curves based on the system's eigenvalues and eigenvectors. In Section \ref{delta_shocks}, we verify that singular solutions satisfy equations \eqref{eq:sys_1} and \eqref{eq:sys_2} in the sense of distributions. There are two standard approaches to this verification. One is to directly postulate a solution involving Dirac delta distributions and substitute it into the weak formulation using test functions. The other, known as the shadow wave method, was introduced by Marko Nedeljkov \cite{Daw-Marko, Marko, Marko2} and involves the construction of a family of smooth approximate solutions with narrow internal layers that concentrate as singular limits. In Section \ref{numerics}, we present state space plots for various values of $a\in \mathbb{R} \setminus \{0\}$, and identify the regions in which the Riemann problem cannot be solved using only shocks, contact discontinuities, and rarefactions. In Section \ref{GSPT}, we apply the method of Geometric Singular Perturbation Theory (GSPT) combined with blow-up techniques to resolve the internal structure of the unbounded self-similar profiles that arise from the Dafermos regularization. Finally, in Section \ref{conclusion}, we summarize our main findings and outline possible directions for future work.

\section{Classical Waves: Contact Discontinuities, Shocks, and Rarefactions} \label{classical_analysis}
In this section, we start by calculating the eigenvalues of the system and examining how the property of strict hyperbolicity varies with the state. Next, we determine the Hugoniot locus, which is the set of points in state space that can be connected to a fixed left state by either a shock wave that meets the Lax admissibility condition or a contact discontinuity. We also identify the $a$-rarefaction curve, representing the set of points in the state space that can be connected through an $a$-rarefaction. In Section \ref{numerics}, we plot the Hugoniot locus and discuss various combinations of classical wave patterns.

\subsection{Hyperbolicity, Genuine Nonlinearity, and Linear Degeneracy} \label{RH_analysis}
We will eventually separate cases based on the value of $a$, but initial calculations can generally proceed. We start by writing the system of \eqref{eq:sys_1} and \eqref{eq:sys_2} in vector form as $\partial_t H + \partial_x G = 0$, where $H(\rho, \u, x, t) = \left(\rho \enspace \rho \u\right)^T$ and 
$$G = \begin{pmatrix}
    \rho \u - \rho \u \fraction{}^a \\
    \rho \u^2 - \rho \u^2 \fraction{}^a 
\end{pmatrix}.$$
For a general vector $ \left(B_1 \enspace B_2\right)^T$  we define
$$D\begin{pmatrix}B_1 \\ B_2\end{pmatrix} = \begin{pmatrix}
    \pdfrac{B_1}{\rho} & \pdfrac{B_1}{\u} \\ \\
    \pdfrac{B_2}{\rho} & \pdfrac{B_2}{\u}
\end{pmatrix}.$$

To determine whether our system is hyperbolic, we seek eigenvalues and eigenvectors, that is, $\lambda$ and $R$ such that $(DG - \lambda DH)R = 0$. Note that 
\begin{align}
    DG - \lambda DH = \begin{pmatrix}
    \u - (a+1) \u \fraction{}^a - \lambda & \rho - \rho \fraction{}^a \\ \\
    \u^2 - (a+1) \u^2 \fraction{}^a - \lambda \u & 2\rho \u - 2\rho \u \fraction{}^a - \lambda \rho
\end{pmatrix}.
\end{align}

Solving $\det(DG - \lambda DH) = 0$ yields
$$0 = \left(\lambda + \u \left((a+1)\fraction{}^a - 1\right)\right)\left(\lambda + \u \left(\fraction{}^a - 1\right)\right).$$
As we will see, the order of these two eigenvalues changes according to the sign of $a\cdot u$ ( $\lambda_a < \lambda_0$ if $a\cdot u > 0$; both quantities are defined precisely below), so we label them to denote their form rather than ``one'' and ``two'' to denote their order. Thus,
\begin{align}
    \lambda_{a} =& \, \u \left(1 - (a+1)\fraction{}^a\right), \label{eq:eigen_a} \\
    \lambda_{0} =& \, \u \left(1 - \fraction{}^a\right). \label{eq:eigen_0}
\end{align}
We note that the eigenvalues are real, $\lambda_a \neq \lambda_0$ for $a \neq 0$ and neither are identically zero when $u \neq 0,$ $\rho \neq \bar{\rho},$ $\rho\neq \bar{\rho}\cdot(\frac{1}{a+1})^{1/a}$ (when $a>-1$), and $\rho \neq 0$; thus, our system is strictly hyperbolic. After a short calculation, we obtain the corresponding eigenvectors,
\begin{align}
    R_{a} =& \begin{pmatrix} 1 \\ \\ 0 \end{pmatrix}, \label{eq:eigenvector_a} \\ 
    R_0 =& \begin{pmatrix}
        \rho\left(1 - \fraction{}^a\right) \label{eq:eigenvector_0} \\ \\
        a \u \fraction{}^a
    \end{pmatrix}.
\end{align}

Finally, observe that
\begin{align}
    D\lambda_{a} \cdot R_{a} =& -\u \frac{(a+1)a \rho^{a-1}}{\rhobar^a}, \\
    D\lambda_0 \cdot R_0 =& -u\frac{a\rho^a}{\rhobar^a} \left(1 - \fraction{}^a\right) + a\fraction{}^a \u \left(1 - \fraction{}^a\right) = 0.
\end{align}
Hence, the $a$-characteristic family is genuinely nonlinear unless $a = -1$ or $u=0$, and the $0$-characteristic family is linearly degenerate.

As a quick digression, we note the already common appearance of $a$ and $a + 1$. These hint that our cases will separate using the sign of these two values: $a < -1$, $a = -1$, $-1 < a < 0$, and $0 < a$.

\subsection{Rankine-Hugoniot Jump Conditions and Related Calculations}
We seek conditions on pairs of points in the $(\rho,\u)-$plane such that there exists a shock wave propagating at speed $s$ and satisfying the Rankine–Hugoniot jump conditions, which connects the two states. However, because we have a system of two variables, we consider both equations \eqref{eq:sys_1} and \eqref{eq:sys_2} and require $s = x'(t), x: \mathbb{R} \rightarrow \mathbb{R}$ such that
\begin{align}
    \begin{cases}
        s[\rho] = \left[\rho \u - \rho \u \fraction{}^a\right], \\
        s[\rho \u] = \left[\rho \u^2 - \rho \u^2 \fraction{}^a\right].
    \end{cases} \label{eq:RH_conditions}
\end{align}
Assuming $[\rho] \neq 0$ and $[\rho \u] \neq 0$, this means that, given a left state $(\rho_L, \u_L)$, we try to find all pairs $(\rho_R, \u_R)$ such that
\begin{align*}
    &\frac{\rho_L \u_L - \rho_L \u_L \fraction{L}^a}{\rho_L - \rho_R} - \frac{\rho_R \u_R - \rho_R \u_R \fraction{R}^a}{\rho_L - \rho_R} \\
    &\quad = \frac{\rho_L \u_L^2 - \rho_L \u_L^2 \fraction{L}^a}{\rho_L \u_L - \rho_R \u_R} - \frac{\rho_R \u_R^2 - \rho_R \u_R^2 \fraction{R}^a}{\rho_L \u_L - \rho_R \u_R}.
\end{align*}
Although this is algebraically impossible in general, we show that by considering certain curves of the form $u = f(u_L, \rho_L, \rho)$ one obtains either a shock or a contact discontinuity for each characteristic family. This allows us to identify all shock and contact discontinuity curves in the $(\rho,u)$-plane.

From this point forward, we denote the right state variables $\rho_R$ and $u_R$ as $\rho$ and $u,$ since we view them as variable points in the state space, with $\rho_L$ and $\u_L$ fixed.

\subsubsection{Shock and Rarefaction Curves of the \texorpdfstring{$a$}{a}-Family: \texorpdfstring{$a \neq -1$}{a != -1}}
Consider $\u_L = \u \neq 0$ and $a \neq -1$. Then \eqref{eq:RH_conditions} is satisfied by $\rho_L$ and $\rho$, and
\begin{align}
    s_a =& \, \u_L \left\{1 - \frac{\left(\rho_L \fraction{L}^a - \rho \fraction{}^a\right)}{\rho_L - \rho}\right\}.
\end{align}
The $a$-shock satisfies the Lax shock admissibility criterion if
\begin{align}
    \lambda_a(\rho, \u) < s_a < \lambda_a(\rho_L, \u_L), \label{eq:Lax_condition}
\end{align}
that is, the shock speed must be between the characteristic speeds on either side of the discontinuity. The condition is equivalent to requiring that characteristics of the $a$-family enter the shock from both sides, ensuring that the shock is compressive. This is equivalent to
\begin{align}
    \begin{cases}
        (a+1)\rho^a > \frac{\left(\rho_L^{a+1} - \rho^{a+1}\right)}{\rho_L - \rho} > (a+1)\rho_L^a & \quad \text{if } \u_L = \u > 0, \\ \\
        (a+1)\rho^a < \frac{\left(\rho_L^{a+1} - \rho^{a+1}\right)}{\rho_L - \rho} < (a+1)\rho_L^a & \quad \text{if } \u_L = \u < 0.
    \end{cases}
\end{align}

Further calculations allow us to obtain Table \ref{tab:shocks_a_family}.
\begin{table}[H]
    \centering
    \begin{tabular}{ | m{2.75cm} | m{2.25cm} | m{2.25cm} | m{2.25cm} | }
        \hline
        & $a < -1$ & $-1 < a < 0$ & $0 < a$ \\
        \hline
        $\u_L = \u > 0$ & $S_a$ to the right & $S_a$ to the left & $S_a$ to the right \\
        \hline
        $\u_L = \u < 0$ & $S_a$ to the left & $S_a$ to the right & $S_a$ to the left \\
        \hline
    \end{tabular}
    \caption{Existence of classical shocks of the $a$ family when $\u = \u_L$}
    \label{tab:shocks_a_family}
\end{table}
Notice that, in verifying the Lax shock admissibility criterion for $a \neq -1$, the inequalities are exactly reversed on the opposite side of the left state, i.e., $\lambda_a(\rho, \u) > s_a > \lambda_a(\rho_L, \u_L)$ for some values of $(\rho, \u)$. Noting that the $a$ family is genuinely nonlinear, we expect a rarefaction curve, a self-similar solution that continuously connects the left and right states when the characteristics fan out. We now proceed to derive its explicit form.

Let $\xi = \frac{x}{t}$ (refer to the beginning of Section \ref{GSPT} for details) so that \eqref{eq:sys_1} and \eqref{eq:sys_2} become
\begin{align}
    \bigl(DG - \xi DH) \begin{pmatrix} \dfrac{\rho}{\xi} \\\\ \dfrac{u}{\xi} \end{pmatrix}= 0.
\end{align}
Since this rarefaction is to be of the $a$ family, we set $\xi = \frac{x}{t} = \lambda_a$ and require \(\begin{pmatrix} \dfrac{\rho}{\xi} \\\\ \dfrac{u}{\xi} \end{pmatrix}\) to equal an eigenvector of the $a$-family, i.e., any multiple of the eigenvector $R_a$ from \eqref{eq:eigenvector_a}. Let this multiple be a function $f(\rho,\u)$. Then
$$\begin{pmatrix} \dfrac{\rho}{\xi} \\\\ \dfrac{u}{\xi} \end{pmatrix} = \begin{pmatrix} f(\rho,\u) \\\\ 0 \end{pmatrix} = f(\rho, \u)R_a.$$
Now $F\bigl(\rho(\xi), \u(\xi)\bigr) = \xi + D_1$ and $\u(\xi) = D_2$ for some $D_1, D_2 \in \mathbb{R}$. For this to connect the left state to the right state, we require
\begin{align*}
    F(\rho_L,\u_L) =& \, \lambda_a(\rho_L,u_L) + D_1, \\
    F(\rho,\u) =& \, \lambda_a(\rho,u) + D_1, \\
    u_L =& \, D_2 = \, u.
\end{align*}
This means $F(\rho_L,\u_L) - F(\rho,\u) = \lambda_a(\rho_L,u_L) - \lambda_a(\rho,u)$. Thus, we see that $F(\rho,\u) = \lambda_a(\rho,\u)$  and $D_1 = 0$ suffice. Rarefaction curves, when they exist, are defined by
\begin{align}
    (\rho, \u)(x,t) = (\rho, \u)\left(\frac{x}{t}\right)= \begin{cases}
        (\rho_L, \u_L) & \xi < \lambda_a(\rho_L,\u_L), \\
        \left(\rhobar \left(\frac{u_L-\xi}{(a+1)\u_L}\right)^{\frac{1}{a}}, \u_L\right) & \lambda_a(\rho_L,\u_L) \leq \xi \leq \lambda_a(\rho,\u), \\
        (\rho,\u_L) & \lambda_a(\rho,\u) < \xi.
    \end{cases} \label{eq:rarefaction}
\end{align}

We summarize this information in Table \ref{tab:rarefaction_a_family}.
\begin{table}[H]
    \centering
    \begin{tabular}{ | m{2.75cm} | m{2.25cm} | m{2.25cm} | m{2.25cm} | }
        \hline
        & $a < -1$ & $-1 < a < 0$ & $0 < a$ \\
        \hline
        $\u_L = \u > 0$ & $R_a$ to the left & $R_a$ to the right & $R_a$ to the left \\
        \hline
        $\u_L = \u < 0$ & $R_a$ to the right & $R_a$ to the left & $R_a$ to the right \\
        \hline
    \end{tabular}
    \caption{Existence of rarefactions of the $a$-family when $\u_L = \u_R$}
    \label{tab:rarefaction_a_family}
\end{table}

\subsubsection{The \texorpdfstring{$a$}{a}-Family: \texorpdfstring{$a = -1$}{a = -1}}
We now return to the special case $a = -1$. Note that the shock speed becomes
\begin{align}
    s_a =& \, \u_L \left\{1 - \frac{\left(\rho_L \fraction{L}^a - \rho \fraction{}^a\right)}{\rho_L - \rho}\right\} \notag \\
    =& \, \u_L \left\{1 - \frac{\rhobar - \rhobar}{\rho_L - \rho}\right\}, \notag \\
    \implies s_a =& \, \u_L = \u_R.
\end{align}
Since $a + 1 = 0$ when $a = -1$, we have $\lambda_a(\rho_L,\u_L) = s = \lambda_a(\rho,\u)$, independent of $\rho_L$ and $\rho$. Thus, there is a contact discontinuity of the $a$-family along the line $\u = \u_L$.

\subsubsection{Contact Discontinuities of the \texorpdfstring{$0$}{0}-Family: \texorpdfstring{$\rho_L \neq \rhobar$}{rhoL != rhobar}} \label{contact_disconts}
We now turn to our second special assumption: that the left and right states lie on the following prescribed curve in the $(\rho,u)$-plane
\begin{align}
    \lambda_0(\rho,\u) = \u - \u \fraction{}^a = \u_L - \u_L \fraction{L}^a = \lambda_0(\rho_L,\u_L).
\end{align}
This implies that the Rankine-Hugoniot jump conditions \eqref{eq:RH_conditions} are satisfied. Thus, it defines a contact discontinuity curve of the $0$-family that is discontinuous on the line $\rho = \rhobar$. We choose the branch of this $C_0$ that is continuous and passes through the left state, and name the other branch the ``mirror'' curve $C_{0,m}$, as it will be an important boundary when we analyze the regions later.

Under this assumption, we derive the curve by expressing $u$ as a function of $\rho$, thus characterizing the relationship that defines the admissible states:
\begin{align}
    &\lambda_0(\rho,\u) = \lambda_0(\rho_L,\u_L), \notag \\
    \Leftrightarrow \u &\left( 1 - \fraction{}^a \right) = \u_L - \u_L \fraction{L}^a, \notag \\
    \implies \u &= \frac{1}{\, \rhobar^a - \rho^a \,} \bigl( \u_L \left( \rhobar^a - \rho_L^a \right) \bigr), \notag \\
    \implies \u &= \frac{1}{\, \rhobar^a - \rho^a \,} \bigl( \u_L \left( \rhobar^a - \rho^a + \rho^a - \rho_L^a \right) \bigr), \notag \\
    \implies \u &= \frac{\, \rho^a - \rho_L^a \,}{\rhobar^a - \rho^a} \u_L + \u_L.
\end{align}

As $\rho \rightarrow \infty \text{ or } \rho \rightarrow 0$, $u(\rho)$ will approach either $0$ or an asymptote that will form another boundary in our region analysis. If $a > 0$,
\begin{gather*}
    \lim_{\rho \rightarrow \infty} \u = - \u_L + \u_L = 0, \\
    \lim_{\rho \rightarrow 0} \u = \u(0) = \frac{- \rho_L^a}{\rhobar^a} \u_L + \u_L = \lambda_0(\rho_L,\u_L).
\end{gather*}
On the other hand, if $a < 0$,
\begin{gather*}
    \lim_{\rho \rightarrow \infty} \u = \frac{- \rho_L^a}{\rhobar^a} \u_L + \u_L = \lambda_0(\rho_L,\u_L), \\
    \lim_{\rho \rightarrow 0} \u = - \u_L + \u_L = 0.
\end{gather*}

\subsubsection{The \texorpdfstring{$0$-Family: $\rho_L = \rhobar$}{0-family: rhoL = rhobar}}
We note that our previous derivation no longer applies in the case $\rho_L = \rhobar$. Returning to the original criterion $\lambda_0(\rho,\u) = \lambda_0(\rho_L,\u_L)$ we obtain
\begin{align}
    \u\left(1 - \fraction{}^a\right) = 0.
\end{align}
This implies $\rho = \rhobar$ when $\u \neq 0$, yielding a contact discontinuity $C_0$ along the line $\rho = \rhobar$.\newline

Having characterized the admissible wave patterns, we now examine whether singular solutions can be interpreted as weak solutions to the system. Our goal is to confirm that such solutions satisfy equations \eqref{eq:sys_1} and \eqref{eq:sys_2} in the distributional sense. Singular solutions are expected to arise in certain regimes, and a detailed analysis of their structure and formation will be presented in a later section.

Two principal methods are commonly used for this purpose. The first involves proposing a solution that includes Dirac delta terms and verifying its validity by inserting the expression into the weak formulation via test functions. The second approach, known as the shadow wave method, was developed by Marko Nedeljkov \cite{Daw-Marko, Marko, Marko2}, and is based on constructing a sequence of smooth approximations with sharply localized internal layers that converge to the singular solution in the limit.

\section{Delta Shocks} \label{delta_shocks}
As will be illustrated in a later section, there exist regions of the $(\rho,\u)$-plane where there are no classical solutions nor solutions that pass through vacuum states ($\rho=0$) or across the critical density threshold $\rho=\bar{\rho}.$ This breakdown leads us to seek unbounded solutions that extend the notion of weak solutions in order to accommodate such singular behavior.

\subsection{Delta Shock and Resulting ODEs} \label{delta_ansatz}
One may justify the presence of delta shocks by postulating a solution that includes a Dirac delta distribution and checking whether it satisfies \eqref{eq:sys_1} and \eqref{eq:sys_2} in the distributional sense.

We define a two-dimensional weighted $\delta$-measure $\omega(s)\delta_S$ supported on a smooth curve $\\S = \left\{\bigl(x(s),t(s)\bigr): a \leq s \leq b\right\}$ by
\begin{align*}
    \brackets{\omega(\cdot)\delta_S, \phi(\cdot, \cdot)} = \int_a^b \omega\bigl(t(s)\bigr) \phi\bigl(x(s), t(s)\bigr) \d s
\end{align*}
for all $\phi \in C_c^\infty\bigl(\mathbb{R}\times(0,\infty)\bigr) = C_c^\infty(\mathbb{R}\times\mathbb{R}_+)$.

Motivated by numerical evidence indicating a Dirac delta concentration in $\rho$ alone, we seek a \textit{delta-shock type} solution of the form
$$\u (x,t) = U_0(x,t), \qquad \rho(x,t) = \rho_0(x,t) + \omega(t)\delta(x-x(t)),$$
where
$$ U_0(x,t) = \begin{cases} \u_{L}(x,t) & x < x(t), \\ \u_{\delta}(t) & x = x(t), \\ \u_{R}(x,t) & x > x(t), \end{cases}
\qquad
\rho_0(x,t) = \begin{cases} \rho_{L}(x,t) & x < x(t), \\ \rho_{R}(x,t) & x > x(t),\end{cases}$$
$S = \left\{\bigl(x(t),t\bigr): 0 \leq t < \infty\right\}$, and $\omega \in C^1(\mathbb{R}_+)$.
These are required to satisfy \eqref{eq:sys_1} and \eqref{eq:sys_2} in the sense of distributions, that is,
\begin{align}
    \brackets{\rho_0, \partial_t \phi} + \brackets{f_1(\rho_0, U_0), \partial_x \phi} = 0, \label{eq:weak_sense_solutions1} \\
    \brackets{\rho_0 U_0, \partial_t \phi} + \brackets{f_2(\rho_0, U_0), \partial_x \phi} = 0. \label{eq:weak_sense_solutions2}
\end{align}
Note that all nonnegative smooth or piecewise constant sequences of functions $(f_n)$ which converge to $\delta$ in the distributional sense obey $(f_n^{a+1}) \rightarrow 0$ when $a < 0$. Therefore, the delta-shock substitution is admissible only for $a < 0$ as this ensures that the highest power of $\rho_0$ that appears in the flux term, namely $a+1$, remains less than one.

Let $\dfrac{x}{t} = \u_{\delta}$ so that, using Green's Theorem and the compact support of $\phi$,
\begin{align*}
    0 =& \, \brackets{\rho_0, \partial_t \phi(x,t)} + \brackets{f_1(\rho_0, U_0), \partial_x \phi(x,t)} \\
    =& \int_0^\infty \int_{-\infty}^{x(t)} \left\{\rho_L \partial_t\phi + f_1(\rho_L,\u_L) \partial_x\phi \right\} \d x\d t \\
    & + \int_0^\infty \int_{x(t)}^\infty \left\{\rho_R \partial_t\phi + f_1(\rho_R,\u_R)\partial_x\phi\right\}\d x\d t \\
    & + \int_0^\infty \bigl\{\omega(t) \partial_t\phi(x(t),t) + \omega(t)\u_{\delta}(t) \partial_x\phi(x(t),t) \bigr\}\d t \\
    =& \oint_{x(t)} -\rho_L \phi(x(t),t) \d x + f_1(\rho_L,\u_L) \phi \d t \\
    & + \oint_{-x(t)} -\rho_R \phi(x(t),t) \d x + f_1(\rho_R,\u_R) \phi \d t \\
    & + \int_0^\infty \omega(t) \d\phi \\
    =& \int_0^\infty -\left\{\rho_L - \rho_R\right\} \phi(x(t),t) x'(t) \d t \\
    & + \int_0^\infty \left\{f_1(\rho_L,\u_L) - f_1(\rho_R,\u_R)\right\} \phi(x(t),t) \d t \\
    & + \phi(x(t),t) \omega(t) \Big\vert_0^\infty - \int_0^\infty \phi(x(t),t) \dfrac{}{t}\bigl(\omega(t)\bigr)\d t, \\
    \implies& -[\rho] x'(t) + [f_1(\rho,\u)] = \dfrac{}{t}\bigl(\omega(t)\bigr).
\end{align*}
Using similar calculations on \eqref{eq:weak_sense_solutions2}, we obtain the following result:
\begin{align*}
    -[\rho \u] x'(t) + [f_2(\rho,\u)] = \dfrac{}{t}\bigl( \omega(t)\u_{\delta}(t) \bigr).
\end{align*}

We impose the initial conditions $x(0) = 0, \omega(0) = 0, \text{ and } \u_{\delta}(0) = s_- \in \mathbb{R}$. Thus, we have
\begin{align}
    \dfrac{x}{t} &=u_{\delta}(t), \label{eq:ODE_1.2} \\
    \dfrac{}{t} \bigl( \omega(t) u_{\delta}(t) \bigr) &= -[\rho \u] x'(t) + \left[\rho \u^2 - \rho \u^2 \fraction{}^a\right], \label{eq:ODE_2.2} \\
    \dfrac{}{t} \bigl( \omega(t) \bigr) &= -[\rho] x'(t) + \left[\rho \u - \rho \u \fraction{}^a\right]. \label{eq:ODE_3.2}
\end{align}
This system corresponds exactly to the shadow wave method (discussed in the following section) where $\u_{\delta}(t) = \eta_0(t) = \eta_1(t)$ and $\omega(t) = \zeta_0(t)+\zeta_1(t)$. Using this connection, we can determine the value of $\u_{\delta}(0) = s_-$ in both cases, as we do at the end of Section \ref{marko_anegative}.

\subsection{Nedeljkov's Shadow Wave Method for Singular and Delta Shocks} \label{marko_shadow_wave}
Following \cite{Daw-Marko, Marko, Marko2}, consider the following weighted shadow wave solution to the equations \eqref{eq:sys_1} and \eqref{eq:sys_2},
\begin{align}
    \left(\rho_\e,\u_\e\right) = \begin{cases}
        \bigl(\rho_0,\u_0\bigr) & x < x(t) - \e, \\
        \bigl(\rho_{0,\e}(t),u_{0,\e}(t)\bigr) & x(t) - \e < x< x(t), \\
        \bigl(\rho_{1,\e}(t),\u_{1,\e}(t)\bigr) & x(t) < x < x(t) + \e, \\
        \bigl(\rho_1,\u_1\bigr) & x(t) + \e < x,
    \end{cases} \label{eq:shadow_wave_sol}
\end{align}
for $\e > 0, x(t) \in C^1\bigl([0,\infty)\bigr)$. 

We expect an unbounded solution in at least one of the two dependent variables, connecting the left state $(\rho_0, u_0)$ to the right state $(\rho,\u)$, while allowing the amplitude to vary across the discontinuity and evolve in time. Thus, we have the piecewise function
\begin{align*}
    \rho_\e(x,t) =& \left\{1 - H(x-x(t)+\e)\right\}\rho_0 \\
    & + \left\{H(x-x(t)+\e) - H(x-x(t))\right\}\rho_{0,\e}(t) \\
    & + \left\{H(x-x(t)) - H(x-x(t)-\e)\right\}\rho_{1,\e}(t) \\
    & + H(x-x(t)-\e)\rho_1,
\end{align*}
and the equation for $\u_\e$ is similar. Let
\begin{align}
    \rho_{0,\e}(t) = \frac{\zeta_0(t)}{\e^k}, \space \rho_{1,\e}(t) = \frac{\zeta_1(t)}{\e^\beta}, \space \u_{0,\e}(t) = \frac{\eta_0(t)}{\e^\gamma}, \space \u_{1,\e}(t) = \frac{\eta_1(t)}{\e^\delta}, \label{eq:shadow_approximations}
\end{align}
where $k, \beta, \gamma, \delta$ will be chosen later such that a solution is admitted. $H$ is the Heaviside step function. Thus, we have
\begin{align*}
    \partial_t\rho_\e &= \delta(x-x(t)+\e)x'(t)\rho_0 \\
    &\quad - \delta(x-x(t)+\e)x'(t)\rho_{0,\e}(t) + \delta(x-x(t))x'(t)\rho_{0,\e}(t) \\
    &\quad + \left\{H(x-x(t)+\e)-H(x-x(t))\right\}\rho_{0,\e}'(t) \\
    &\quad - \delta(x-x(t))x'(t)\rho_{1,\e}(t) + \delta(x-x(t)-\e)x'(t)\rho_{1,\e}(t) \\
    &\quad + \left\{H(x-x(t))-H(x-x(t)-\e)\right\}\rho_{1,\e}'(t) \\
    &\quad - \delta(x-x(t)-\e)x'(t)\rho_1.
\end{align*}

For a test function $\phi \in C_c^\infty(\mathbb{R}\times\R_+)$,
\begin{align*}
    \brackets{\partial_t\rho_\e,\phi} &= \iint_{\mathbb{R}^2_+}\partial_t\rho_\e (x,t)\cdot\phi(x,t)\d x\d t \\
    & = \int_0^\infty \{ \rho_0 x'(t) \phi(x(t)-\e,t) - \rho_{0,\e}(t) x'(t) \phi(x(t)-\e,t) \\
    &\quad + \rho_{0,\e}(t) x'(t) \phi(x(t),t) - \rho_{1,\e}(t) x'(t) \phi(x(t),t) \\
    &\quad + \rho_{1,\e}(t) x'(t) \phi(x(t)+\e) - \rho_1 x'(t) \phi(x(t)+\e,t) \} \d t \\
    &\quad + \int_0^\infty \int_{x(t)-\e}^{x(t)} \rho_{0,\e}'(t) \phi(x,t)\d x\d t \\
    &\quad + \int_0^\infty \int_{x(t)}^{x(t)+\e} \rho_{1,\e}'(t) \phi(x,t)\d x\d t.
\end{align*}

Using the approximation $\phi(x\pm\e,t) = \phi(x,t) \pm \e\partial_x\phi(x,t) + \mathcal{O}(\e^2)$ and the Mean Value Theorem  for $\int_{x(t)}^{x(t)+\e} \phi(x,t)\d x$, we obtain
\begin{align*}
    \brackets{\partial_t\rho_\e,\phi} &\approx \int_0^\infty \left\{x'(t)(\rho_0 - \rho_1) + \e(\rho_{0,\e}'(t) + \rho_{1,\e}'(t))\right\}\phi(x(t),t)\d t \\
    &\quad + \int_0^\infty \left\{\e x'(t)(\rho_{0,\e}(t)-\rho_0 + \rho_{1,\e}(t)-\rho_1)\right\}\partial_x\phi(x(t),t) \d t \\
    & = \brackets{\left(x'(t)(\rho_0-\rho_1) + \e(\rho_{0,
    \e}' + \rho_{1,\e}')\right)\delta(x-x(t)),\phi(x,t)} \\
    &\quad + \brackets{-\e x'(t)(\rho_{0,\e}-\rho_0 + \rho_{1,\e}-\rho_1)\delta'(x-x(t)),\phi(x,t)}.
\end{align*}
Note that $f \approx g$ means that $\frac{f-g}{\e}$ converges to zero as $\e \rightarrow 0$. We continue with the flux:
\begin{align*}
    &\brackets{\partial_x\left(\rho_\e\u_\e\left(1-\fraction{}^a\right)\right),\phi(x,t)} \\
    & = - \brackets{\rho_\e \u_\e \left(1-\fraction{}^a\right),\partial_x\phi(x,t)} \\
    & = -\int_0^\infty \int_{-\infty}^{x(t)-\e} \rho_0 \u_0 \left(1-\fraction{0}^a\right)\partial_x\phi \d x\d t \\
    &\quad - \int_0^\infty \int_{x(t)-\e}^{x(t)} \rho_{0,\e} \u_{0,\e} \left(1-\fraction{0,\e}^a\right) \partial_x\phi \d x\d t \\
    &\quad - \int_0^\infty \int_{x(t)}^{x(t)+\e} \rho_{1,\e} \u_{1,\e} \left(1-\fraction{1,\e}^a\right) \partial_x\phi \d x\d t \\
    &\quad - \int_0^\infty \int_{x(t)+\e}^{\infty} \rho_1 \u_1 \left(1-\fraction{1}^a\right) \partial_x\phi \d x\d t \\
    &\approx \int_0^\infty -f_1(\rho_0,\u_0)\bigl(\phi(x(t),t) - \e\partial_x\phi(x(t),t)\bigr) + f_1(\rho_{0,\e},\u_{0,\e})\bigl(\phi(x(t),t) - \e\partial_x\phi(x(t),t)\bigr) \\
    &\qquad\quad\,\, - f_1(\rho_{0,\e},\u_{0,\e})\phi(x(t),t) + f_1(\rho_{1,\e},\u_{1,\e})\phi(x(t),t) \\
    &\qquad\quad\,\, - f_1(\rho_{1,\e},\u_{1,\e})\bigl(\phi(x(t),t) + \e\partial_x\phi(x(t),t)\bigr) + f_1(\rho_1,\u_1)\bigl(\phi(x(t),t) + \e\partial_x\phi(x(t),t)\bigr)\d t \\
    & = \brackets{-\left\{f_1(\rho_{0},\u_0) - f_1(\rho_{1},\u_1)\right\}\delta(x-x(t)),\phi(x,t)} \\
    &\quad + \brackets{\e\left\{f_1(\rho_{0,\e},\u_{0,\e}) - f_1(\rho_0,\u_0) + f_1(\rho_{1,\e},\u_{1,\e}) - f_1(\rho_1,\u_1)\right\}\delta'(x-x(t)),\phi(x,t) }.
\end{align*}
Since the sum of these two terms must equal $0$ for all $\e > 0$ according to \eqref{eq:sys_1}, we require
\begin{gather}
    x'(t)[\rho] + \lim_{\e\rightarrow0}\left\{\e\rho_{0,\e}'(t)+\e\rho_{1,\e}'(t)\right\} -\left[f_1(\rho,\u)\right] = 0, \label{eq:shadow_A} \\
    - x'(t)\lim_{\e\rightarrow0}\left\{\e\rho_{0,\e}(t)+\e\rho_{1,\e}(t)\right\} + \lim_{\e\rightarrow0}\left\{\e f_1\bigl(\rho_{0,\e}(t),u_{0,\e}(t)\bigr) + \e f_1\bigl(\rho_{1,\e}(t),\u_{1,\e}(t)\bigr)\right\} = 0. \label{eq:shadow_B}
\end{gather}

Using similar calculations for $\brackets{\partial_t(\rho \u),\phi}$ and $\brackets{\partial_x(f_2(\rho,\u)),\phi}$, we obtain
\begin{gather}
    x'(t)[\rho \u] + \lim_{\e\rightarrow0}\left\{\e(\rho_{0,\e} \u_{0,\e})'(t)+\e(\rho_{1,\e} \u_{1,\e})'(t)\right\} - \left[f_2(\rho,\u)\right] = 0,\label{eq:shadow_C} \\
    - x'(t)\lim_{\e\rightarrow0}\left\{\e(\rho_{0,\e} \u_{0,\e})(t)+\e(\rho_{1,\e} \u_{1,\e})(t)\right\} + \lim_{\e\rightarrow0}\left\{\e f_2\bigl(\rho_{0,\e}(t),u_{0,\e}(t)\bigr) + \e f_2\bigl(\rho_{1,\e}(t),\u_{1,\e}(t)\bigr)\right\} = 0. \label{eq:shadow_D}
\end{gather}
Thus, \eqref{eq:sys_1} and \eqref{eq:sys_2} are satisfied by the shadow wave ansatz \eqref{eq:shadow_wave_sol} if and only if \eqref{eq:shadow_A} - \eqref{eq:shadow_D} hold. We now seek to determine values of the parameters $k, \beta, \gamma, \delta$ such that each of the limits is finite, thereby ensuring that the equations are satisfied.

\subsubsection{Delta Solutions for \texorpdfstring{$a < 0$}{a < 0}} \label{marko_anegative}
For $a < 0$, there are many combinations of exponents which yield a solution to \eqref{eq:shadow_A} - \eqref{eq:shadow_D}, but the one we observe numerically (see \cite{REU2025_1}) is described here. It is also the same as that produced by the previous method for delta shocks examined in Section \ref{delta_ansatz}.

Let $k = \beta = 1$ and $\gamma = \delta = 0$. Using \eqref{eq:shadow_approximations}, $\kappa_1(t) = x'(t)[\rho] - [f_1(\rho,\u)]$, and $\kappa_2(t) = x'(t)[\rho \u] - [f_2(\rho,\u)]$, \eqref{eq:shadow_A} - \eqref{eq:shadow_D} become
\begin{align*}
    -\kappa_1(t) =& \, \zeta_0'(t) + \zeta_1'(t), \\
    x'(t)\bigl(\zeta_0(t) + \zeta_1(t)\bigr) =& \, \zeta_0 \eta_0(t) + \zeta_1(t) \eta_1(t), \\
    -\kappa_2(t) =& \, (\zeta_0 \eta_0)'(t) + (\zeta_1 \eta_1)'(t), \\
    x'(t)\bigl((\zeta_0\eta_0)(t) + (\zeta_1\eta_1)(t)\bigr) =& \, (\zeta_0\eta_0)(t)\eta_0(t) + (\zeta_1\eta_1)(t) \eta_1(t).
\end{align*}
Let $\zeta(t) := \zeta_0(t) + \zeta_1(t)$ and $\eta(t) = \eta_0(t) = \eta_1(t)$, so we obtain
\begin{align}
    x'(t) =& \, \eta(t), \label{eq:ODE_1.1} \\
    \zeta'(t) =& -\kappa_1(t), \label{eq:ODE_2.1} \\
    (\zeta \eta)'(t) =& -\kappa_2(t). \label{eq:ODE_3.1}
\end{align}
This system corresponds exactly to the one obtained using the method proposed in Section \ref{delta_ansatz}.

Note that the initial value problem specified in \eqref{eq:initial_conditions} also provides the initial conditions for the present system. Because the unbounded behavior is anticipated at the discontinuity, we set $x(0) = 0$. Additionally, the initial conditions are bounded, which requires $\zeta(0) = 0$. Finally, observe that $x'(t) = \eta(t) \equiv \frac{\kappa_2(t)}{\kappa_1(t)} \equiv \frac{\kappa_2(0)}{\kappa_1(0)}$ solves the system of ODEs. This corresponds to the solutions we numerically observed in \cite{REU2025_1} where $\zeta(t)$ grows linearly over time. Assuming a continuous $\eta(t)$, let us determine the form of $s_- = x'(0)$: let $K_1(t) = \int_0^t \kappa_1(s)\d s$ and $K_2(t) = \int_0^t \kappa_2(s)\d s$.
Then
\begin{align}
    \zeta(t) &= -K_1(t) = -x(t)[\rho] + [f_1(\rho,\u)]t, \notag \\
    \zeta(t) \eta(t) &= -K_2(t) = -x(t)[\rho \u] + [f_2(\rho,\u)]t \notag, \\
    \implies x'(0) = s_- &= \lim_{\e \rightarrow 0}\frac{\zeta(\e) \eta(\e)}{\zeta(\e)} = \lim_{\e \rightarrow 0} \frac{-K_2(\e)}{-K_1(\e)} \notag \\
    &= \lim_{\e \rightarrow 0} \frac{x'(\e)[\rho \u] - [f_2(\rho,\u)]}{x'(\e)[\rho] - [f_1(\rho,\u)]} \notag \\
    &= \frac{s_- [\rho \u] - \left[\rho \u^2 \left(1 - \fraction{}^a\right)\right]}{s_- [\rho] - \left[\rho \u \left(1 - \fraction{}^a\right) \right]} \notag, \\
    \implies s_-^2 [\rho] - s_- &\left([\rho \u] + \left[\rho \u \left(1 - \fraction{}^a\right) \right] \right) + \left[\rho \u^2 \left(1 - \fraction{}^a\right)\right] = 0. \label{eq:eta}
\end{align}

\subsubsection{Delta Solutions for \texorpdfstring{$a > 0$}{a > 0}} \label{marko_apositive}
For $a > 0$, we can also identify multiple combinations that could yield admissible solutions. However, the configuration observed in our numerical simulations satisfies the Rankine–Hugoniot condition only in the second component of \eqref{eq:RH_conditions} while a deficit remains in the first component; that is, $\kappa_1(t) \not\equiv 0 \equiv \kappa_2(t)$. Furthermore, numerical evidence indicates that the delta concentration in $\rho$ occurs when $u = 0$, suggesting that $\lim_{\e\rightarrow0} \u_{i,\e} = 0$.

Let $k = \beta = 1$, $\gamma = \delta = -a$. Then \eqref{eq:shadow_A} - \eqref{eq:shadow_D} become
\begin{align*}
    -\kappa_1(t) =& \, \zeta_0'(t) + \zeta_1'(t), \\
    -x'(t)\bigl(\zeta_0(t) + \zeta_1(t)\bigr) =& \, \zeta_0(t) \eta_0(t) \left(\frac{\zeta_0(t)}{\rhobar}\right)^a + \zeta_1(t) \eta_1(t) \left(\frac{\zeta_1(t)}{\rhobar}\right)^a, \\
    \kappa_2(t) =& \, 0 \Leftrightarrow x'(t) = \frac{[f_2(\rho,\u)]}{[\rho \u]}, \\
    0 =& \, 0.
\end{align*}
Note now that $[\rho \u] \neq 0$ outside of the regions in which we can connect the left and right states with combinations of classical curves. Let $\zeta(t) := \zeta_0(t) + \zeta_1(t)$ and $\eta(t) = \eta_0(t) = \eta_1(t)$ so that we obtain
\begin{align}
    x'(t) \equiv& \, s_+ = \frac{[f_2(\rho,\u)]}{[\rho \u]}, \label{eq:ODE_4} \\
    \zeta'(t) =& \, -\kappa_1 = -s_+ [\rho] + [f_1(\rho,\u)], \label{eq:ODE_5} \\
    \eta(t) =& -s_+ \left(\frac{\rhobar}{\frac{1}{2}\zeta(t)}\right)^a. \label{eq:ODE_6}
\end{align}
As we are still working within the framework of the Riemann problem, we impose the same type of initial conditions as in the previous case, namely, $\zeta(0) = 0$ and $x(0) = 0$.

Our solution to \eqref{eq:shadow_wave_sol}, therefore, in the limit $\e \rightarrow 0$ is of the form
$$\u(x,t) = \begin{cases} \u_{L} & x < x(t), \\ 0 & x = x(t), \\ \u_{R} & x > x(t), \end{cases}
\qquad
\rho(x,t) = \begin{cases} \rho_{L} & x < x(t), \\ \left(-\frac{[f_2(\rho,\u)]}{[\rho \u]} [\rho] + [f_1(\rho,\u)]\right)t\cdot\delta\bigl(x - x(t)\bigr) & x = x(t), \\
\rho_{R} & x > x(t).\end{cases}$$

\begin{figure}
    \centering
    \includegraphics[width=0.5\linewidth]{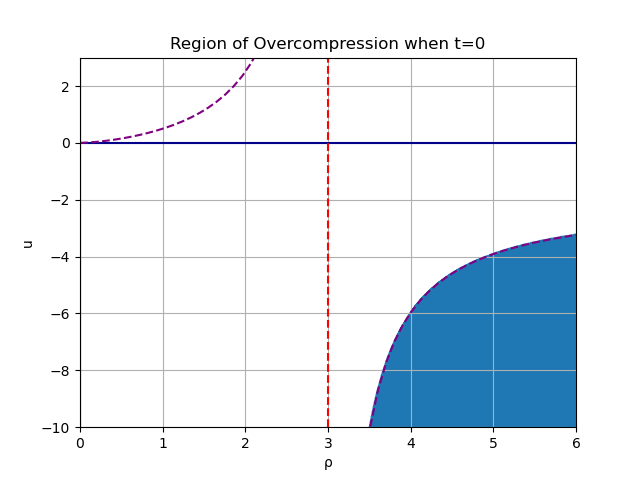}
    \caption{Overcompressive Region for Case 4}
    \label{fig:overcompressive}
\end{figure}

We seek delta-shocks connecting a given left state $(\rho_L, u_L)$ with a right state $(\rho_R,u_R)$ that are overcompressive, meaning that all characteristic curves run into the delta-shock curve from both sides. Therefore, we use the following inequality as our admissibility criterion
\begin{align}
    \max\left\{\lambda_a(U), \lambda_0(U)\right\} < s_\pm < \min\left\{\lambda_a(U_L), \lambda_0(U_L)\right\}. \label{eq:overcompressive}
\end{align}
We can now compare the initial speed $x'(0) = s_\pm$ with our overcompressive admissibility condition \eqref{eq:overcompressive} and determine where delta-shock solutions are expected (see Figure \ref{fig:overcompressive}). We also note that reducing \eqref{eq:ODE_1.2} - \eqref{eq:ODE_3.2} to solving for $\u_\delta(t)$ yields a nonlinear second-order ODE which cannot be solved analytically. Although other solution profiles may exist, our numerical experiments consistently yield the constant-speed solution described above.

In the following section, we present plots for various values of $a$, highlighting the resulting wave structures and their dependence on the system's parameters. 

\section{Regions of Classical and Nonclassical Solutions: 24 Cases} \label{numerics}
\subsection{Numerical Preliminaries}
The Local Lax-Friedrichs (LLF) scheme is a well-known numerical approximation that is used to reconstruct the flux of hyperbolic conservation laws. In essence, the scheme begins with a discretization of the spatial and temporal domains into a mesh grid, after which the solution is numerically approximated at each mesh point in fixed time. By letting $\Delta x$ and $\Delta t$ be the corresponding cell sizes, we can write any point in our grid as $(x_i, t_j)=(i\Delta x, j \Delta t)$, and the solution $H=(\rho \enspace \rho \u)^T$ can be written similarly at any point as $H_i^j=H(x_i,t_j)$. To find the solution $H$ at the point $(x_i, t_{j+1})$, one only needs to have constructed the values at $(x_{i-1}, t_j)$ and $(x_{j+1},t_j)$. The relevant flux reconstruction is given by the equation:
\begin{align}
H_{i}^{j+1} = \frac{1}{2}(H_{i-1}^j + H_{i+1}^j)+\frac{CFL}{2\lambda}(G_{i+1}^{j}-G_{i-1}^{j})
\end{align}
where $CFL=\lambda \frac{\Delta t}{\Delta x}$ is the Courant number and $\lambda := \max_{i} |\lambda_i|$ is the greater  of the two characteristic speeds. The Courant number measures the numerical stability of the scheme and must satisfy the inequality
\begin{align}CFL \leq \frac{1}{2}.\end{align}
By choosing $\Delta x=1$, this can be ensured by simply requiring $\lambda \Delta t \leq \frac{1}{2}$ throughout the procedure. Thus, our spatial grid size is always fixed, but our temporal steps may differ greatly depending on the eigenvalues calculated at each iteration. The CFL condition guarantees that the scheme converges to the physically correct weak solution satisfying the Lax Entropy Condition \cite{Tadmor}. Additionally, for a more in-depth treatment and explanation of the LLF scheme, consult \cite{Lev_1,Lev_2}.

Because the system's eigenvalues can greatly influence the values of $\Delta t$, certain choices of physical parameters prove easier to analyze numerically than others. The problem can be simplified into a multitude of cases depending on the system parameters and initial states. For our simulation, we chose the following conditions:
\begin{itemize}
    \item $\rhobar=5$;
    \item $a\in \{-1.5,-1,-0.5,0.5\}$ for each of the four cases of $a$ we will consider;
    \item $\rho_L \in \{3, 5, 8\}$ and $\u_L = \pm4$, depending on the ordering $\rho_L\text{ and }\rhobar$  and the sign of $\u_L$.
\end{itemize}
Note that we have introduced the change of variables in our code $m = \rho\u$ so that our solution vector may be written more simply as $H = (\rho \enspace m)^T$. Furthermore, our left and right states are represented in our figures as $L=(\rho_L \enspace \u_L)^T$ and $R=(\rho \enspace \u)^T$, respectively. Lastly, to ensure accuracy, data was renormalized every 100 steps within an error bound of $10^{-7}$ to minimize numerical noise.

\subsection{Interpretation of Numerical Outputs}
Because our solution is self-similar when there is no time-dependence in the flux functions (see Culver et al. \cite{REU2025_1} for a treatment of the system with a time-dependent source term), we can express the solution in terms of $x/t$. When interpreting the numerical figures, we look at the graphs of $\rho$ and $\u$ against $x/t$ together and compare how $\rho$ and $u$ change or remain constant at the same $x/t$ values. Our original plots showed 20 iterations (each with 1000 steps) to illustrate behavior over time, while the right column displayed the final iteration, the iteration most closely portraying the weak solution. We used an increasing thickness to differentiate the graphs at different times. A single graph may show multiple features in sequence, such as an $R_a$ followed by a vacuum and then a $C_0$. We expect the following behaviors from our analysis in Sections \ref{classical_analysis} and \ref{delta_shocks}:
\begin{itemize}
    \item If $\rho$ changes (increasing or decreasing) while $u$ remains constant:
    \begin{itemize}
        \item For $a \neq -1$, this indicates an $S_a$, a shock of the $a$-family or $a$-shock, or an $R_a$, an $a$-rarefaction:
        \begin{itemize}
            \item If $\rho$ changes abruptly, there is a jump discontinuity and thus an $a$-shock.
            \item If $\rho$ changes continuously, there is an $R_a$.
        \end{itemize}
        \item For $a = -1$, this indicates a $C_a$, an $a$-contact discontinuity.
    \end{itemize}
    \item If $\rho_L \neq \rhobar$, then $\rho$ and $u$ changing simultaneously indicates a $C_0$. However, if $\rho_L = \rhobar$, then $C_0$ appears as a vertical line, where $\rho = \rho_L = \rhobar$ remains constant while $u$ changes.
    \item If $\rho$ reaches zero, this is called a vacuum state.
    \item If $u$ is changing by a small amount, but $\rho$ spikes to a large value (resembling a Dirac delta), this indicates a delta shock.
\end{itemize}

\subsection{Discussion of Regions}
The equations for our shocks, rarefactions, and contact discontinuities are given by
\begin{gather}
    S(\rho_L,\u_L): \u = \u_L, \\
    R(\rho_L,\u_L):\u=\u_L, \\
    C(\rho_L,\u_L): \u = \u_L \left(\frac{\rho^a - \rho_L^a}{\rhobar^a - \rho^a}  +1 \right).
\end{gather}
The regions are further defined by the $\rho$-axis and the max bound of the overcompressive region calculated numerically in Python following the shadow wave method. Take Case 4 (which will be described in greater depth later on) as an example where $C_0,_m$ bounds the region where delta shock solutions occur, marked by the shaded area. 
We will find that it is often the case that these overcompressive regions are bounded by $C_{0,m}$ and by the $\rho$-axis on occasion. Some limiting behavior of $C(\rho,\u)$ as $\rho \rightarrow 0$ or $\infty$ also occurs and separates regions of solutions, including the asymptotes discussed at the end of Section \ref{contact_disconts} and, when $a < 0$, $C_{0,\ell}$ the limiting contact discontinuity at $\rho_L \rightarrow \infty$, given by the equation
\begin{gather}
    C_{\ell}(\rho_L,\u_L): u = \u_L \left(\frac{- \rho_L^a}{\rhobar^a - \rho^a}  +1 \right).
\end{gather}

Recall the formulation of the Riemann problem in equation \eqref{eq:initial_conditions}. We have a left state $(\rho_L, \u_L)$ and a right state $(\rho, \u)$, and we seek a solution path between them in state space. For any given case, we can split the state space into regions of different solutions and analyze them accordingly. Of our twenty-four cases, we observe a total of six types of regions, each displaying a characteristic wave or combination of waves. Listed here, they are
\begin{itemize}
    \item Region $I$: $S_a$ and $C_0$\newline
    The solution to the Riemann problem in this region is classical. The left state is connected to the right by means of a shock and contact discontinuity, the order of which depends on the case being considered.
    \item Region $II$: $R_a$ and $C_0$\newline
    The transition from the left state to the right state occurs via a rarefaction wave followed by a contact discontinuity. The rarefaction continuously connects states that do not satisfy the Lax Entropy Condition, as the characteristic speed increases across the wave. The order of the waves is again dependent on the eigenvalues in the specific case.
    \item Region $III$: $C_a$ and $C_0$\newline
    When $a = -1$, the shock and rarefactions of the $a$-family become a single contact discontinuity along the line $u = \u_L$. The order of the two contacts discontinuities depends on whether $\ul{L} > 0$ or not.
    \item Region $IV$: Delta Shock $S_\delta$\newline
    The right state exists in a region that cannot be reached by classical means. The solution becomes unbounded in the $\rho$ variable when the left and right states satisfy the overcompressibility condition. There may also be curves of the $a$-family which supplement the delta shock. A lack of subscript will denote regions where only a delta shock connects the left and right states.
    \item Region $V$: Vacuum\newline
    The right state can only be reached by allowing the density of the system to hit zero, creating a vacuum state. In some of the figures that follow, the vacuum state may be divided into subregions based on the waves that combine to produce the solution. This will be denoted by the subscript for the region: for instance, ``Region $V_{C_0VC_aC_0}$'' refers to a contact discontinuity into a vacuum followed by two more contact discontinuities. In this example, two contact discontinuities of the same $0$-family are allowed because we pass through the degeneracy $\rho = 0$.
    \item Region $VI$: States that pass through degeneracy\newline
    When $u = 0$, the system loses hyperbolicity and we have to ``reset'' the Riemann problem. Thus, we are able to employ two classical curves of the same family on either side of the degeneracy, which will arise as a curve of the $a$-family followed by a $C_0$ along the line $\rho = \rhobar$ (when $\lambda_0 = 0$) and finishing with another curve of the $a$-family.
\end{itemize}
Note that in the region descriptions throughout this paper, the subscript indicates which wave family appears first in the path. For example, $I_a$ means that in the region $I$, the left state first connects to the middle state via the $a$-family curve, specifically through a shock $S_a$, followed by a contact discontinuity $C_0$.

\subsubsection{\texorpdfstring{$a<-1$}{a< -1}}
We consider first
\begin{itemize}
    \item Case 1: $u_L<0$ and $\rho_L<\rhobar$,
    \item Case 2: $u_L<0$ and $\rho_L=\rhobar$, and 
    \item Case 3: $u_L<0$ and $\rho_L>\rhobar$,
\end{itemize}
and we discuss the various regions and the types of wave combinations that comprise the solution to the Riemann Problem, detailing the path taken to reach the right state.

For Cases 1 - 3, we have the following regions:
\begin{itemize}
    \item Region $I_a$: The left state is connected to the right by an $a$-shock followed by a $0$-contact discontinuity. 
    \item Region $II_a$: The left state is connected to the right by an $a$-rarefaction curve followed by a $0$-contact discontinuity. 
    \item Region $IV$: This region follows an overcompressive delta shock which transitions it to a right state.
    \item Region $VI$: The left state will either take an $R_a$ or $S_a$, in Cases 1 and 3, respectively, to reach $\rhobar$, unless the left state is already at $\rhobar$ to begin with. Then it will go up $\rho = \rhobar$ using the $C_0$ curve of Case 2. Finally, it will then take either an $S_a$ or an $R_a$ from $(\rhobar,\u)$ to reach the right state $(\rho,\u)$ as in Case 5 (presented shortly).
\end{itemize}

\begin{figure}[H]
    \begin{minipage}[t]{.5\textwidth}
        \vspace{0pt}
        \centering
        \begin{tikzpicture}[scale = 1,
        declare function={
            a = -1.5;       
            rhowithbar = 5; 
            rholeft = 3;    
            uleft = -4;      
            A = 0;          
            cont_disc_rholeftnotrhowithbar(\t)
                = uleft + (uleft - A)*(exp(ln(rholeft)*(a)) - exp(ln(\t)*(a)))/(exp(ln(\t)*(a)) - exp(ln(rhowithbar)*(a)));
            cont_disc_limit
                = uleft + (uleft - A)*(exp(ln(rholeft)*(a)))/(- exp(ln(rhowithbar)*(a)));
            cont_disc_limiting_curve(\t)
                = uleft + (uleft - A)*(exp(ln(1000)*(a)) - exp(ln(\t)*(a)))/(exp(ln(\t)*(a)) - exp(ln(rhowithbar)*(a)));
        }
    ]
    \tikzstyle{arrow} = [thick,->,>=stealth]

    \begin{axis}[
        axis lines = center,
        xmin=0, xmax=15,
        ymin=-10, ymax=10,
        xlabel=$\rho$,
        ylabel=$\u$,
        x label style = {at={(axis description cs:0.95,0.4)},anchor=south},
        y label style = {at={(axis description cs:0.0,0.95)},anchor=west},
        xticklabel=\empty,
        yticklabel=\empty,
        major tick style={draw=none}
    ]

    \node at (axis cs:rholeft,uleft) {\textbullet};

    \draw[dashed] (axis cs:rhowithbar,10) -- (axis cs:rhowithbar,-10)
        node[pos=0.97, right] {$\overline{\rho}$};

    \addplot [
        domain=0.01:(rhowithbar-1.12),
        samples=100,
        color=black,
        solid
    ] {cont_disc_rholeftnotrhowithbar(x)}
        node[pos=0.9, left] {\footnotesize $C_{0}$};
    \addplot [
        domain=(rhowithbar+2.5):15,
        samples=100,
        color=black,
        solid
    ] {cont_disc_rholeftnotrhowithbar(x)}
        node[pos=0.4, left] {\footnotesize $C_{0,m}$};
    \addplot [
        domain=(rhowithbar+2):15,
        samples=100,
        color=black,
        solid
    ] {cont_disc_limiting_curve(x)}
        node[pos=0.2, left] {\footnotesize $C_{0,\ell}$};
    \draw[dotted] (axis cs:0,cont_disc_limit) -- (axis cs:15,cont_disc_limit);

    \draw[solid] (axis cs:rholeft,uleft) -- (axis cs:15,uleft)
        node[pos=0.95, above] {\footnotesize $R_a$};
    \draw[solid] (axis cs:rholeft,uleft) -- (axis cs:0,uleft)
        node[pos=0.8, above] {\footnotesize $S_a$};

    \node at (axis cs:10,4.5) {\footnotesize $VI$};
    \draw[black,-latex,very thin] (axis cs:10.5,5) -- (axis cs:12,7.5);
    \node at (axis cs:2.5,4.5) {\footnotesize $VI$};
    \node at (axis cs:1.5,-5.5) {\footnotesize $I_a$};
    \node at (axis cs:12,-8) {\footnotesize $IV$};
    \node at (axis cs:4,-2) {\footnotesize $II_a$};
    \draw[black,-latex,very thin] (axis cs:4,-2.75) -- (axis cs:4.25, -6);
    \draw[black,-latex,very thin] (axis cs:4.5,-2) -- (axis cs:10,-2);
    \draw[black,-latex,very thin] (axis cs:4.5,-2.4) -- (axis cs:7.5,-6);
    \end{axis}
    
\end{tikzpicture}
    \end{minipage}\hfill
    \begin{minipage}[t]{.5\textwidth}
        \vspace{0pt}
        \centering
        \begin{tikzpicture}[scale = 1,
        declare function={
            a = -1.5;       
            rhowithbar = 5; 
            rholeft = 5;    
            uleft = -4;      
            A = 0;          
            cont_disc_rholeftnotrhowithbar(\t)
                = uleft + (uleft - A)*(exp(ln(rholeft)*(a)) - exp(ln(\t)*(a)))/(exp(ln(\t)*(a)) - exp(ln(rhowithbar)*(a)));
            cont_disc_limit
                = uleft + (uleft - A)*(exp(ln(rholeft)*(a)))/(- exp(ln(rhowithbar)*(a)));
            cont_disc_limiting_curve(\t)
                = uleft + (uleft - A)*(exp(ln(1000)*(a)) - exp(ln(\t)*(a)))/(exp(ln(\t)*(a)) - exp(ln(rhowithbar)*(a)));
        }
    ]

    \tikzstyle{arrow} = [thick,->,>=stealth]

    \begin{axis}[
        axis lines = center,
        xmin=0, xmax=15,
        ymin=-10, ymax=10,
        xlabel=$\rho$,
        ylabel=$\u$,
        x label style = {at={(axis description cs:0.95,0.4)},anchor=south},
        y label style = {at={(axis description cs:0.0,0.95)},anchor=west},
        xticklabel=\empty,
        yticklabel=\empty,
        major tick style={draw=none}
    ]

    \node at (axis cs:rholeft,uleft) {\textbullet};

    \draw[dashed] (axis cs:rhowithbar,10) -- (axis cs:rhowithbar,-10);

    \addplot [
        domain=(rhowithbar+2):15,
        samples=100,
        color=black,
        solid
    ] {cont_disc_limiting_curve(x)}
        node[pos=0.2, left] {\footnotesize $C_{0,\ell}$};
    \draw[solid] (axis cs:rhowithbar,10) -- (axis cs:rhowithbar,-10)
        node[pos=0.97, right] {\footnotesize $C_0$};;

    \draw[solid] (axis cs:rholeft,uleft) -- (axis cs:15,uleft)
        node[pos=0.95, above] {\footnotesize $R_a$};
    \draw[solid] (axis cs:rholeft,uleft) -- (axis cs:0,uleft)
        node[pos=0.85, above] {\footnotesize $S_a$};

    \node at (axis cs:2.5,-7.5) {\footnotesize $I_a$};
        \draw[black,-latex,very thin] (axis cs:2.5,-6.8) -- (axis cs:2.5,-2);
        
    \node at (axis cs:11.5,-8) {\footnotesize $IV$};

    \node at (axis cs:2.5,5) {\footnotesize $VI$};
    
    \node at (axis cs:10,5) {\footnotesize $VI$};

    \node at (axis cs:10, -2) {\footnotesize $II_a$};
    
    \draw[black,-latex,very thin] (axis cs:9.5,-2.5) -- (axis cs:7,-6);
  
    \end{axis}
\end{tikzpicture}
    \end{minipage}
    \!\!\!\!\!\!\!\!\begin{minipage}[b]{.5\textwidth}
        \vspace{0pt}
        \centering
        \begin{tikzpicture}[scale = 1,
        declare function={
            a = -1.5;       
            rhowithbar = 5; 
            rholeft = 8;    
            uleft = -4;      
            A = 0;          
            cont_disc_rholeftnotrhowithbar(\t)
                = uleft + (uleft - A)*(exp(ln(rholeft)*(a)) - exp(ln(\t)*(a)))/(exp(ln(\t)*(a)) - exp(ln(rhowithbar)*(a)));
            cont_disc_limit
                = uleft + (uleft - A)*(exp(ln(rholeft)*(a)))/(- exp(ln(rhowithbar)*(a)));
            cont_disc_limiting_curve(\t)
                = uleft + (uleft - A)*(exp(ln(1000)*(a)) - exp(ln(\t)*(a)))/(exp(ln(\t)*(a)) - exp(ln(rhowithbar)*(a)));
        }
    ]

    \tikzstyle{arrow} = [thick,->,>=stealth]

    \begin{axis}[
        axis lines = center,
        xmin=0, xmax=15,
        ymin=-10, ymax=10,
        xlabel=$\rho$,
        ylabel=$\u$,
        x label style = {at={(axis description cs:0.95,0.4)},anchor=south},
        y label style = {at={(axis description cs:0.0,0.95)},anchor=west},
        xticklabel=\empty,
        yticklabel=\empty,
        major tick style={draw=none}
    ]

    \node at (axis cs:rholeft,uleft) {\textbullet};

    \draw[dashed] (axis cs:rhowithbar,10) -- (axis cs:rhowithbar,-10)
        node[pos=0.97, right] {$\overline{\rho}$};

    \addplot [
        domain=0.01:(rhowithbar-0.58),
        samples=100,
        color=black,
        solid
    ] {cont_disc_rholeftnotrhowithbar(x)}
        node[pos=0.8, left] {\footnotesize $C_{0,m}$};
    \addplot [
        domain=(rhowithbar+0.8):15,
        samples=100,
        color=black,
        solid
    ] {cont_disc_rholeftnotrhowithbar(x)}
        node[pos=0.3, left] {\footnotesize $C_{0}$};
    \addplot [
        domain=(rhowithbar+2):15,
        samples=100,
        color=black,
        solid
    ] {cont_disc_limiting_curve(x)}
        node[pos=0.2, left] {\footnotesize $C_{0,\ell}$};
    \draw[dotted] (axis cs:0,cont_disc_limit) -- (axis cs:15,cont_disc_limit);

    \draw[solid] (axis cs:rholeft,uleft) -- (axis cs:15,uleft)
        node[pos=0.95, above] {\footnotesize $R_a$};
    \draw[solid] (axis cs:rholeft,uleft) -- (axis cs:0,uleft)
        node[pos=0.85, above] {\footnotesize $S_a$};

    \node at (axis cs:11.5,-8.5) {\footnotesize $IV$};
    
    \node at (axis cs:2,5) {\footnotesize $VI$};
        \draw[black,-latex,very thin] (axis cs:2.1,4.25) -- (axis cs:4,1.25);

    \node at (axis cs:10,5) {\footnotesize $VI$};
    
    \node at (axis cs:2.5,-7) {\footnotesize $I_a$};
        \draw[black,-latex,very thin] (axis cs:2.5,-6.4) -- (axis cs:2.5,-1.5);
        \draw[black,-latex,very thin] (axis cs:2.75,-6.5) -- (axis cs:8,-1.5);
        \draw[black,-latex,very thin] (axis cs:2.85,-7) -- (axis cs:5.75,-7);

    \node at (axis cs:10,-5) {\footnotesize $II_a$};

    \end{axis}
\end{tikzpicture}
    \end{minipage}%
    \hfill
    \begin{minipage}[b]{.45\textwidth}
        \caption{\\
            Case 1 (above left): $\lambda_a < 0 < \lambda_0$ \\
            Case 2 (above right): $\lambda_a < \lambda_0 = 0$ \\
            Case 3 (left): $\lambda_a < \lambda_0 < 0$
        }
        \label{fig:cases1-3}
        \vspace{2cm}
    \end{minipage}%
\end{figure}
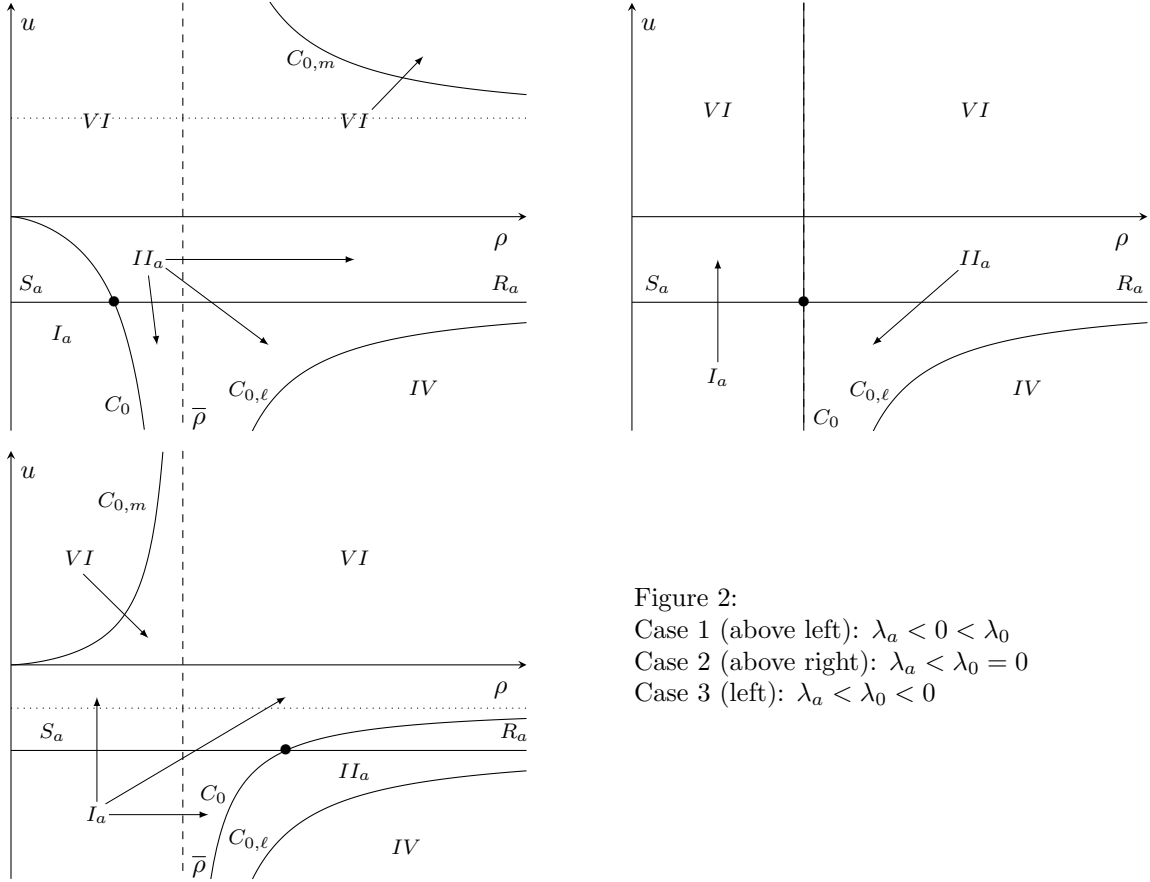

\quad As shown above, the solution to the Riemann problem will involve multiple different cases (as we use the word) for the various middle states as the values of $\rho$ and $\u$ evolve. For instance, suppose we have chosen our initial conditions so that we are in Case 1 with $\rho > \rhobar$ and $\u < 0$. We then have that the solution will follow the rarefaction curve $R_a$ into Case 3 before following the contact discontinuity $C_0$. We have carefully chosen the numbering of our regions so that a right state, say, in Region $I_a$ is always achieved in the manner of $S_a \rightarrow C_0$. In Figure \ref{fig:cases1-3}, we present the state spaces for Cases 1 - 3:

We then turn to the remaining three regions of $(\rho_L,\u_L)$:
\begin{itemize}
    \item Case 4: $\u_L>0$ and $\rho_L<\rhobar$, 
    \item Case 5: $\u_L>0$ and $\rho_L=\rhobar$, and
    \item Case 6: $\u_L>0$ and $\rho_L>\rhobar$.
\end{itemize}

These are comparable to Cases 1, 2, and 3 except that $\u_L>0$. In addition, a vacuum region exists in these cases where $\rho$ hits zero. The regions for the solution to the Riemann problem are as follows:
\begin{itemize}
    \item Regions $I_0$, $II_0$, and $IV$: The same wave connections as in Cases 1 - 3, but the sequence in which these waves occur is reversed. That is, in Region $I_0$, the left state is connected to the right in the manner of $C_0 + S_a$, instead of $S_a + C_0$.
        \begin{itemize}
            \item Region $IV_{\delta}$: The left state first follows an over-compressive $\delta$-shock $S_{\delta}$. It then follows an $a$-rarefaction to the right state.
        \end{itemize}
    \item Region $V$: This occurs when the path from the left state to the right state passes through a vacuum state. In cases where a vacuum arises, the path between states often involves a sequence of waves. The following are typical wave sequences that occur when a vacuum is present:
        \begin{itemize}
            \item $V_{C_{0}VC_{0}}$: The solution consists of a $C_0$ from the left state to a vacuum, followed by another $C_0$ from the vacuum to the right state as in Case 1.
            \item $V_{C_{0}VR_{a}C_{0}}$: The solution consists of a $C_0$ from the left state to a vacuum. It then goes from the vacuum to an $R_a$, followed by another $C_0$ as in Case 1, going to Case 3.
            \item $V_{R_{a}VC_{0}}$ - The solution consists of an $R_a$ from the left state to a vacuum. It then goes from the vacuum to a $C_0$ as in Case 1.
        \end{itemize}
\end{itemize}

Below, we present the state spaces for Cases 4 through 6:
\begin{figure}[H]
    \begin{minipage}[t]{.5\textwidth}
        \vspace{0pt}
        \centering
        \begin{tikzpicture}[scale = 1,
        declare function={
            a = -1.5;       
            rhowithbar = 5; 
            rholeft = 3;    
            uleft = 4;      
            A = 0;          
            cont_disc_rholeftnotrhowithbar(\t)
                = uleft + (uleft - A)*(exp(ln(rholeft)*(a)) - exp(ln(\t)*(a)))/(exp(ln(\t)*(a)) - exp(ln(rhowithbar)*(a)));
            cont_disc_limit
                = uleft + (uleft - A)*(exp(ln(rholeft)*(a)))/(- exp(ln(rhowithbar)*(a)));
            cont_disc_limiting_curve(\t)
                = uleft + (uleft - A)*(exp(ln(1000)*(a)) - exp(ln(\t)*(a)))/(exp(ln(\t)*(a)) - exp(ln(rhowithbar)*(a)));
        }
    ]

    \tikzstyle{arrow} = [thick,->,>=stealth]

    \begin{axis}[
        axis lines = center,
        xmin=0, xmax=15,
        ymin=-10, ymax=10,
        xlabel=$\rho$,
        ylabel=$\u$,
        x label style = {at={(axis description cs:0.95,0.4)},anchor=south},
        y label style = {at={(axis description cs:0.0,0.95)},anchor=west},
        xticklabel=\empty,
        yticklabel=\empty,
        major tick style={draw=none}
    ]

    \node at (axis cs:rholeft,uleft) {\textbullet};

    \draw[dashed] (axis cs:rhowithbar,10) -- (axis cs:rhowithbar,-10)
        node[pos=0.97, right] {$\overline{\rho}$};;

    \addplot [
        domain=0.01:(rhowithbar-1.12),
        samples=100,
        color=black,
        solid
    ] {cont_disc_rholeftnotrhowithbar(x)}
        node[pos=0.9, left] {\footnotesize $C_{0}$};
    \addplot [
        domain=(rhowithbar+2.5):15,
        samples=100,
        color=black,
        solid
    ] {cont_disc_rholeftnotrhowithbar(x)}
        node[pos=0.4, above left] {\footnotesize $C_{0,m}$};
    \draw[dotted] (axis cs:0,cont_disc_limit) -- (axis cs:15,cont_disc_limit);

    \draw[solid] (axis cs:rholeft,uleft) -- (axis cs:15,uleft)
        node[pos=0.95, above] {\footnotesize $S_a$};
    \draw[solid] (axis cs:rholeft,uleft) -- (axis cs:0,uleft)
        node[pos=0.7, above] {\footnotesize $R_a$};

    \node at (axis cs:1.75,6.5) {\footnotesize $II_0$};
    
    \node at (axis cs:10,6.5) {\footnotesize $I_0$};
        \draw[black,-latex,very thin] (axis cs:9.6,6.2) -- (axis cs:3.75,2.8);
        \draw[black,-latex,very thin] (axis cs:9.85,5.9) -- (axis cs:9.85,2);
        \draw[black,-latex,very thin] (axis cs:9.55,6.5) -- (axis cs:4.2, 6.5);

    \node at (axis cs:12,-8.5) {\footnotesize $IV$};
    
    \node at (axis cs:2.5, -2.5) {\footnotesize $V_{C_{0}VC_{0}}$};
        \draw[black,-latex,very thin] (axis cs:2.5,-3) -- (axis cs:2.5,-7.5);

    \node at (axis cs:10, -2.5) {\footnotesize $V_{C_{0}VR_{a}C_{0}}$};
        \draw[black,-latex,very thin] (axis cs:10,-3) -- (axis cs:6,-6);
    \end{axis}
\end{tikzpicture}
    \end{minipage}
    \hfill
    \begin{minipage}[t]{.5\textwidth}
        \vspace{0pt}
        \centering
        \begin{tikzpicture}[
        declare function={
            a = -1.5;       
            rhowithbar = 5; 
            rholeft = 5;    
            uleft = 4;      
            A = 0;          
            cont_disc_rholeftnotrhowithbar(\t)
                = uleft + (uleft - A)*(exp(ln(rholeft)*(a)) - exp(ln(\t)*(a)))/(exp(ln(\t)*(a)) - exp(ln(rhowithbar)*(a)));
            cont_disc_limit
                = uleft + (uleft - A)*(exp(ln(rholeft)*(a)))/(- exp(ln(rhowithbar)*(a)));
            cont_disc_limiting_curve(\t)
                = uleft + (uleft - A)*(exp(ln(1000)*(a)) - exp(ln(\t)*(a)))/(exp(ln(\t)*(a)) - exp(ln(rhowithbar)*(a)));
        }
    ]

    \tikzstyle{arrow} = [thick,->,>=stealth]

    \begin{axis}[
        axis lines = center,
        xmin=0, xmax=15,
        ymin=-10, ymax=10,
        xlabel=$\rho$,
        ylabel=$\u$,
        x label style = {at={(axis description cs:0.95,0.4)},anchor=south},
        y label style = {at={(axis description cs:0.0,0.95)},anchor=west},
        xticklabel=\empty,
        yticklabel=\empty,
        major tick style={draw=none}
    ]

    \node at (axis cs:rholeft,uleft) {\textbullet};

    \draw[dashed] (axis cs:rhowithbar,10) -- (axis cs:rhowithbar,-10);

    \draw[solid] (axis cs:rhowithbar,10) -- (axis cs:rhowithbar,-10)
        node[pos=0.97, right] {\footnotesize $C_0$};

    \draw[solid] (axis cs:rholeft,uleft) -- (axis cs:15,uleft)
        node[pos=0.95, above] {\footnotesize $S_a$};
    \draw[solid] (axis cs:rholeft,uleft) -- (axis cs:0,uleft)
        node[pos=0.85, above] {\footnotesize $R_a$};

    \node at (axis cs:2.5,7) {\footnotesize $II_0$};
        \draw[black,-latex,very thin] (axis cs:2.5,6.5) -- (axis cs:2.5,2);
    
    \node at (axis cs:10,7) {\footnotesize $I_0$};
        \draw[black,-latex,very thin] (axis cs:9.9,6.5) -- (axis cs:9.9,2);
    
    \node at (axis cs:2.5,-5) {\footnotesize $V_{R_{a}VC_{0}}$};
    
    \node at (axis cs:10,-5) {\footnotesize $IV$};
  \end{axis}
\end{tikzpicture}
    \end{minipage}
    \!\!\!\!\!\!\!\!\begin{minipage}[b]{.5\textwidth}
        \vspace{0pt}
        \centering
        \begin{tikzpicture}[
        declare function={
            a = -1.5;       
            rhowithbar = 5; 
            rholeft = 8;    
            uleft = 4;      
            A = 0;          
            cont_disc_rholeftnotrhowithbar(\t)
                = uleft + (uleft - A)*(exp(ln(rholeft)*(a)) - exp(ln(\t)*(a)))/(exp(ln(\t)*(a)) - exp(ln(rhowithbar)*(a)));
            cont_disc_limit
                = uleft + (uleft - A)*(exp(ln(rholeft)*(a)))/(- exp(ln(rhowithbar)*(a)));
            cont_disc_limiting_curve(\t)
                = uleft + (uleft - A)*(exp(ln(1000)*(a)) - exp(ln(\t)*(a)))/(exp(ln(\t)*(a)) - exp(ln(rhowithbar)*(a)));
        }
    ]
    \tikzstyle{arrow} = [thick,->,>=stealth]
    \begin{axis}[
        axis lines = center,
        xmin=0, xmax=15,
        ymin=-10, ymax=10,
        xlabel=$\rho$,
        ylabel=$\u$,
        x label style = {at={(axis description cs:0.95,0.4)},anchor=south},
        y label style = {at={(axis description cs:0.0,0.95)},anchor=west},
        xticklabel=\empty,
        yticklabel=\empty,
        major tick style={draw=none},
        view={0}{90}
    ]

    \node at (axis cs:rholeft,uleft) {\textbullet};

    \draw[dashed] (axis cs:rhowithbar,10) -- (axis cs:rhowithbar,-10)
        node[pos=0.97, right] {$\overline{\rho}$};

    \addplot [
        domain=0.01:(rhowithbar-0.58),
        samples=100,
        color=black,
        solid
    ] {cont_disc_rholeftnotrhowithbar(x)}
        node[pos=0.9, left] {\footnotesize $C_{0,m}$};
    \addplot [
        domain=(rhowithbar+0.8):15,
        samples=100,
        color=black,
        solid
    ] {cont_disc_rholeftnotrhowithbar(x)}
        node[pos=0.05, right] {\footnotesize $C_{0}$};
    \draw[dotted] (axis cs:0,cont_disc_limit) -- (axis cs:15,cont_disc_limit);

    \draw[solid] (axis cs:rholeft,uleft) -- (axis cs:15,uleft)
        node[pos=0.9, above] {\footnotesize $S_a$};
    \draw[solid] (axis cs:rholeft,uleft) -- (axis cs:0,uleft)
        node[pos=0.85, above] {\footnotesize $R_a$};

    \addplot [
        domain=0.01:15,
        samples=100,
        color=black,
        solid
    ] {(cont_disc_limit+0.29)*exp(-1/(0.0205*(x+5)*(x+5)))-0.29}
        node[pos=0.85, below] {\footnotesize $S_\delta$};

    \node at (axis cs:11, 7.5) {\footnotesize $I_0$};
        \draw[black,-latex,very thin] (axis cs:11.05, 7) -- (axis cs:12.5,3.25);
    
    \node at (axis cs:2.5, 7.5) {\footnotesize $II_0$};
        \draw[black,-latex,very thin] (axis cs:2.5, 7) -- (axis cs:2.5,2.5);
        \draw[black,-latex,very thin] (axis cs:2.85, 7) -- (axis cs:6,2.6);
        \draw[black,-latex,very thin] (axis cs:3,7.45) -- (axis cs:5.75,7.4);
    
    \node at (axis cs:2, -5) {\footnotesize $V_{R_{a}VC_{0}}$};
    
    \node at (axis cs:11.1,-5) {\footnotesize $IV$};
    \draw[black,-latex,very thin] (axis cs:10.5, -5) -- (axis cs:4.25,-2.5);
    \draw[black,-latex,very thin] (axis cs:11.1, -4.4) -- (axis cs:7.6,1);
    \draw[black,-latex,very thin] (axis cs:10.7, -4.6) -- (axis cs:4,0.8);

    \node at (axis cs:1.75,1.25) {\footnotesize $IV_{\delta}$};
    \end{axis}
\end{tikzpicture}
    \end{minipage}
    \hfill
    \begin{minipage}[b]{.45\textwidth}
        \caption{\\
            Case 4 (above left): $\lambda_a > 0 > \lambda_0$ \\
            Case 5 (above right): $\lambda_a > \lambda_0 = 0$ \\
            Case 6 (left): $\lambda_a > \lambda_0 > 0$
        }
        \label{fig:cases4-6}
        \vspace{2cm}
    \end{minipage}
\end{figure}
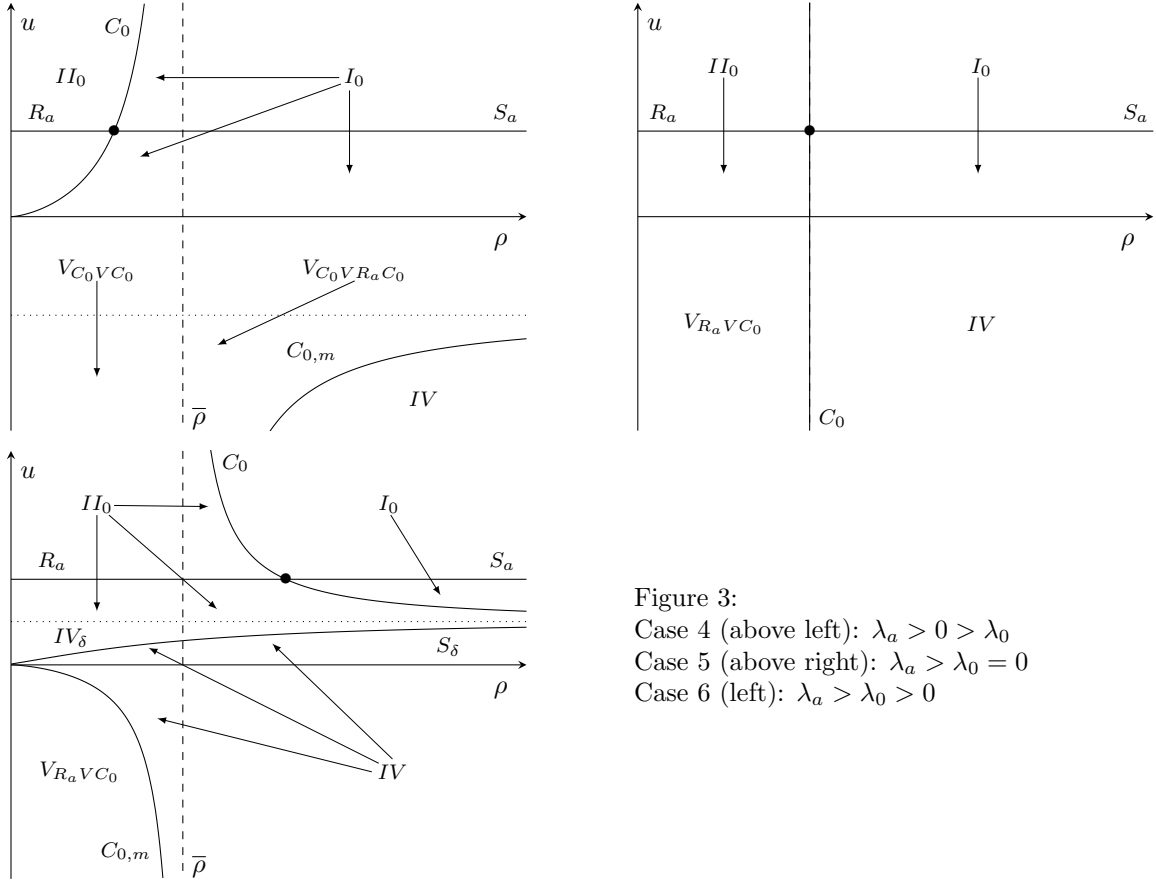

\subsubsection{\texorpdfstring{$-1<a<0$}{-1 < a < 0}}
We note that the structure of the $-1<a<0$ cases are similar to the structure of the $a<-1$ cases. The only difference is that the portions of the line $\u =\u_L$ have now reversed roles from $S_a$ to $R_a$ and vice versa. We identify three more cases: 
\begin{itemize}
    \item Case 13: $u_L<0$ and $\rho_L<\rhobar$,
    \item Case 14: $u_L<0$ and $\rho_L=\rhobar$, and
    \item Case 15: $u_L<0$ and $\rho_L>\rhobar$.
\end{itemize}
We now discuss the solutions to the Riemann Problems in the various regions and detail the path that is taken to reach the right state. For Cases 13 - 15, we have the following regions:
\begin{itemize}
    \item Region $I_a$: The left state is connected to the right by an $a$-shock followed by a $0$-contact discontinuity.
    \item Region $II_a$: The left state is connected to the right state by a rarefaction curve $R_a$ followed by a $C_0$.
    \item Region $IV$: The left state follows an overcompressive delta shock $S_\delta$ as $\rho$ blows up.
    \item Region $VI$: The left state will either take an $a$-shock or rarefaction curve to first reach $\rhobar$, from Case 13 and Case 15, respectively, going into Case 14. It will then travel along $\rhobar$ via a $0$-contact discontinuity, finishing with an $a$-shock or $a$-rarefaction curve to the right state.
\end{itemize}
Below, we present the state spaces for Cases 13 through 15.

\begin{figure}[H]
    \begin{minipage}[t]{.5\textwidth}
        \vspace{0pt}
        \centering
        \begin{tikzpicture}[
        declare function={
            a = -.5;       
            rhowithbar = 5; 
            rholeft = 3;    
            uleft = -4;      
            A = 0;          
            cont_disc_rholeftnotrhowithbar(\t)
                = uleft + (uleft - A)*(exp(ln(rholeft)*(a)) - exp(ln(\t)*(a)))/(exp(ln(\t)*(a)) - exp(ln(rhowithbar)*(a)));
            cont_disc_limit
                = uleft + (uleft - A)*(exp(ln(rholeft)*(a)))/(- exp(ln(rhowithbar)*(a)));
            cont_disc_limiting_curve(\t)
                = uleft + (uleft - A)*(exp(ln(1000)*(a)) - exp(ln(\t)*(a)))/(exp(ln(\t)*(a)) - exp(ln(rhowithbar)*(a)));
        }
    ]

    \tikzstyle{arrow} = [thick,->,>=stealth]

    \begin{axis}[
        axis lines = center,
        xmin=0, xmax=15,
        ymin=-10, ymax=10,
        xlabel=$\rho$,
        ylabel=$\u$,
        x label style = {at={(axis description cs:1.02,0.45)},anchor=south},
        y label style = {at={(axis description cs:0.0,0.95)},anchor=west},
        xticklabel=\empty,
        yticklabel=\empty,
        major tick style={draw=none}
    ]

    \node at (axis cs:rholeft,uleft) {\textbullet};

    \draw[dashed] (axis cs:rhowithbar,10) -- (axis cs:rhowithbar,-10)
        node[pos=0.97, right] {$\overline{\rho}$};

    \addplot [
        domain=0.01:(rhowithbar-1),
        samples=100,
        color=black,
        solid
    ] {cont_disc_rholeftnotrhowithbar(x)}
        node[pos=0.9, left] {\footnotesize $C_{0}$};
    \addplot [
        domain=(rhowithbar+1):15,
        samples=100,
        color=black,
        solid
    ] {cont_disc_rholeftnotrhowithbar(x)-1}
        node[pos=0.35, left] {\footnotesize $C_{0,m}$};
    \addplot [
        domain=(rhowithbar+2):15,
        samples=100,
        color=black,
        solid
    ] {cont_disc_limiting_curve(x+3)+3}
        node[pos=0.4, above] {\footnotesize $C_{0,\ell}$};
    \draw[dotted] (axis cs:0,cont_disc_limit) -- (axis cs:15,cont_disc_limit);

    \draw[solid] (axis cs:rholeft,uleft) -- (axis cs:15,uleft)
        node[pos=0.95, above] {\footnotesize $S_a$};
    \draw[solid] (axis cs:rholeft,uleft) -- (axis cs:0,uleft)
        node[pos=0.85, above, xshift=5] {\footnotesize $R_a$};

    \node at (axis cs:12,-8) {\footnotesize $IV$};
    
    \node at (axis cs:1.75,-7) {\footnotesize $II_a$};
    
    \node at (axis cs:4,-1.75) {\footnotesize $I_a$};
        \draw[black,-latex,very thin] (axis cs:4,-2.55) -- (axis cs:4.25, -6);
        \draw[black,-latex,very thin] (axis cs:4.5,-2) -- (axis cs:10,-2);
        \draw[black,-latex,very thin] (axis cs:4.4,-2.4) -- (axis cs:7.5,-6);
    
    \node at (axis cs:2.5,5) {\footnotesize $VI$};
        \draw[black,-latex,very thin] (axis cs:3,5) -- (axis cs:10,5);
        \draw[black,-latex,very thin] (axis cs:3,4.5) -- (axis cs:7.5,2);
    \end{axis}
\end{tikzpicture}
    \end{minipage}
    \hfill
    \begin{minipage}[t]{.5\textwidth}
        \vspace{0pt}
        \centering
        \begin{tikzpicture}[
        declare function={
            a = -.5;       
            rhowithbar = 5; 
            rholeft = 5;    
            uleft = -4;      
            A = 0;          
            cont_disc_rholeftnotrhowithbar(\t)
                = uleft + (uleft - A)*(exp(ln(rholeft)*(a)) - exp(ln(\t)*(a)))/(exp(ln(\t)*(a)) - exp(ln(rhowithbar)*(a)));
            cont_disc_limit
                = uleft + (uleft - A)*(exp(ln(rholeft)*(a)))/(- exp(ln(rhowithbar)*(a)));
            cont_disc_limiting_curve(\t)
                = uleft + (uleft - A)*(exp(ln(1000)*(a)) - exp(ln(\t)*(a)))/(exp(ln(\t)*(a)) - exp(ln(rhowithbar)*(a)));
        }
    ]

    \tikzstyle{arrow} = [thick,->,>=stealth]

    \begin{axis}[
        axis lines = center,
        xmin=0, xmax=15,
        ymin=-10, ymax=10,
        xlabel=$\rho$,
        ylabel=$\u$,
        x label style = {at={(axis description cs:1.02,.45)},anchor=south},
        y label style = {at={(axis description cs:0.0,0.95)},anchor=west},
        xticklabel=\empty,
        yticklabel=\empty,
        major tick style={draw=none}
    ]
    \node at (axis cs:rholeft,uleft) {\textbullet};

    \draw[dashed] (axis cs:rhowithbar,10) -- (axis cs:rhowithbar,-10);

    \addplot [
        domain=(rhowithbar+2):15,
        samples=100,
        color=black,
        solid
    ] {cont_disc_limiting_curve(x+2)+3}
        node[pos=0.45, above] {\footnotesize $C_{0,\ell}$};
    \draw[solid] (axis cs:rhowithbar,10) -- (axis cs:rhowithbar,-10)
        node[pos=0.97, right] {\footnotesize $C_0$};;

    \draw[solid] (axis cs:rholeft,uleft) -- (axis cs:15,uleft)
        node[pos=0.95, above] {\footnotesize $S_a$};
    \draw[solid] (axis cs:rholeft,uleft) -- (axis cs:0,uleft)
        node[pos=0.85, above] {\footnotesize $R_a$};

     \node at (axis cs:2.5,-7.5) {\footnotesize $II_a$};
        \draw[black,-latex,very thin] (axis cs:2.5,-6.8) -- (axis cs:2.5,-2);
        
    \node at (axis cs:12.5,-8.5) {\footnotesize $IV$};

    \node at (axis cs:2.5,5) {\footnotesize $VI$};
    
    \node at (axis cs:10,5) {\footnotesize $VI$};

    \node at (axis cs:10, -2) {\footnotesize $I_a$};
    
    \draw[black,-latex,very thin] (axis cs:9.5,-2.5) -- (axis cs:7,-6);
  
    \end{axis}
\end{tikzpicture}
    \end{minipage}
    \!\!\!\!\!\!\!\!\begin{minipage}[b]{.5\textwidth}
        \vspace{0pt}
        \centering
        \begin{tikzpicture}[
        declare function={
            a = -.5;       
            rhowithbar = 5; 
            rholeft = 8;    
            uleft = -4;      
            A = 0;          
            cont_disc_rholeftnotrhowithbar(\t)
                = uleft + (uleft - A)*(exp(ln(rholeft)*(a)) - exp(ln(\t)*(a)))/(exp(ln(\t)*(a)) - exp(ln(rhowithbar)*(a)));
            cont_disc_limit
                = uleft + (uleft - A)*(exp(ln(rholeft)*(a)))/(- exp(ln(rhowithbar)*(a)));
            cont_disc_limiting_curve(\t)
                = uleft + (uleft - A)*(exp(ln(1000)*(a)) - exp(ln(\t)*(a)))/(exp(ln(\t)*(a)) - exp(ln(rhowithbar)*(a)));
        }
    ]

    \tikzstyle{arrow} = [thick,->,>=stealth]

    \begin{axis}[
        axis lines = center,
        xmin=0, xmax=15,
        ymin=-10, ymax=10,
        xlabel=$\rho$,
        ylabel=$\u$,
        x label style = {at={(axis description cs:1.02,0.45)},anchor=south},
        y label style = {at={(axis description cs:0.0,0.95)},anchor=west},
        xticklabel=\empty,
        yticklabel=\empty,
        major tick style={draw=none}
    ]

    \node at (axis cs:rholeft,uleft) {\textbullet};

    \draw[dashed] (axis cs:rhowithbar,10) -- (axis cs:rhowithbar,-10)
        node[pos=0.97, right] {$\overline{\rho}$};

    \addplot [
        domain=0.01:(rhowithbar-0.58),
        samples=100,
        color=black,
        solid
    ] {cont_disc_rholeftnotrhowithbar(x)}
        node[pos=0.7, left] {\footnotesize $C_{0,m}$};
    \addplot [
        domain=(rhowithbar+0.8):15,
        samples=100,
        color=black,
        solid
    ] {cont_disc_rholeftnotrhowithbar(x)}
        node[pos=0.8, above] {\footnotesize $C_{0}$};
    \addplot [
        domain=(rhowithbar+.5):15,
        samples=100,
        color=black,
        solid
    ] {cont_disc_limiting_curve(x+.5)+4}
        node[pos=.8, above, xshift= -2] {\footnotesize $C_{0,\ell}$};
    \draw[dotted] (axis cs:0,cont_disc_limit) -- (axis cs:15,cont_disc_limit);

    \draw[solid] (axis cs:rholeft,uleft) -- (axis cs:15,uleft)
        node[pos=0.95, above] {\footnotesize $S_a$};
    \draw[solid] (axis cs:rholeft,uleft) -- (axis cs:0,uleft)
        node[pos=0.85, above] {\footnotesize $R_a$};

    \node at (axis cs:12.5,-8.5) {\footnotesize $IV$};
    
    \node at (axis cs:2,5) {\footnotesize $VI$};
        \draw[black,-latex,very thin] (axis cs:2.1,4.25) -- (axis cs:4,1.25);

    \node at (axis cs:10,5) {\footnotesize $VI$};
    
    \node at (axis cs:2.5,-7) {\footnotesize $II_a$};
        \draw[black,-latex,very thin] (axis cs:2.5,-6.4) -- (axis cs:2.5,-2);
        \draw[black,-latex,very thin] (axis cs:2.75,-6.5) -- (axis cs:8,-2);
        \draw[black,-latex,very thin] (axis cs:3,-7) -- (axis cs:5.75,-7);

    \node at (axis cs:10,-5) {\footnotesize $I_a$};
    
    \end{axis}
\end{tikzpicture}
    \end{minipage}
    \hfill
    \begin{minipage}[b]{.45\textwidth}
        \caption{\\
            Case 13 (above left): $\lambda_a < \lambda_0$ and $0< \lambda_0$ \\
            Case 14 (above right): $\lambda_a < \lambda_0 = 0$ \\
            Case 15 (left): $\lambda_a < \lambda_0 < 0$
        }
        \label{fig:cases13-15}
        \vspace{2cm}
    \end{minipage}
\end{figure}

Cases 16 - 18 are similar:
\begin{itemize}
    \item Case 16: $u_L>0$ and $\rho_L<\rhobar$,
    \item Case 17: $u_L>0$ and $\rho_L=\rhobar$, and
    \item Case 18: $u_L>0$ and $\rho_L>\rhobar$.
\end{itemize}
These are the versions of Cases 13, 14, and 15 where $\u_L>0$. Additionally, a vacuum region exists in these cases where $\rho$ hits zero. The regions for the solution to the Riemann problem are as follows:
\begin{itemize}
    \item Regions $I_0$, $II_0$, and $IV$ are identical to the Cases 13 - 15 except that the order of the shock/rarefaction and contact discontinuity is reversed in $I_0$ and $II_0$.
    \item Region $V$: Just like in Cases 4 – 6, this occurs when the path from the left state to the right state passes through a vacuum (i.e., when $\rho = 0$), and the transition often involves a sequence of waves. The only difference from Cases 4 - 6 is that, instead of following an $R_a$, the left state in these cases follows an $S_a$.
\end{itemize}

Below are the state spaces for Cases 16 - 18:
\begin{figure}[H]
    \ \ \begin{minipage}[t]{.5\textwidth}
        \vspace{0pt}
        \centering
        \begin{tikzpicture}[
        declare function={
            a = -.5;       
            rhowithbar = 5; 
            rholeft = 3;    
            uleft = 4;      
            A = 0;          
            cont_disc_rholeftnotrhowithbar(\t)
                = uleft + (uleft - A)*(exp(ln(rholeft)*(a)) - exp(ln(\t)*(a)))/(exp(ln(\t)*(a)) - exp(ln(rhowithbar)*(a)));
            cont_disc_limit
                = uleft + (uleft - A)*(exp(ln(rholeft)*(a)))/(- exp(ln(rhowithbar)*(a)));
            cont_disc_limiting_curve(\t)
                = uleft + (uleft - A)*(exp(ln(.01)*(a)) - exp(ln(\t)*(a)))/(exp(ln(\t)*(a)) - exp(ln(rhowithbar)*(a)));
        }
    ]

    \tikzstyle{arrow} = [thick,->,>=stealth]

    \begin{axis}[
        axis lines = center,
        xmin=0, xmax=15,
        ymin=-11.5, ymax=10,
        xlabel=$\rho$,
        ylabel=$\u$,
        x label style = {at={(axis description cs:1.02,0.5)},anchor=south},
        y label style = {at={(axis description cs:0.0,0.95)},anchor=west},
        xticklabel=\empty,
        yticklabel=\empty,
        major tick style={draw=none}
    ]

    \node at (axis cs:rholeft,uleft) {\textbullet};

    \draw[dashed] (axis cs:rhowithbar,10) -- (axis cs:rhowithbar,-10)
        node[pos=1, below] {$\overline{\rho}$};;

    \addplot [
        domain=0.01:(rhowithbar-.5),
        samples=100,
        color=black,
        solid
    ] {cont_disc_rholeftnotrhowithbar(x)}
        node[pos=0.7, left] {\footnotesize $C_{0}$};
    \addplot [
        domain=(rhowithbar+1):15,
        samples=100,
        color=black,
        solid
    ] {cont_disc_rholeftnotrhowithbar(x)}
        node[pos=0.35, left] {\footnotesize $C_{0,m}$};

    \draw[dotted] (axis cs:0,cont_disc_limit) -- (axis cs:15,cont_disc_limit);
    \draw[dotted] (axis cs:0,cont_disc_limit) -- (axis cs:15,cont_disc_limit);

    \draw[solid] (axis cs:rholeft,uleft) -- (axis cs:15,uleft)
        node[pos=0.95, above] {\footnotesize $R_a$};
    \draw[solid] (axis cs:rholeft,uleft) -- (axis cs:0,uleft)
        node[pos=.4, above] {\footnotesize $S_a$};

    \node at (axis cs:1.75,6.5) {\footnotesize $I_0$};
    
    \node at (axis cs:10,6.5) {\footnotesize $II_0$};
        \draw[black,-latex,very thin] (axis cs:9.65,6.15) -- (axis cs:3.75,3);
        \draw[black,-latex,very thin] (axis cs:9.85,5.9) -- (axis cs:9.85,2);
        \draw[black,-latex,very thin] (axis cs:9.5,6.5) -- (axis cs:4.2, 6.5);

    \node at (axis cs:12,-8.5) {\footnotesize $IV$};
    
    \node at (axis cs:2.5, -5) {\footnotesize $V_{C_{0}VC_{0}}$};
        \draw[black,-latex,very thin] (axis cs:2.5,-4.4) -- (axis cs:2.5,-0.5);

    \node at (axis cs:10, -2.5) {\footnotesize $V_{C_{0}VS_{a}C_{0}}$};
        \draw[black,-latex,very thin] (axis cs:10,-2) -- (axis cs:10,-0.5);

    \end{axis}
\end{tikzpicture}
    \end{minipage}
    \hfill
    \begin{minipage}[t]{.5\textwidth}
        \vspace{0pt}
        \centering
        \begin{tikzpicture}[
        declare function={
            a = -.5;       
            rhowithbar = 5; 
            rholeft = 5;    
            uleft = 4;      
            A = 0;          
            cont_disc_rholeftnotrhowithbar(\t)
                = uleft + (uleft - A)*(exp(ln(rholeft)*(a)) - exp(ln(\t)*(a)))/(exp(ln(\t)*(a)) - exp(ln(rhowithbar)*(a)));
            cont_disc_limit
                = uleft + (uleft - A)*(exp(ln(rholeft)*(a)))/(- exp(ln(rhowithbar)*(a)));
            cont_disc_limiting_curve(\t)
                = uleft + (uleft - A)*(exp(ln(1000)*(a)) - exp(ln(\t)*(a)))/(exp(ln(\t)*(a)) - exp(ln(rhowithbar)*(a)));
        }
    ]

    \tikzstyle{arrow} = [thick,->,>=stealth]

    \begin{axis}[
        axis lines = center,
        xmin=0, xmax=15,
        ymin=-10, ymax=10,
        xlabel=$\rho$,
        ylabel=$\u$,
        x label style = {at={(axis description cs:1.02,0.45)},anchor=south},
        y label style = {at={(axis description cs:0.0,0.95)},anchor=west},
        xticklabel=\empty,
        yticklabel=\empty,
        major tick style={draw=none}
    ]

    \node at (axis cs:rholeft,uleft) {\textbullet};

    \draw[dashed] (axis cs:rhowithbar,10) -- (axis cs:rhowithbar,-10);

    \draw[solid] (axis cs:rhowithbar,10) -- (axis cs:rhowithbar,-10)
        node[pos=0.97, right] {\footnotesize $C_0$};;

    \draw[solid] (axis cs:rholeft,uleft) -- (axis cs:15,uleft)
        node[pos=0.95, above] {\footnotesize $R_a$};
    \draw[solid] (axis cs:rholeft,uleft) -- (axis cs:0,uleft)
        node[pos=0.85, above] {\footnotesize $S_a$};

    \node at (axis cs:2.5,7) {\footnotesize $I_0$};
        \draw[black,-latex,very thin] (axis cs:2.5,6.5) -- (axis cs:2.5,2);
    
    \node at (axis cs:10,7) {\footnotesize $II_0$};
        \draw[black,-latex,very thin] (axis cs:9.9,6.5) -- (axis cs:9.9,2);
    
    \node at (axis cs:2.5,-5) {\footnotesize $V_{S_{a}VC_{0}}$};
    
    \node at (axis cs:10,-5) {\footnotesize $IV$};
    \end{axis}
\end{tikzpicture}
    \end{minipage}
    \!\!\!\!\!\!\!\!\begin{minipage}[b]{.5\textwidth}
        \vspace{0pt}
        \centering
        \begin{tikzpicture}[
        declare function={
            a = -0.5;       
            rhowithbar = 5; 
            rholeft = 8;    
            uleft = 4;      
            A = 0;          
            cont_disc_rholeftnotrhowithbar(\t)
                = uleft + (uleft - A)*(exp(ln(rholeft)*(a)) - exp(ln(\t)*(a)))/(exp(ln(\t)*(a)) - exp(ln(rhowithbar)*(a)));
            cont_disc_limit
                = uleft + (uleft - A)*(exp(ln(rholeft)*(a)))/(- exp(ln(rhowithbar)*(a)));
            cont_disc_limiting_curve(\t)
                = uleft + (uleft - A)*(exp(ln(1000)*(a)) - exp(ln(\t)*(a)))/(exp(ln(\t)*(a)) - exp(ln(rhowithbar)*(a)));
        }
    ]

    \tikzstyle{arrow} = [thick,->,>=stealth]

    \begin{axis}[
        axis lines = center,
        xmin=0, xmax=15,
        ymin=-10, ymax=10,
        xlabel=$\rho$,
        ylabel=$\u$,
        x label style = {at={(axis description cs:0.95,0.4)},anchor=south},
        y label style = {at={(axis description cs:0.0,0.95)},anchor=west},
        xticklabel=\empty,
        yticklabel=\empty,
        major tick style={draw=none}
    ]

    \node at (axis cs:rholeft,uleft) {\textbullet};

    \draw[dashed] (axis cs:rhowithbar,10) -- (axis cs:rhowithbar,-10)
        node[pos=0.97, right] {$\overline{\rho}$};;

    \addplot [
        domain=(rhowithbar+0.20):15,
        samples=100,
        color=black,
        solid
    ] {cont_disc_rholeftnotrhowithbar(x)}
        node[pos=0.55, right] {\footnotesize $C_{0}$};
    \addplot [
        domain=0.01:(rhowithbar-0.20),
        samples=100,
        color=black,
        solid
    ] {cont_disc_rholeftnotrhowithbar(x)}
        node[pos=0.25, left] {\footnotesize $C_{0,m}$};
    \draw[dotted] (axis cs:0,cont_disc_limit) -- (axis cs:15,cont_disc_limit);

    \draw[solid] (axis cs:rholeft,uleft) -- (axis cs:15,uleft)
        node[pos=0.95, above] {\footnotesize $R_a$};
    \draw[solid] (axis cs:rholeft,uleft) -- (axis cs:0,uleft)
        node[pos=0.85, above] {\footnotesize $S_a$};

    \node at (axis cs:11, 7.5) {\footnotesize $II_0$};
        \draw[black,-latex,very thin] (axis cs:11.05, 7) -- (axis cs:12.5,3.25);
    
    \node at (axis cs:2.5, 7.5) {\footnotesize $I_0$};
        \draw[black,-latex,very thin] (axis cs:2.5, 7) -- (axis cs:2.5,2.5);
        \draw[black,-latex,very thin] (axis cs:2.85, 7) -- (axis cs:6,2.6);
        \draw[black,-latex,very thin] (axis cs:2.9,7.4) -- (axis cs:5.75,7.4);
    
    \node at (axis cs:2, -5) {\footnotesize $V_{S_{a}VC_{0}}$};
    
    \node at (axis cs:11.1,-5) {\footnotesize $IV$};
    \draw[black,-latex,very thin] (axis cs:10.5, -5) -- (axis cs:3.9,-2.5);
    \draw[black,-latex,very thin] (axis cs:11.1, -4.4) -- (axis cs:7.6,.5);
    \draw[black,-latex,very thin] (axis cs:10.5, -4.55) -- (axis cs:3.9,.5);

    \end{axis}
\end{tikzpicture}
    \end{minipage}
    \hfill
    \begin{minipage}[b]{.45\textwidth}
        \caption{\\
            Case 16: $\lambda_a > \lambda_0$ and $0 > \lambda_0$ \\
            Case 17: $\lambda_a > \lambda_0 = 0$ \\
            Case 18: $\lambda_a > \lambda_0 > 0$
        }
        \label{fig:cases16-18}
        \vspace{2cm}
    \end{minipage}
\end{figure}
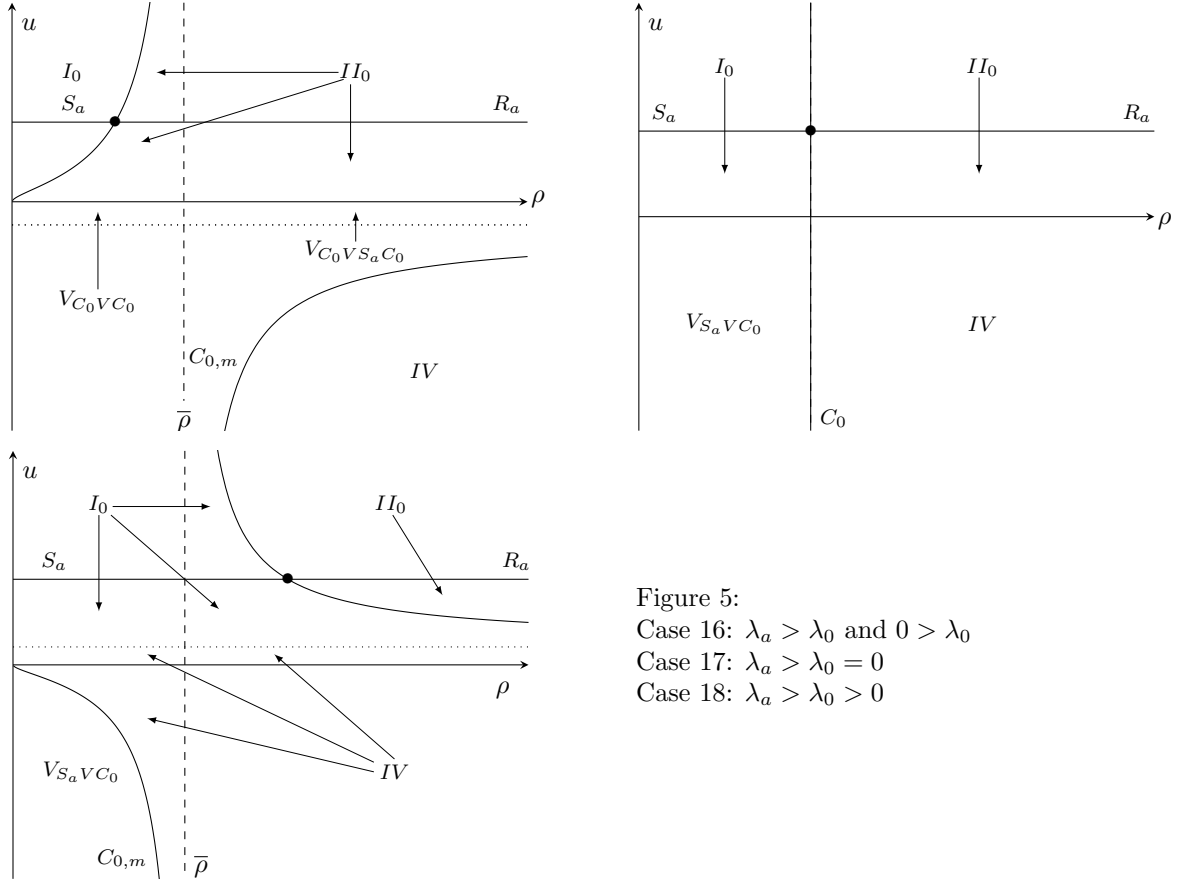

\subsubsection{\texorpdfstring{$a=-1$}{a=-1}}
We again distinguish three cases:
\begin{itemize}
    \item Case 7: $u_L<0$ and $\rho_L<\rhobar$,
    \item Case 8: $u_L<0$ and $\rho_L=\rhobar$, and
    \item Case 9: $u_L<0$ and $\rho_L>\rhobar$.
\end{itemize}
We now discuss the various regions for which different solutions exist and detail the path that is taken to reach the right state.
For Cases 7 - 9, we have the following regions:
\begin{itemize}
    \item Region $III_a$: The left state is connected to the right by an $a$-contact discontinuity followed by a $0$-contact discontinuity.
    \item Region $IV$: The overcompressive region is similar to other negative values of $a$.
    \item Region $VI$: The left state will take an $a$-contact discontinuity to first reach $\rhobar$ from Case 7 or Case 9, going into Case 8. It will then travel along $\rhobar$ via a $0$-contact discontinuity, followed by an $a$-contact discontinuity to the right state.
\end{itemize}
Below are the state spaces for Cases 7 - 9:

\begin{figure}[H]
    \begin{minipage}[t]{.5\textwidth}
        \vspace{0pt}
        \centering
        \begin{tikzpicture}[
        declare function={
            a = -1;       
            rhowithbar = 5; 
            rholeft = 3;    
            uleft = -4;      
            A = 0;          
            cont_disc_rholeftnotrhowithbar(\t)
                = uleft + (uleft - A)*(exp(ln(rholeft)*(a)) - exp(ln(\t)*(a)))/(exp(ln(\t)*(a)) - exp(ln(rhowithbar)*(a)));
            cont_disc_limit
                = uleft + (uleft - A)*(exp(ln(rholeft)*(a)))/(- exp(ln(rhowithbar)*(a)));
            cont_disc_limiting_curve(\t)
                = uleft + (uleft - A)*(exp(ln(1000)*(a)) - exp(ln(\t)*(a)))/(exp(ln(\t)*(a)) - exp(ln(rhowithbar)*(a)));
        }
    ]

    \tikzstyle{arrow} = [thick,->,>=stealth]

    \begin{axis}[
        axis lines = center,
        xmin=0, xmax=15,
        ymin=-10, ymax=10,
        xlabel=$\rho$,
        ylabel=$\u$,
        x label style = {at={(axis description cs:0.95,0.4)},anchor=south},
        y label style = {at={(axis description cs:0.0,0.95)},anchor=west},
        xticklabel=\empty,
        yticklabel=\empty,
        major tick style={draw=none}
    ]

    \node at (axis cs:rholeft,uleft) {\textbullet};

    \draw[dashed] (axis cs:rhowithbar,10) -- (axis cs:rhowithbar,-10)
        node[pos=0.97, right] {$\overline{\rho}$};

    \addplot [
        domain=0.01:(rhowithbar-1.05),
        samples=100,
        color=black,
        solid
    ] {cont_disc_rholeftnotrhowithbar(x)}
        node[pos=0.7, left] {\footnotesize $C_{0}$};
    \addplot [
        domain=(rhowithbar+1.8):15,
        samples=100,
        color=black,
        solid
    ] {cont_disc_rholeftnotrhowithbar(x)}
        node[pos=0.2, left] {\footnotesize $C_{0,m}$};
    \addplot [
        domain=(rhowithbar+3.3):15,
        samples=100,
        color=black,
        solid
    ] {cont_disc_limiting_curve(x)}
        node[pos=0.2, left] {\footnotesize $C_{0,\ell}$};
    \draw[dotted] (axis cs:0,cont_disc_limit) -- (axis cs:15,cont_disc_limit);

    \draw[solid] (axis cs:rholeft,uleft) -- (axis cs:15,uleft)
        node[pos=0.95, above] {\footnotesize $C_a$};
    \draw[solid] (axis cs:rholeft,uleft) -- (axis cs:0,uleft);

    \node at (axis cs:12.5,-9) {\footnotesize $IV$};
    
    \node at (axis cs:2.5,6) {\footnotesize $VI$};
        \draw[black,-latex,very thin] (axis cs:3,6) -- (axis cs:12.5,6);
        \draw[black,-latex,very thin] (axis cs:3,5.75) -- (axis cs:9,2);
        
    \node at (axis cs:10,-3) {\footnotesize $III_a$};
        \draw[black,-latex,very thin] (axis cs:9.25,-3) -- (axis cs:1.5,-3);
        \draw[black,-latex,very thin] (axis cs:9.3,-3.35) -- (axis cs:1.5,-8.5);
    \end{axis}
\end{tikzpicture}
    \end{minipage}
    \hfill
    \begin{minipage}[t]{.5\textwidth}
        \vspace{0pt}
        \centering
        \begin{tikzpicture}[
        declare function={
            a = -1;       
            rhowithbar = 5; 
            rholeft = 5;    
            uleft = -4;      
            A = 0;          
            cont_disc_rholeftnotrhowithbar(\t)
                = uleft + (uleft - A)*(exp(ln(rholeft)*(a)) - exp(ln(\t)*(a)))/(exp(ln(\t)*(a)) - exp(ln(rhowithbar)*(a)));
            cont_disc_limit
                = uleft + (uleft - A)*(exp(ln(rholeft)*(a)))/(- exp(ln(rhowithbar)*(a)));
            cont_disc_limiting_curve(\t)
                = uleft + (uleft - A)*(exp(ln(1000)*(a)) - exp(ln(\t)*(a)))/(exp(ln(\t)*(a)) - exp(ln(rhowithbar)*(a)));
        }
    ]

    \tikzstyle{arrow} = [thick,->,>=stealth]

    \begin{axis}[
        axis lines = center,
        xmin=0, xmax=15,
        ymin=-10, ymax=10,
        xlabel=$\rho$,
        ylabel=$\u$,
        x label style = {at={(axis description cs:0.95,0.4)},anchor=south},
        y label style = {at={(axis description cs:0.0,0.95)},anchor=west},
        xticklabel=\empty,
        yticklabel=\empty,
        major tick style={draw=none}
    ]

    \node at (axis cs:rholeft,uleft) {\textbullet};

    \draw[dashed] (axis cs:rhowithbar,10) -- (axis cs:rhowithbar,-10);

    \addplot [
        domain=(rhowithbar+3.2):15,
        samples=100,
        color=black,
        solid
    ] {cont_disc_limiting_curve(x)}
        node[pos=0.2, left] {\footnotesize $C_{0,\ell}$};
    \draw[solid] (axis cs:rhowithbar,10) -- (axis cs:rhowithbar,-10)
        node[pos=0.97, right] {\footnotesize $C_0$};;

    \draw[solid] (axis cs:rholeft,uleft) -- (axis cs:15,uleft)
        node[pos=0.95, above] {\footnotesize $C_a$};
    \draw[solid] (axis cs:rholeft,uleft) -- (axis cs:0,uleft);

    \node at (axis cs:12.5,-9) {\footnotesize $IV$};
    
    \node at (axis cs:2.5,5) {\footnotesize $III_0$};
    
    \node at (axis cs:10,5) {\footnotesize $III_0$};
    
    \node at (axis cs:2.5,-2) {\footnotesize $III_a$};
        \draw[black,-latex,very thin] (axis cs:2.5,-2.5) -- (axis cs:2.5,-7.5);
        \draw[black,-latex,very thin] (axis cs:3.25,-2) -- (axis cs:10,-2);
        \draw[black,-latex,very thin] (axis cs:3.25,-2.4) -- (axis cs:7,-5.85);
    \end{axis}
\end{tikzpicture}
    \end{minipage}
    \!\!\!\!\!\!\!\begin{minipage}[b]{.5\textwidth}
        \vspace{0pt}
        \centering
        \begin{tikzpicture}[
        declare function={
            a = -1;       
            rhowithbar = 5; 
            rholeft = 8;    
            uleft = -4;      
            A = 0;          
            cont_disc_rholeftnotrhowithbar(\t)
                = uleft + (uleft - A)*(exp(ln(rholeft)*(a)) - exp(ln(\t)*(a)))/(exp(ln(\t)*(a)) - exp(ln(rhowithbar)*(a)));
            cont_disc_limit
                = uleft + (uleft - A)*(exp(ln(rholeft)*(a)))/(- exp(ln(rhowithbar)*(a)));
            cont_disc_limiting_curve(\t)
                = uleft + (uleft - A)*(exp(ln(1000)*(a)) - exp(ln(\t)*(a)))/(exp(ln(\t)*(a)) - exp(ln(rhowithbar)*(a)));
        }
    ]

    \tikzstyle{arrow} = [thick,->,>=stealth]

    \begin{axis}[
        axis lines = center,
        xmin=0, xmax=15,
        ymin=-10, ymax=10,
        xlabel=$\rho$,
        ylabel=$\u$,
        x label style = {at={(axis description cs:0.95,0.4)},anchor=south},
        y label style = {at={(axis description cs:0.0,0.95)},anchor=west},
        xticklabel=\empty,
        yticklabel=\empty,
        major tick style={draw=none}
    ]

    \node at (axis cs:rholeft,uleft) {\textbullet};

    \draw[dashed] (axis cs:rhowithbar,10) -- (axis cs:rhowithbar,-10)
        node[pos=0.97, right] {$\overline{\rho}$};

    \addplot [
        domain=0.01:(rhowithbar-0.58),
        samples=100,
        color=black,
        solid
    ] {cont_disc_rholeftnotrhowithbar(x)}
        node[pos=0.7, left] {\footnotesize $C_{0,m}$};
    \addplot [
        domain=(rhowithbar+0.8):15,
        samples=100,
        color=black,
        solid
    ] {cont_disc_rholeftnotrhowithbar(x)}
        node[pos=0.3, left] {\footnotesize $C_{0}$};
    \addplot [
        domain=(rhowithbar+3.3):15,
        samples=100,
        color=black,
        solid
    ] {cont_disc_limiting_curve(x)}
        node[pos=0.2, left] {\footnotesize $C_{0,\ell}$};
    \draw[dotted] (axis cs:0,cont_disc_limit) -- (axis cs:15,cont_disc_limit);

    \draw[solid] (axis cs:rholeft,uleft) -- (axis cs:15,uleft)
        node[pos=0.95, below] {\footnotesize $C_a$};
    \draw[solid] (axis cs:rholeft,uleft) -- (axis cs:0,uleft);

    \node at (axis cs:12.5,-9) {\footnotesize $IV$};
    
   \node at (axis cs:2,5) {\footnotesize $VI$};
        \draw[black,-latex,very thin] (axis cs:2.7,5) -- (axis cs:12.5,5);
        \draw[black,-latex,very thin] (axis cs:2.6,4.85) -- (axis cs:4.5,2);
        
    \node at (axis cs:2,-1.5) {\footnotesize $III_a$};
        \draw[black,-latex,very thin] (axis cs:2,-2) -- (axis cs:2,-7.5);
        \draw[black,-latex,very thin] (axis cs:3,-1.5) -- (axis cs:10,-1.5);
       \draw[black,-latex,very thin] (axis cs:3,-2) -- (axis cs:8,-7);  
    \end{axis}
\end{tikzpicture}
    \end{minipage}
    \hfill
    \begin{minipage}[b]{.45\textwidth}
        \caption{\\
            Case 7 (above left): $\lambda_a < 0 < \lambda_0$ \\
            Case 8 (abover right): $\lambda_a < \lambda_0 = 0$ \\
            Case 9 (left)   : $\lambda_a < \lambda_0 < 0$
        }
        \label{fig:cases7-9}
        \vspace{2cm}
    \end{minipage}
\end{figure}

We then turn to the following three cases:
\begin{itemize}
    \item Case 10: $\u_L>0$ and $\rho_L<\rhobar$, 
    \item Case 11: $\u_L>0$ and $\rho_L=\rhobar$, and
    \item Case 12: $\u_L>0$ and $\rho_L>\rhobar$.
\end{itemize}
These are analogous to Cases 7, 8, and 9 except that $\u_L>0$, and also to Case 4, 5, 6, and Cases 16, 17, 18. The regions for the solution to the Riemann problem are as follows:
\begin{itemize}
    \item Region $III_0$: The order of the contact discontinuities has now switched to $C_0$ followed by $C_a$.
    \item Region $IV$: The same behavior as the standard overcompressive regions, as in Cases 16 - 18.
    \item Region $V$: Just like in Cases 4 - 6 and Cases 16 - 18 except the $a$-shocks or $a$-rarefactions are now $a$-contact discontinuities.
\end{itemize}

Below are the state spaces for Cases 10 - 12:

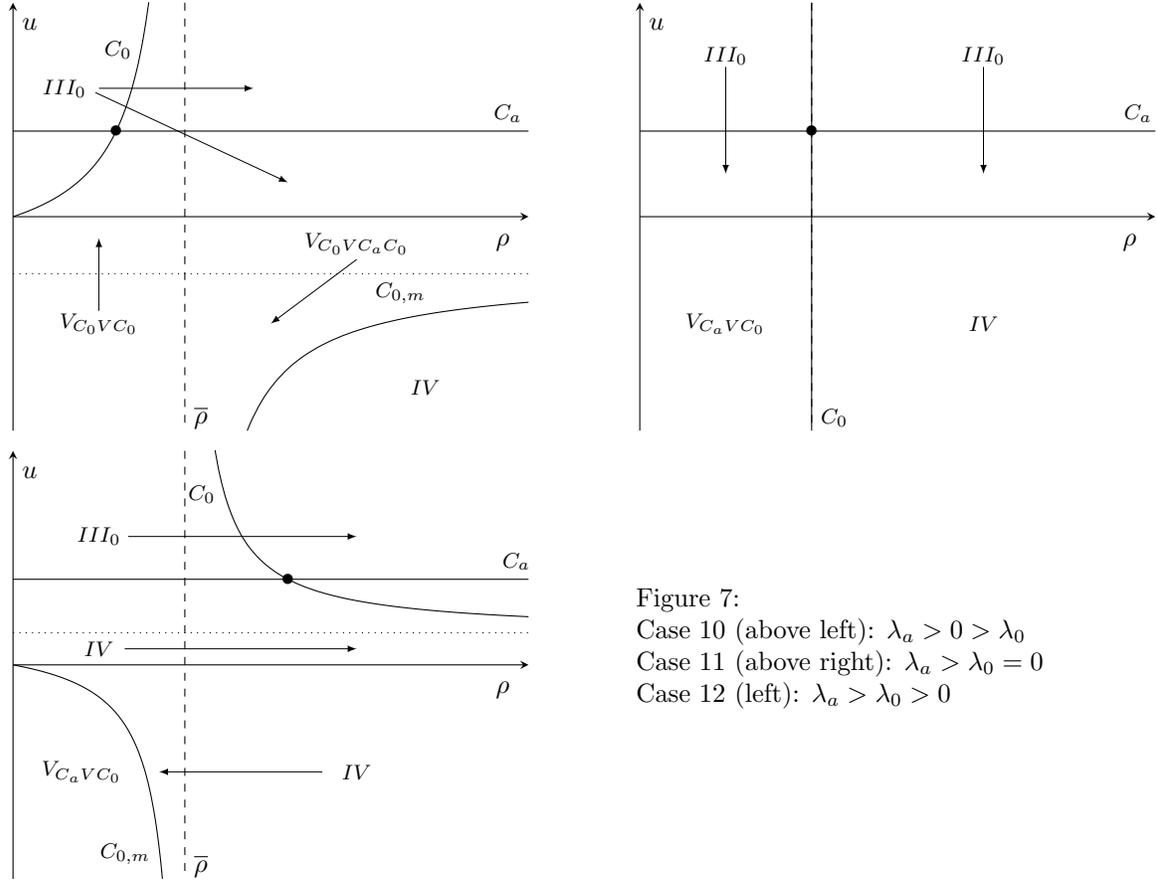
\begin{figure}[H]
    \begin{minipage}[t]{.5\textwidth}
        \vspace{0pt}
        \centering
        \begin{tikzpicture}[
        declare function={
            a = -1;       
            rhowithbar = 5; 
            rholeft = 3;    
            uleft = 4;      
            A = 0;          
            cont_disc_rholeftnotrhowithbar(\t)
                = uleft + (uleft - A)*(exp(ln(rholeft)*(a)) - exp(ln(\t)*(a)))/(exp(ln(\t)*(a)) - exp(ln(rhowithbar)*(a)));
            cont_disc_limit
                = uleft + (uleft - A)*(exp(ln(rholeft)*(a)))/(- exp(ln(rhowithbar)*(a)));
            cont_disc_limiting_curve(\t)
                = uleft + (uleft - A)*(exp(ln(1000)*(a)) - exp(ln(\t)*(a)))/(exp(ln(\t)*(a)) - exp(ln(rhowithbar)*(a)));
        }
    ]

    \tikzstyle{arrow} = [thick,->,>=stealth]

    \begin{axis}[
        axis lines = center,
        xmin=0, xmax=15,
        ymin=-10, ymax=10,
        xlabel=$\rho$,
        ylabel=$\u$,
        x label style = {at={(axis description cs:0.95,0.4)},anchor=south},
        y label style = {at={(axis description cs:0.0,0.95)},anchor=west},
        xticklabel=\empty,
        yticklabel=\empty,
        major tick style={draw=none}
    ]

    \node at (axis cs:rholeft,uleft) {\textbullet};

    \draw[dashed] (axis cs:rhowithbar,10) -- (axis cs:rhowithbar,-10)
        node[pos=0.97, right] {$\overline{\rho}$};;

    \addplot [
        domain=0.01:(rhowithbar-1.05),
        samples=100,
        color=black,
        solid
    ] {cont_disc_rholeftnotrhowithbar(x)}
        node[pos=0.8, left] {\footnotesize $C_{0}$};
    \addplot [
        domain=(rhowithbar+1.815):15,
        samples=100,
        color=black,
        solid
    ] {cont_disc_rholeftnotrhowithbar(x)}
        node[pos=0.75, above left] {\footnotesize $C_{0,m}$};
    \draw[dotted] (axis cs:0,cont_disc_limit) -- (axis cs:15,cont_disc_limit);

    \draw[solid] (axis cs:rholeft,uleft) -- (axis cs:15,uleft)
        node[pos=0.95, above] {\footnotesize $C_a$};
    \draw[solid] (axis cs:rholeft,uleft) -- (axis cs:0,uleft);

    \node at (axis cs:12,-8) {\footnotesize $IV$};
    
    \node at (axis cs:2.5,-5) {\footnotesize $V_{C_{0}VC_{0}}$};
        \draw[black,-latex,very thin] (axis cs:2.5,-4.4) -- (axis cs:2.5,-1);

    \node at (axis cs:10, -1.25) {\footnotesize $V_{C_{0}VC_{a}C_{0}}$};
        \draw[black,-latex,very thin] (axis cs:10,-2) -- (axis cs:7.5,-5);
    
    \node at (axis cs:1.5,6) {\footnotesize $III_0$};
        \draw[black,-latex,very thin] (axis cs:2.5,6) -- (axis cs:7,6);
        \draw[black,-latex,very thin] (axis cs:2.4,5.8) -- (axis cs:8,1.6);
        
    \end{axis}
\end{tikzpicture}
    \end{minipage}
    \hfill
    \begin{minipage}[t]{.5\textwidth}
        \vspace{0pt}
        \centering
        \begin{tikzpicture}[
        declare function={
            a = -1;       
            rhowithbar = 5; 
            rholeft = 5;    
            uleft = 4;      
            A = 0;          
            cont_disc_rholeftnotrhowithbar(\t)
                = uleft + (uleft - A)*(exp(ln(rholeft)*(a)) - exp(ln(\t)*(a)))/(exp(ln(\t)*(a)) - exp(ln(rhowithbar)*(a)));
            cont_disc_limit
                = uleft + (uleft - A)*(exp(ln(rholeft)*(a)))/(- exp(ln(rhowithbar)*(a)));
            cont_disc_limiting_curve(\t)
                = uleft + (uleft - A)*(exp(ln(1000)*(a)) - exp(ln(\t)*(a)))/(exp(ln(\t)*(a)) - exp(ln(rhowithbar)*(a)));
        }
    ]

    \tikzstyle{arrow} = [thick,->,>=stealth]

    \begin{axis}[
        axis lines = center,
        xmin=0, xmax=15,
        ymin=-10, ymax=10,
        xlabel=$\rho$,
        ylabel=$\u$,
        x label style = {at={(axis description cs:0.95,0.4)},anchor=south},
        y label style = {at={(axis description cs:0.0,0.95)},anchor=west},
        xticklabel=\empty,
        yticklabel=\empty,
        major tick style={draw=none}
    ]

    \node at (axis cs:rholeft,uleft) {\textbullet};

    \draw[dashed] (axis cs:rhowithbar,10) -- (axis cs:rhowithbar,-10);

    \draw[solid] (axis cs:rhowithbar,10) -- (axis cs:rhowithbar,-10)
        node[pos=0.97, right] {\footnotesize $C_0$};;

    \draw[solid] (axis cs:rholeft,uleft) -- (axis cs:15,uleft)
        node[pos=0.95, above] {\footnotesize $C_a$};
    \draw[solid] (axis cs:rholeft,uleft) -- (axis cs:0,uleft);

    \node at (axis cs:10,-5) {\footnotesize $IV$};
    
    \node at (axis cs:2.5,-5) {\footnotesize $V_{C_{a}VC_{0}}$};
    
    \node at (axis cs:2.5, 7.5) {\footnotesize $III_0$};
         \draw[black,-latex,very thin] (axis cs:2.5,7) -- (axis cs:2.5,2);
    
    \node at (axis cs:10,7.5) {\footnotesize $III_0$};
        \draw[black,-latex,very thin] (axis cs:10,7) -- (axis cs:10,2);
    \end{axis}
\end{tikzpicture}
    \end{minipage}
    \begin{minipage}[b]{.5\textwidth}
        \vspace{0pt}
        \centering
        \begin{tikzpicture}[
        declare function={
            a = -1;       
            rhowithbar = 5; 
            rholeft = 8;    
            uleft = 4;      
            A = 0;          
            cont_disc_rholeftnotrhowithbar(\t)
                = uleft + (uleft - A)*(exp(ln(rholeft)*(a)) - exp(ln(\t)*(a)))/(exp(ln(\t)*(a)) - exp(ln(rhowithbar)*(a)));
            cont_disc_limit
                = uleft + (uleft - A)*(exp(ln(rholeft)*(a)))/(- exp(ln(rhowithbar)*(a)));
            cont_disc_limiting_curve(\t)
                = uleft + (uleft - A)*(exp(ln(1000)*(a)) - exp(ln(\t)*(a)))/(exp(ln(\t)*(a)) - exp(ln(rhowithbar)*(a)));
        }
    ]
    \tikzstyle{arrow} = [thick,->,>=stealth]
    \begin{axis}[
        axis lines = center,
        xmin=0, xmax=15,
        ymin=-10, ymax=10,
        xlabel=$\rho$,
        ylabel=$\u$,
        x label style = {at={(axis description cs:0.95,0.4)},anchor=south},
        y label style = {at={(axis description cs:0.0,0.95)},anchor=west},
        xticklabel=\empty,
        yticklabel=\empty,
        major tick style={draw=none},
        view={0}{90}
    ]

    \node at (axis cs:rholeft,uleft) {\textbullet};

    \draw[dashed] (axis cs:rhowithbar,10) -- (axis cs:rhowithbar,-10)
        node[pos=0.97, right] {$\overline{\rho}$};

    \addplot [
        domain=0.01:(rhowithbar-0.58),
        samples=100,
        color=black,
        solid
    ] {cont_disc_rholeftnotrhowithbar(x)}
        node[pos=0.8, left] {\footnotesize $C_{0,m}$};
    \addplot [
        domain=(rhowithbar+0.8):15,
        samples=100,
        color=black,
        solid
    ] {cont_disc_rholeftnotrhowithbar(x)}
        node[pos=0.2, left] {\footnotesize $C_{0}$};
    \draw[dotted] (axis cs:0,cont_disc_limit) -- (axis cs:15,cont_disc_limit);

    \draw[solid] (axis cs:rholeft,uleft) -- (axis cs:15,uleft)
        node[pos=0.95, above] {\footnotesize $C_a$};
    \draw[solid] (axis cs:rholeft,uleft) -- (axis cs:0,uleft);

    \node at (axis cs:2,-5) {\footnotesize $V_{C_{a}VC_{0}}$};
    
    \node at (axis cs:2.5,0.75) {\footnotesize $IV$};
        \draw[black,-latex,very thin] (axis cs:3.25,0.75) -- (axis cs:10,0.75);
    
    \node at (axis cs:10,-5) {\footnotesize $IV$};
        \draw[black,-latex,very thin] (axis cs:9,-5) -- (axis cs:4.25,-5);

    \node at (axis cs:2.5,6) {\footnotesize $III_0$};
        \draw[black,-latex,very thin] (axis cs:3.35,6) -- (axis cs:10,6);
    
    \end{axis}
\end{tikzpicture}
    \end{minipage}
    \hfill
    \begin{minipage}[b]{.45\textwidth}
        \caption{\\
            Case 10 (above left): $\lambda_a > 0 > \lambda_0$ \\
            Case 11 (above right): $\lambda_a > \lambda_0 = 0$ \\
            Case 12 (left): $\lambda_a > \lambda_0 > 0$
        }
        \label{fig:cases10-12}
        \vspace{2cm}
    \end{minipage}
\end{figure}

\subsubsection{\texorpdfstring{$a>0$}{a>0}}
The last six cases are when $a > 0$. 
\begin{itemize}
    \item Case 19: $u_L<0$ and $\rho_L<\rhobar$,
    \item Case 20: $u_L<0$ and $\rho_L=\rhobar$, and
    \item Case 21: $u_L<0$ and $\rho_L>\rhobar$.
\end{itemize}
The various regions for the solution to the Riemann problem are
\begin{itemize}
    \item Region $I_0$: The left state is connected to the right by a $0$-contact discontinuity followed by an $a$-shock.
    \item Region $II_0$: The left state is connected to the right by a $0$-contact discontinuity followed by an $a$-rarefaction.
    \item Region $IV$: The left state is connected to the right by an overcompressive delta shock.
    \item Region $V$: The vacuum states can occur in multiple ways, first either following $C_0$ or $S_a$, and then staying on the line $\rho = 0$ for some nonzero time. When the solution can, it uses a $0$-contact discontinuity, and otherwise an $a$-shock or $a$-rarefaction.
\end{itemize}
The following are the state spaces for Cases 19 - 21:

\begin{figure}[H]
    \!\begin{minipage}[t]{.5\textwidth}
        \vspace{0pt}
        \centering
        \begin{tikzpicture}[
        declare function={
            a = .5;       
            rhowithbar = 5; 
            rholeft = 3;    
            uleft = -4;      
            A = 0;          
            cont_disc_rholeftnotrhowithbar(\t)
                = uleft + (uleft - A)*(exp(ln(rholeft)*(a)) - exp(ln(\t)*(a)))/(exp(ln(\t)*(a)) - exp(ln(rhowithbar)*(a)));
            cont_disc_limit
                = uleft + (uleft - A)*(exp(ln(rholeft)*(a)))/(- exp(ln(rhowithbar)*(a)));
            cont_disc_limiting_curve(\t)
                = uleft + (uleft - A)*(exp(ln(1000)*(a)) - exp(ln(\t)*(a)))/(exp(ln(\t)*(a)) - exp(ln(rhowithbar)*(a)));
            rho1 = (rhowithbar)/(exp((1/a)*ln(a+1)));
        }
    ]
    \tikzstyle{arrow} = [thick,->,>=stealth]
    \begin{axis}[
        axis lines = center,
        xmin=0, xmax=15,
        ymin=-10, ymax=10,
        xlabel=$\rho$,
        ylabel=$\u$,
        x label style = {at={(axis description cs:1.02,0.45)},anchor=south},
        y label style = {at={(axis description cs:0.0,0.95)},anchor=west},
        xticklabel=\empty,
        yticklabel=\empty,
        major tick style={draw=none},
        view={0}{90}
    ]

    \node at (axis cs:rholeft,uleft) {\textbullet};

    \draw[dashed] (axis cs:rhowithbar,10) -- (axis cs:rhowithbar,-10)
        node[pos=0.97, right] {$\overline{\rho}$};

    \addplot [
        domain=0.01:(rhowithbar-0.58),
        samples=100,
        color=black,
        solid
    ] {cont_disc_rholeftnotrhowithbar(x)}
        node[pos=0.6, left] {\footnotesize $C_{0}$};
    \addplot [
        domain=(rhowithbar+0.8):15,
        samples=100,
        color=black,
        solid
    ] {cont_disc_rholeftnotrhowithbar(x)}
        node[pos=0.15, right] {\footnotesize $C_{0,m}$};
    \draw[dotted] (axis cs:0,cont_disc_limit) -- (axis cs:15,cont_disc_limit);

    \draw[solid] (axis cs:rholeft,uleft) -- (axis cs:15,uleft)
        node[pos=0.95, above] {\footnotesize $R_a$};
    \draw[solid] (axis cs:rholeft,uleft) -- (axis cs:0,uleft)
        node[pos=0.85, above] {\footnotesize $S_a$};

    \node at (axis cs:11,6) {\footnotesize $IV$};
    
    \node at (axis cs:1.75,-6) {\footnotesize $I_0$};
    
    \node at (axis cs:10,-2) {\footnotesize $II_0$};
        \draw[black,-latex,very thin] (axis cs:9.5,-2) -- (axis cs:3,-2);
        \draw[black,-latex,very thin] (axis cs:9.5,-2.3) -- (axis cs:4.5,-7.5);
        
    \node at (axis cs:2.5,5) {\footnotesize $V_{C_{0}VC_{0}}$};

    \node at (axis cs:7.5,2) {\footnotesize $V_{C_{0}VS_{a}C_{0}}$};
    
    \node at (axis cs:5,-.5) {\footnotesize $V_{C_{0}VR_{a}}$};
    \end{axis}
\end{tikzpicture}
    \end{minipage}
    \hfill
    \begin{minipage}[t]{.5\textwidth}
        \vspace{0pt}
        \centering
        \begin{tikzpicture}[
        declare function={
            a = .5;       
            rhowithbar = 5; 
            rholeft = 5;    
            uleft = -4;      
            A = 0;          
            cont_disc_rholeftnotrhowithbar(\t)
                = uleft + (uleft - A)*(exp(ln(rholeft)*(a)) - exp(ln(\t)*(a)))/(exp(ln(\t)*(a)) - exp(ln(rhowithbar)*(a)));
            cont_disc_limit
                = uleft + (uleft - A)*(exp(ln(rholeft)*(a)))/(- exp(ln(rhowithbar)*(a)));
            cont_disc_limiting_curve(\t)
                = uleft + (uleft - A)*(exp(ln(1000)*(a)) - exp(ln(\t)*(a)))/(exp(ln(\t)*(a)) - exp(ln(rhowithbar)*(a)));
        }
    ]

    \tikzstyle{arrow} = [thick,->,>=stealth]

    \begin{axis}[
        axis lines = center,
        xmin=0, xmax=15,
        ymin=-10, ymax=10,
        xlabel=$\rho$,
        ylabel=$\u$,
        x label style = {at={(axis description cs:1.02,0.45)},anchor=south},
        y label style = {at={(axis description cs:0.0,0.95)},anchor=west},
        xticklabel=\empty,
        yticklabel=\empty,
        major tick style={draw=none}
    ]

    \node at (axis cs:rholeft,uleft) {\textbullet};

    \draw[dashed] (axis cs:rhowithbar,10) -- (axis cs:rhowithbar,-10);

    \draw[solid] (axis cs:rhowithbar,10) -- (axis cs:rhowithbar,-10)
        node[pos=0.97, right] {\footnotesize $C_0$};;

    \draw[solid] (axis cs:rholeft,uleft) -- (axis cs:15,uleft)
        node[pos=0.95, above] {\footnotesize $R_a$};
    \draw[solid] (axis cs:rholeft,uleft) -- (axis cs:0,uleft)
        node[pos=0.85, above] {\footnotesize $S_a$};

    \node at (axis cs:10,5) {\footnotesize $IV$};
    
    \node at (axis cs:2.5,5) {\footnotesize $V_{S_{a}VC_{0}}$};
    
    \node at (axis cs:2.5,-2) {\footnotesize $I_0$};
        \draw[black,-latex,very thin] (axis cs:2.5,-2.5) -- (axis cs:2.5,-7);
    
    \node at (axis cs:10,-2) {\footnotesize $II_0$};
        \draw[black,-latex,very thin] (axis cs:10,-2.5) -- (axis cs:10,-7);
    \end{axis}
\end{tikzpicture}
    \end{minipage}
    \begin{minipage}[b]{.49\textwidth}
        \vspace{0pt}
        \centering
        \begin{tikzpicture}[
        declare function={
            a = .5;       
            rhowithbar = 5; 
            rholeft = 8;    
            uleft = -4;      
            A = 0;          
            cont_disc_rholeftnotrhowithbar(\t)
                = uleft + (uleft - A)*(exp(ln(rholeft)*(a)) - exp(ln(\t)*(a)))/(exp(ln(\t)*(a)) - exp(ln(rhowithbar)*(a)));
            cont_disc_limit
                = uleft + (uleft - A)*(exp(ln(rholeft)*(a)))/(- exp(ln(rhowithbar)*(a)));
            rho1 = (rhowithbar)/(exp((1/a)*ln(a+1)));
            rho1_bound = cont_disc_rholeftnotrhowithbar(rho1);
        }
    ]

    \tikzstyle{arrow} = [thick,->,>=stealth]

    \begin{axis}[
        axis lines = center,
        xmin=0, xmax=15,
        ymin=-10, ymax=10,
        xlabel=$\rho$,
        ylabel=$\u$,
        x label style = {at={(axis description cs:1.02,0.45)},anchor=south},
        y label style = {at={(axis description cs:0.0,0.95)},anchor=west},
        xticklabel=\empty,
        yticklabel=\empty,
        major tick style={draw=none}
    ]

    \node at (axis cs:rholeft,uleft) {\textbullet};

    \draw[dashed] (axis cs:rhowithbar,10) -- (axis cs:rhowithbar,-10)
        node[pos=0.97, right] {$\overline{\rho}$};

    \addplot [
        domain=0.01:(rhowithbar-0.58),
        samples=100,
        color=black,
        solid
    ] {cont_disc_rholeftnotrhowithbar(x)}
        node[pos=0.5, left] {\footnotesize $C_{0,m}$};
    \addplot [
        domain=(rhowithbar+0.8):15,
        samples=100,
        color=black,
        solid
    ] {cont_disc_rholeftnotrhowithbar(x)}
        node[pos=0.3, right] {\footnotesize $C_{0}$};
    \draw[dotted] (axis cs:0.01,cont_disc_limit) -- (axis cs:15,cont_disc_limit);

    \draw[solid] (axis cs:rholeft,uleft) -- (axis cs:15,uleft)
        node[pos=0.95, above] {\footnotesize $R_a$};
    \draw[solid] (axis cs:rholeft,uleft) -- (axis cs:0,uleft)
        node[pos=0.85, above] {\footnotesize $S_a$};

    \node at (axis cs:1.5,6) {\footnotesize $V_{S_{a}VC_{0}}$};
    
    \node at (axis cs:10,5) {\footnotesize $IV$};
        \draw[black,-latex,very thin] (axis cs:9.35,4.7) -- (axis cs:2.5,0.5);
        \draw[black,-latex,very thin] (axis cs:10,4.4) -- (axis cs:10,0.5);
        
    \node at (axis cs:11,-7) {\footnotesize $II_0$};
    
    \node at (axis cs:2.5,-2) {\footnotesize $I_0$};
        \draw[black,-latex,very thin] (axis cs:2.5,-2.5) -- (axis cs:2.5,-7.5);
        \draw[black,-latex,very thin] (axis cs:2.85,-2.25) -- (axis cs:6,-5);
        \draw[black,-latex,very thin] (axis cs:2.95,-2) -- (axis cs:8,-2);

    \end{axis}
\end{tikzpicture}
    \end{minipage}
    \hfill
    \begin{minipage}[b]{.49\textwidth}
        \caption{\\
            Case 19: $\lambda_a>\lambda_0$ and $0>\lambda_0$ \\
            Case 20: $\lambda_a>\lambda_0=0$ \\
            Case 21: $\lambda_a>\lambda_0>0$
        }
        \label{fig:cases19-21}
        \vspace{2cm}
    \end{minipage}
\end{figure}
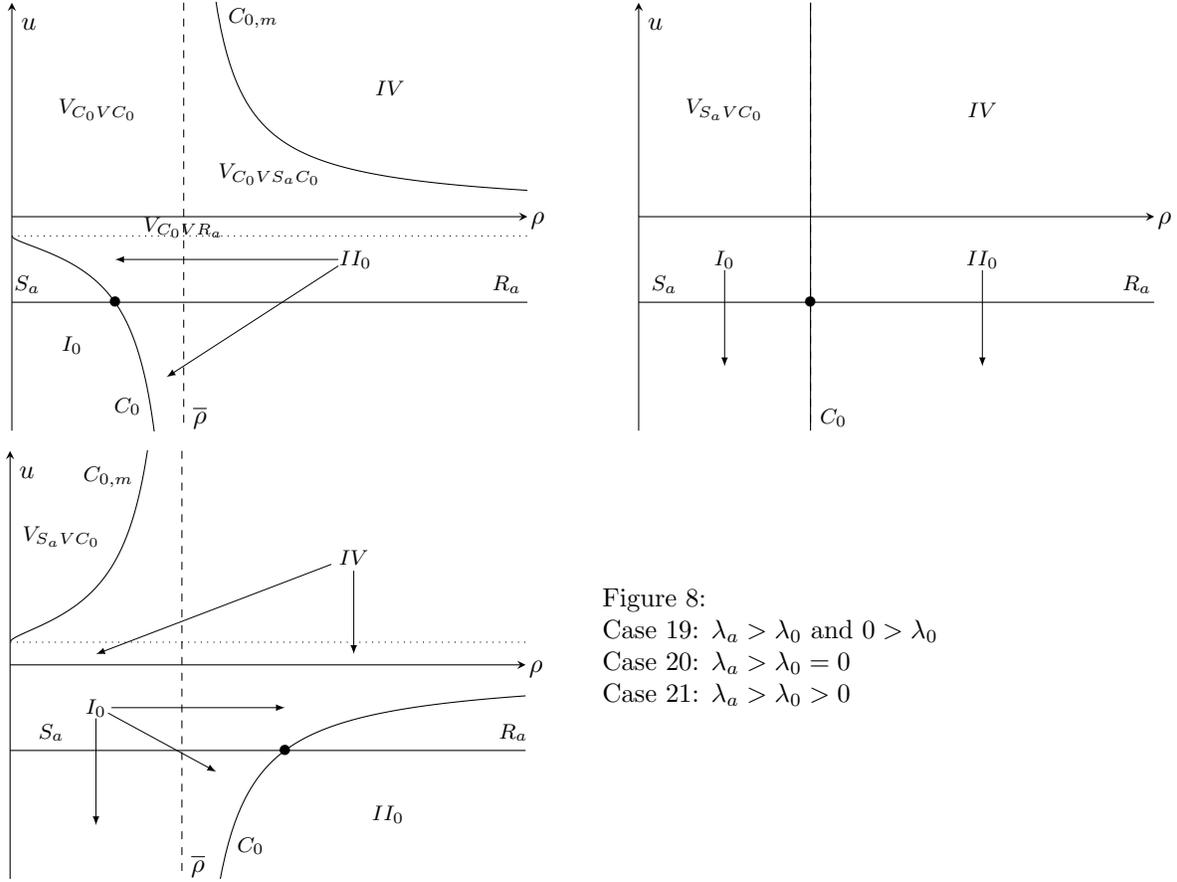

Cases 22 through 24 are similar to the preceding three except that $\u_L>0$:
\begin{itemize}
    \item Case 22: $\u_L>0$ and $\rho_L<\rhobar$,
    \item Case 23: $\u_L>0$ and $\rho_L=\rhobar$, and
    \item Case 24: $\u_L>0$ and $\rho_L>\rhobar$.
\end{itemize}
The regions which exhibit a solution to the Riemann problem are as follows:
\begin{itemize}
    \item Region $I_a$: The right state is reached by an $a$-shock followed by a $0$-contact discontinuity;
    \item Region $II_a$: The right state is reached by an $a$-rarefaction followed by a $0$-contact discontinuity;
    \item Region $V$: The left follows an $a$-rarefaction to reach a vacuum state, then a $0$-contact discontinuity from the vacuum state to reach the right state.
    \item Region $VI$: The right state is reached by an $a$-shock from Case 22 and an $a$-rarefaction from Case 24, both into Case 23. It then follows a $0$-contact discontinuity into Case 20.
\end{itemize}

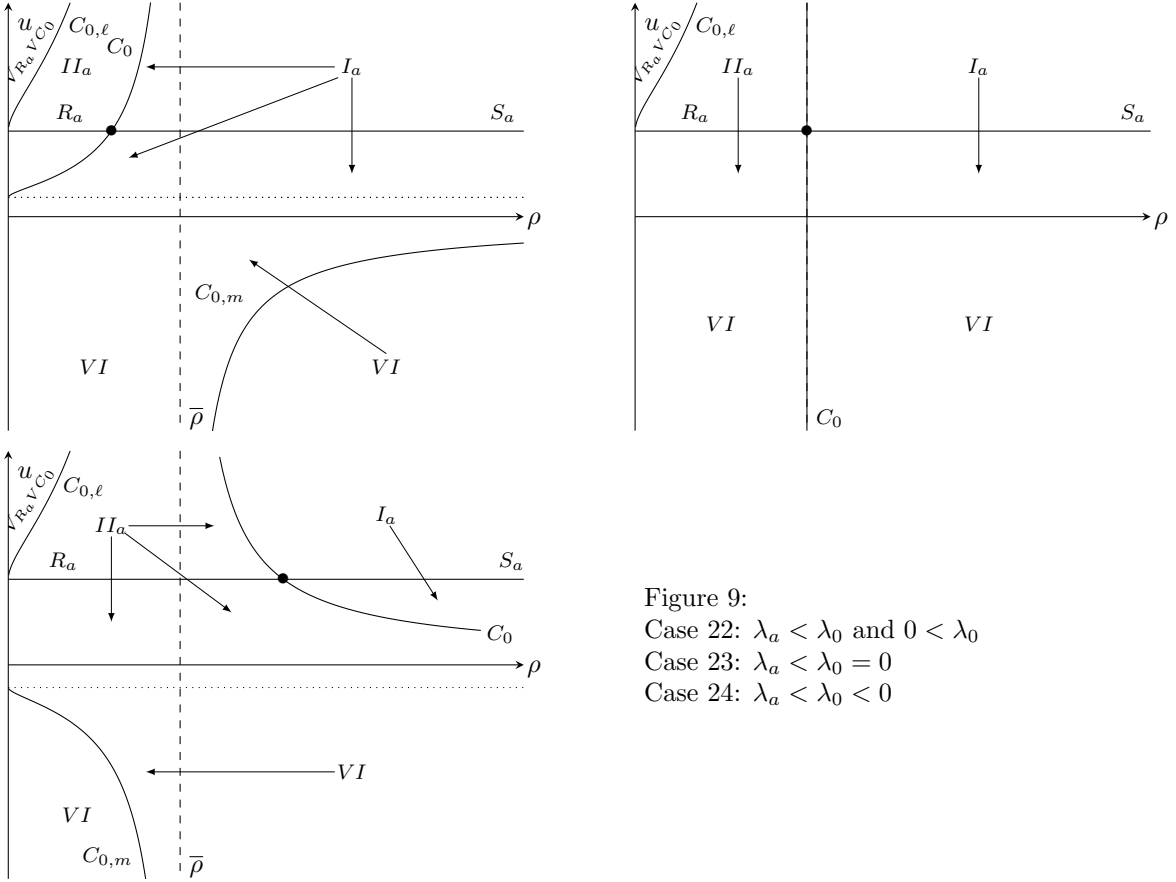
\begin{figure}[H]
    \begin{minipage}[t]{.5\textwidth}
        \vspace{0pt}
        \centering
        \begin{tikzpicture}[
        declare function={
            a = .5;       
            rhowithbar = 5; 
            rholeft = 3;    
            uleft = 4;      
            A = 0;          
            cont_disc_rholeftnotrhowithbar(\t)
                = uleft + (uleft - A)*(exp(ln(rholeft)*(a)) - exp(ln(\t)*(a)))/(exp(ln(\t)*(a)) - exp(ln(rhowithbar)*(a)));
            cont_disc_limit
                = uleft + (uleft - A)*(exp(ln(rholeft)*(a)))/(- exp(ln(rhowithbar)*(a)));
            cont_disc_limiting_curve(\t)
                = uleft + (uleft - A)*(0 - exp(ln(\t)*(a)))/(exp(ln(\t)*(a)) - exp(ln(rhowithbar)*(a)));
        }
    ]

    \tikzstyle{arrow} = [thick,->,>=stealth]

    \begin{axis}[
        axis lines = center,
        xmin=0, xmax=15,
        ymin=-10, ymax=10,
        xlabel=$\rho$,
        ylabel=$\u$,
        x label style = {at={(axis description cs:1.02,0.45)},anchor=south},
        y label style = {at={(axis description cs:0.0,0.95)},anchor=west},
        xticklabel=\empty,
        yticklabel=\empty,
        major tick style={draw=none}
    ]

    \node at (axis cs:rholeft,uleft) {\textbullet};

    \draw[dashed] (axis cs:rhowithbar,10) -- (axis cs:rhowithbar,-10)
        node[pos=.97, right] {$\overline{\rho}$};;

    \addplot [
        domain=0.01:(rhowithbar-.86),
        samples=100,
        color=black,
        solid
    ] {cont_disc_rholeftnotrhowithbar(x)}
        node[pos=0.8, left] {\footnotesize $C_{0}$};
    \addplot [
        domain=(rhowithbar+.7):15,
        samples=100,
        color=black,
        solid
    ] {cont_disc_rholeftnotrhowithbar(x)}
        node[pos=0.5, above left] {\footnotesize $C_{0,m}$};

    \addplot [
        domain=0.01:(rhowithbar-3.2),
        samples=100,
        color=black,
        solid
   ] {cont_disc_limiting_curve(x)}
    node[pos=0.8, right] {\footnotesize $C_{0,\ell}$};

    \draw[dotted] (axis cs:0,cont_disc_limit) -- (axis cs:15,cont_disc_limit);
    \draw[dotted] (axis cs:0,cont_disc_limit) -- (axis cs:15,cont_disc_limit);

    \draw[solid] (axis cs:rholeft,uleft) -- (axis cs:15,uleft)
        node[pos=0.95, above] {\footnotesize $S_a$};
    \draw[solid] (axis cs:rholeft,uleft) -- (axis cs:0,uleft)
        node[pos=.4, above] {\footnotesize $R_a$};

    \node[rotate=60, anchor=center] at (axis cs:.6,7.75) {\scriptsize $V_{R_aVC_{0}}$};
    
    \node at (axis cs:2,7) {\footnotesize $II_a$};
    
    \node at (axis cs:10,7) {\footnotesize $I_a$};
        \draw[black,-latex,very thin] (axis cs:10,6.5) -- (axis cs:10,2);
        \draw[black,-latex,very thin] (axis cs:9.5,7) -- (axis cs:4,7);
        \draw[black,-latex,very thin] (axis cs:9.6,6.5) -- (axis cs:3.5,2.75);
        
    \node at (axis cs:2.5,-7) {\footnotesize $VI$};
    
    \node at (axis cs:11,-7) {\footnotesize $VI$};
    \draw[black,-latex,very thin] (axis cs:11,-6.4) -- (axis cs:7,-2);
    
    \end{axis}

\end{tikzpicture}
    \end{minipage}
    \hfill
    \begin{minipage}[t]{.5\textwidth}
        \vspace{0pt}
        \centering
        \begin{tikzpicture}[
        declare function={
            a = .5;       
            rhowithbar = 5; 
            rholeft = 5;    
            uleft = 4;      
            A = 0;          
            cont_disc_rholeftnotrhowithbar(\t)
                = uleft + (uleft - A)*(exp(ln(rholeft)*(a)) - exp(ln(\t)*(a)))/(exp(ln(\t)*(a)) - exp(ln(rhowithbar)*(a)));
            cont_disc_limit
                = uleft + (uleft - A)*(exp(ln(rholeft)*(a)))/(- exp(ln(rhowithbar)*(a)));
            cont_disc_limiting_curve(\t)
                = uleft + (uleft - A)*(0 - exp(ln(\t)*(a)))/(exp(ln(\t)*(a)) - exp(ln(rhowithbar)*(a)));
        }
    ]

    \tikzstyle{arrow} = [thick,->,>=stealth]

    \begin{axis}[
        axis lines = center,
        xmin=0, xmax=15,
        ymin=-10, ymax=10,
        xlabel=$\rho$,
        ylabel=$\u$,
        x label style = {at={(axis description cs:1.02,0.45)},anchor=south},
        y label style = {at={(axis description cs:0.0,0.95)},anchor=west},
        xticklabel=\empty,
        yticklabel=\empty,
        major tick style={draw=none}
    ]

    \node at (axis cs:rholeft,uleft) {\textbullet};

    \draw[dashed] (axis cs:rhowithbar,10) -- (axis cs:rhowithbar,-10);

    \addplot [
        domain=0.01:(rhowithbar-3.2),
        samples=100,
        color=black,
        solid
   ] {cont_disc_limiting_curve(x)}
    node[pos=0.8, right] {\footnotesize $C_{0,\ell}$};
    
    \draw[solid] (axis cs:rhowithbar,10) -- (axis cs:rhowithbar,-10)
        node[pos=0.97, right] {\footnotesize $C_0$};;

    \draw[solid] (axis cs:rholeft,uleft) -- (axis cs:15,uleft)
        node[pos=0.95, above] {\footnotesize $S_a$};
    \draw[solid] (axis cs:rholeft,uleft) -- (axis cs:0,uleft)
        node[pos=0.65, above] {\footnotesize $R_a$};

    \node at (axis cs:10,7) {\footnotesize $I_a$};
        \draw[black,-latex,very thin] (axis cs:10,6.5) -- (axis cs:10,2);
    
    \node[rotate=60, anchor=center] at (axis cs:.6,7.75) {\scriptsize $V_{R_aVC_0}$};
    
    \node at (axis cs:3,7) {\footnotesize $II_a$};
        \draw[black,-latex,very thin] (axis cs:3,6.5) -- (axis cs:3,2);
    
    \node at (axis cs:2.5,-5) {\footnotesize $VI$};
    
    \node at (axis cs:10,-5) {\footnotesize $VI$};
    \end{axis}
\end{tikzpicture}
    \end{minipage}
    \begin{minipage}[b]{.5\textwidth}
        \vspace{0pt}
        \centering
        \begin{tikzpicture}[
        declare function={
            a = .5;       
            rhowithbar = 5; 
            rholeft = 8;    
            uleft = 4;      
            A = 0;          
            cont_disc_rholeftnotrhowithbar(\t)
                = uleft + (uleft - A)*(exp(ln(rholeft)*(a)) - exp(ln(\t)*(a)))/(exp(ln(\t)*(a)) - exp(ln(rhowithbar)*(a)));
            cont_disc_limit
                = uleft + (uleft - A)*(exp(ln(rholeft)*(a)))/(- exp(ln(rhowithbar)*(a)));
            cont_disc_limiting_curve(\t)
                = uleft + (uleft - A)*(0 - exp(ln(\t)*(a)))/(exp(ln(\t)*(a)) - exp(ln(rhowithbar)*(a)));
        }
    ]

    \tikzstyle{arrow} = [thick,->,>=stealth]

    \begin{axis}[
        axis lines = center,
        xmin=0, xmax=15,
        ymin=-10, ymax=10,
        xlabel=$\rho$,
        ylabel=$\u$,
        x label style = {at={(axis description cs:1.02,0.45)},anchor=south},
        y label style = {at={(axis description cs:0.0,0.95)},anchor=west},
        xticklabel=\empty,
        yticklabel=\empty,
        major tick style={draw=none}
    ]

    \node at (axis cs:rholeft,uleft) {\textbullet};

    \draw[dashed] (axis cs:rhowithbar,10) -- (axis cs:rhowithbar,-10)
        node[pos=0.97, right] {$\overline{\rho}$};

    \addplot [
        domain=0.01:(rhowithbar-1),
        samples=100,
        color=black,
        solid
    ] {cont_disc_rholeftnotrhowithbar(x)}
        node[pos=.9, left] {\footnotesize $C_{0,m}$};
    \addplot [
        domain=(rhowithbar+1.15):13.75,
        samples=100,
        color=black,
        solid
    ] {cont_disc_rholeftnotrhowithbar(x)}
        node[pos=1.1, left] {\footnotesize $C_{0}$};
    \addplot [
        domain=0.01:(rhowithbar-3.2),
        samples=100,
        color=black,
        solid
    ] {cont_disc_limiting_curve(x)}
        node[pos=0.7, right] {\footnotesize $C_{0,\ell}$};
    \draw[dotted] (axis cs:0,cont_disc_limit) -- (axis cs:15,cont_disc_limit);

    \draw[solid] (axis cs:rholeft,uleft) -- (axis cs:15,uleft)
        node[pos=0.95, above] {\footnotesize $S_a$};
    \draw[solid] (axis cs:rholeft,uleft) -- (axis cs:0,uleft)
        node[pos=0.8, above] {\footnotesize $R_a$};

    \node at (axis cs:11,7) {\footnotesize $I_a$};
        \draw[black,-latex,very thin] (axis cs:11.1,6.5) -- (axis cs:12.5,3);
    
    \node[rotate=60, anchor=center] at (axis cs:.6,7.75) {\scriptsize $V_{R_aVC_0}$};
    
    \node at (axis cs:3,6.5) {\footnotesize $II_a$};
        \draw[black,-latex,very thin] (axis cs:3,6) -- (axis cs:3,2);
        \draw[black,-latex,very thin] (axis cs:3.5,6.5) -- (axis cs:6,6.5);
        \draw[black,-latex,very thin] (axis cs:3.4,6.15) -- (axis cs:6.5,2.45);
        
    \node at (axis cs:2,-7) {\footnotesize $VI$};
    
    \node at (axis cs:10,-5) {\footnotesize $VI$};
        \draw[black,-latex,very thin] (axis cs:9.5,-5) -- (axis cs:4,-5);
    \end{axis}
\end{tikzpicture}
    \end{minipage}
    \hfill
    \begin{minipage}[b]{.45\textwidth}
        \caption{\\
            Case 22: $\lambda_a<\lambda_0$ and $0<\lambda_0$ \\
            Case 23: $\lambda_a<\lambda_0=0$ \\
            Case 24: $\lambda_a < \lambda_0< 0$
        }
        \label{fig:cases22-24}
        \vspace{2cm}
    \end{minipage}
\end{figure}

\section{Geometric Singular Perturbation Theory} \label{GSPT}
In this section, we use geometric singular perturbation theory (GSPT), a central tool in the study of singularly perturbed systems, to establish the existence of self-similar viscous profiles for sufficiently small $\e > 0$, thus confirming the formation of overcompressive delta shocks we have described in the previous sections. In our case, the solution consists of two distinct parts: an outer region, comprising the constant states $(\rho_L, u_L)$ and $(\rho_R, u_R)$, and an inner region that captures the transition layer on a rescaled fast time scale.

Normal hyperbolicity ensures that directions transverse to an invariant manifold dominate the dynamics and allows us to apply Fenichel’s theory \cite{Fe}. However, in the vicinity of specific critical points, this condition fails, complicating the analysis. Following Schecter's work \cite{Sc}, we introduced a blow-up transformation that regularizes the dynamics near the degenerate set, revealing hidden hyperbolic behavior. Our analysis begins by studying the system in the singular limit $\e = 0$, then uses the Exchange Lemma to guarantee the persistence of connecting orbits in the perturbed system when $\e > 0$ is small. By verifying the required transversality and non-degeneracy conditions, we show that the singular configuration, a concatenation of outer and inner layers, extends to a smooth orbit in the full system. This construction provides a rigorous justification for the appearance of overcompressive delta shocks as singular limits of smooth viscous profiles. 

Following \cite{Ke_4, Sc}, we employ the Dafermos regularization, which preserves self-similar solutions to the traveling wave of the system. Specifically, we analyze the resulting equations in self-similar coordinates to study the existence and structure of viscous profiles. The regularized equations take the form
\begin{numcases}{}
    $$\rho_t + \left[\rho \u - \rho \u \fraction{}^a\right]_x = \e t \rho_{xx},$$ \\
    $$(\rho \u)_t + \left[\rho \u^2 - \rho \u^2 \fraction{}^a\right]_x = \e t (\rho \u)_{xx}.$$
\end{numcases}

Let $\xi = \frac{x}{t}$. Then, for any self-similar function $U(x,t) = U(\xi)$, we have $U_x = U_\xi \dfrac{\xi}{x} = \frac{1}{t}U_\xi$, $U_{xx} = \dfrac{}{x}\bigl(U_x\bigr)\dfrac{\xi}{x} = \frac{1}{t^2}U_{\xi \xi}$, and $U_t = U_{\xi} \dfrac{\xi}{t} = -\frac{x}{t^2} U_\xi$. Substituting these into the above regularized system and multiplying by $t$, we have
\begin{align}
    \begin{cases}
        -\xi \rho_\xi + \left[\rho \u - \rho \u \fraction{}^a\right]_\xi = \e \rho_{\xi \xi}, \\
        -\xi (\rho \u)_\xi + \left[\rho \u^2 - \rho \u^2 \fraction{}^a\right]_\xi = \e (\rho \u)_{\xi \xi}.
    \end{cases}
\end{align}

The initial conditions from the Riemann problem translate, in the self-similar framework, to the following far-field conditions:
\begin{align}
    \begin{cases}
        (\rho, \u)(-\infty) = (\rho_L, \u_L), \\
        (\rho, \u)(+\infty) = (\rho_R, \u_R).
    \end{cases}
\end{align}
From this point on, we proceed to work in vector notation and use $m = \rho \u$. Let $U = (\rho \enspace m)^T$ and note that we now have
$$F(U) = \begin{pmatrix} F_1(U) \\ F_2(U) \end{pmatrix} = \begin{pmatrix}
    m - m \fraction{}^a \\
    \frac{m^2}{\rho} - \frac{m^2}{\rho} \fraction{}^a
\end{pmatrix}.$$
For convenience, we introduce $V = \e \left(\dfrac{\rho}{\xi} \enspace \dfrac{m}{\xi} \right)^T$ and the slow time $\vartheta = \xi - s_{\text{singular}}$, treating $\xi$ as an additional state variable. This increases the dimension but yields an autonomous system. Thus, we obtain, letting $'$ represent differentiation with respect to $\vartheta$,
\begin{align}
    \begin{cases}
        \e U' = V, \\
        \e V' = \bigl(DF(U) - \xi I\bigr)V, \\
        \xi' = 1,
    \end{cases}
\end{align}
where $DF$ is the Fréchet derivative of $F$ with respect to the variables $\rho$ and $m$.

As this is singular when $\e \rightarrow 0$ (higher-order derivatives vanish and the system drops to lower order), we introduce $\tau$, the fast time such that $\vartheta = \e \tau$, and let $\dot{}$ represent differentiation with respect to $\tau$. Then we have
\begin{align}
    \begin{cases}
        \dot{U} = V, \\
        \dot{V} = \bigl(DF - \xi I\bigr) V, \\
        \dot{\xi} = \e,
    \end{cases} \label{eq:fast_time_sys_eneq0}
\end{align}
subject to the boundary conditions
\begin{align}
    \begin{cases}
        (U,V,\xi)(-\infty) = (U_L, 0 , -\infty), \\
        (U,V,\xi)(+\infty) = (U_L, 0 , +\infty).
    \end{cases}
\end{align}

At $\e = 0$, this becomes
\begin{align}
    \begin{cases}
        \dot{U_1} = V_1, \\
        \dot{U_2} = V_2, \\
        \dot{V_1} = \left(-\frac{\, a m \rho^{a-1} \,}{\rhobar^a} - \xi\right) V_1 + \left(1 - \fraction{}^a\right) V_2, \\
        \dot{V_2} = \left(-\frac{m^2}{\rho^2} \left(1 - \fraction{}^a\right) - \frac{a m^2}{\rho} \frac{\rho^{a-1}}{\rhobar^a}\right) V_1 + \left(\frac{2m}{\rho}\left(1 - \fraction{}^a\right) - \xi\right) V_2, \\
        \dot{\xi} = 0.
    \end{cases} \label{eq:fast_time_sys}
\end{align}
Thus, the 3-dimensional space $S = \{(U, V, \xi) : V = 0\}$ consists of equilibria when $\e = 0$ and is an invariant subspace under \eqref{eq:fast_time_sys} for every $\e$. In fact, for $\e>0$, the dynamics in this space reduce to linear motion in the increasing $\xi$ direction.

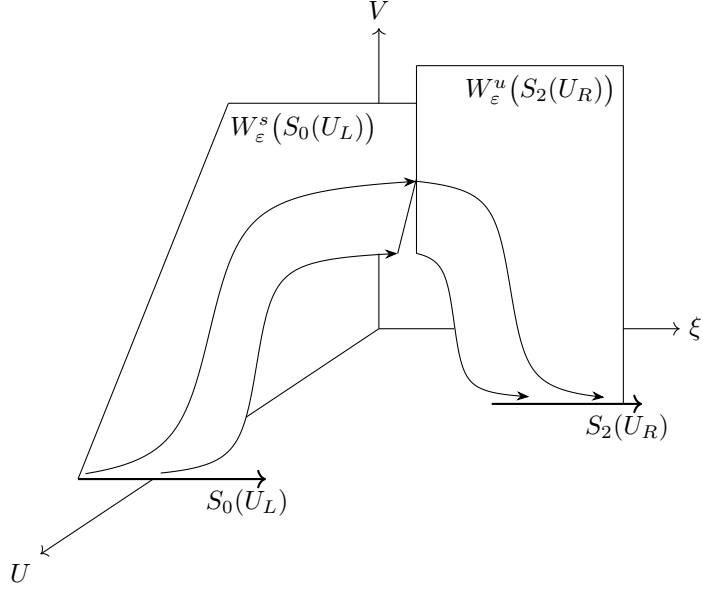
\begin{figure}[h]
    \centering
    \begin{tikzpicture}[scale=2]
        \draw (0,0)--(.5,0);
        \draw[->] (1.625,0)--(2,0) node[right] {$\xi$};
        \draw (0,0)--(0,.5);
        \draw[->] (0,1.5)--(0,2) node[above] {$V$};
        \draw (0,0)--(-.88215,-.5881);
        \draw[->] (-1.5,-1)--(-2.25,-1.5) node[below left] {$U$};

        \draw[->, thick] (-2,-1)--(-.75,-1) node[pos=.9, below] {$S_0(U_L)$};
        \draw[->, thick] (.75,-.5)--(1.75,-.5) node[pos=.9, below] {$S_2(U_R)$};

        \draw (1.625,-.5)--(1.625,1.75);
        \draw (1.625,1.75)--(.25,1.75) node[pos=.4, below] {$W_\e^u\bigl(S_2(U_R)\bigr)$};
        \draw (.25,1.75)--(.25,.5);
    
        \draw (.25,1.5)--(-1,1.5) node[pos=.6, below] {$W_\e^s\bigl(S_0(U_L)\bigr)$};
        \draw (-1,1.5)--(-2,-1);
        \draw (.125,.5)--(.25,1);

        \draw[-Stealth, samples=100] plot[domain=.25:1] ({\x}, {-(pi/180)*atan(64*(\x-.5)/3)/(3) + .04});
        \draw[-Stealth, samples=100] plot[domain=.25:1.5] ({\x}, {-(pi/180)*atan(12*(\x-.875))/(2) + .26208});

        \draw[-Stealth, samples=100] plot[domain=-1.45:.125] ({\x}, {(pi/180)*atan(12*(\x+1.625/2))/2 - .24});
        \draw[-Stealth, samples=100] plot[domain=-1.95:.25] ({\x}, {(pi/180)*atan(6*(\x+2.25/2))/1.45 - 0.019});
    \end{tikzpicture}
    \caption{The objective of finding a heteroclinic orbit connecting the left state to the right state.}
    \label{fig:inner_soln_at_intersection}
\end{figure}

To analyze the structure near the equilibria, we linearize the full system around the fixed points $V=0$, $\e = 0$. The Jacobian matrix at these points reveals eigenvalues $\lambda = 0$ with multiplicity three, corresponding to free directions in the components $U$ and $\xi$. The remaining eigenvalues are real, nonzero when $a\neq 0$ and are given by $\lambda = \lambda_a - \xi$ and $\lambda = \lambda_0 - \xi$.  Therefore, the linearization admits a full set of eigenvectors, ensuring the existence of a complete eigenbasis. Using the eigenvalues of \eqref{eq:sys_1} - \eqref{eq:sys_2}, we identify two subsets of $S$: for $\delta > 0$, we define the 3-dimensional manifolds
\begin{align*}
    S_0 =& \, \left\{ \bigl(U, V, \xi\bigr): \lVert U \rVert \leq \frac{1}{\delta}, V = 0, \xi \leq \min\{\lambda_a, \lambda_0\} - \delta \right\}, \\
    S_2 =& \, \left\{ \bigl(U, V, \xi\bigr): \lVert U \rVert \leq \frac{1}{\delta}, V = 0, \max\{\lambda_a, \lambda_0\} + \delta \leq \xi \right\}.
\end{align*}
When $\lambda = -\xi + \lambda_a(U)$ and $\lambda = -\xi + \lambda_0(U)$, both eigenvalues are positive in $S_0$, and there is an unstable manifold of dimension 2 in $S_0$. In addition, both eigenvalues are negative in $S_2$, and there is a stable manifold of dimension 2 in $S_2.$ Then $(U_L, 0, -\infty)$ is an $\alpha$-limit of points in $S_0$, and $(U_R, 0, +\infty)$ is an $\omega$-limit of points in $S_2$. Given a left state $U_L$, we define the 1-dimensional invariant set
$$S_0(U_L) = \, \left\{ \bigl(U, V, \xi\bigr): U = U_L, V = 0, \xi < \min\{\lambda_a, \lambda_0\} \text{ at } U_L \right\}.$$
This line of equilibria admits a 3-dimensional unstable manifold $W_\e^u\big(S_0(U_L)\big)$, which is a perturbation of the unstable manifold
$$W_0^u\bigl(S_0(U_L)\bigr) = \, \left\{ \bigl(U, V, \xi\bigr): U \in \Omega_\xi, V = V(U), \xi < \min\{\lambda_a, \lambda_0\} \text{ at } U_L \right\},$$
where $\Omega_\xi$ is an open subset of $U$-space depending on $\xi$ and $U_L$. The function $V(U)$ is obtained by solving \eqref{eq:fast_time_sys_eneq0}. Similarly, for a given right state $U_R$, we define the 1-dimensional set
$$S_2(U_R) = \, \left\{ \bigl(U, V, \xi\bigr): U = U_R, V = 0, \max\{\lambda_a, \lambda_0\} < \xi \text{ at } U_R \right\},$$
which has a corresponding 3-dimensional stable manifold $W_\e^s\big(S_2(U_R)\big)$, the perturbation of
$$W_0^s\bigl(S_2(U_R)\bigr) = \, \left\{ \bigl(U, V, \xi\bigr): U \in \Omega_\xi, V = V(U), \max\{\lambda_a, \lambda_0\} < \xi \text{ at } U_R \right\}.$$
Our objective is to demonstrate that the perturbed invariant manifolds $W_\epsilon^u\big(S_0(U_L)\big)$ and $W_\epsilon^s\big(S_2(U_R)\big)$ intersect transversely in the extended 5-dimensional phase space. (See Figure \ref{fig:inner_soln_at_intersection}.) Any trajectory contained in this intersection corresponds to a heteroclinic orbit connecting the left state $U_L$ as $\tau \to -\infty$ to the right state $U_R$ as $\tau \to +\infty$, thus yielding a self-similar viscous profile.

To isolate the deviation from the Rankine–Hugoniot condition across the internal structure of the profile, we introduce the variable
\begin{align}
    W = -V + F(U) - \xi U,
\end{align}
so that $W = 0$ exactly when the classical Rankine-Hugoniot condition is satisfied for the shock speed $\xi = x'(t)$. This change of variables isolates the jump structure and facilitates the analysis of the profile's internal structure. We also now treat $\e$ as a state variable, and our system becomes
\begin{align}
    \begin{cases}
        \dot{U} = \bigl(F(U) - \xi U\bigr) - W, \\
        \dot{W} = -\e U, \\
        \dot{\xi} = \e, \\
        \dot{\e} = 0.
    \end{cases} \label{eq:fast_time_sys_jump_accounted}
\end{align}
We reformulate the problem in an extended 4-dimensional phase space with coordinates $(U, W, \xi, \epsilon)$. For each fixed $\e > 0$, the subspace $\e = \text{const}.$ is invariant. Within this setting, we define the following:
\begin{itemize}
    \item Two normally hyperbolic invariant manifolds:
    \begin{align*}
        T_0 &= \left\{(U, W, \xi, \epsilon) : \|U\| \leq \frac{1}{\delta},\, W = F(U) - \xi U,\, \xi \leq \min\{\lambda_a, \lambda_0\} - \delta \right\}, \\
        T_2 &= \left\{(U, W, \xi, \epsilon) : \|U\| \leq \frac{1}{\delta},\, W = F(U) - \xi U,\, \max\{\lambda_a, \lambda_0\}+\delta \leq \xi  \right\}.
    \end{align*}
    
    \item The reduced 1-dimensional invariant sets for the boundary states:
    \begin{align*}
        T^\epsilon_0(U_L) &= \left\{(U, W, \xi, \epsilon) : U = U_L,\, W = F(U_L) - \xi U_L,\,  \xi < \min\{\lambda_a, \lambda_0\} \text{ at } U_L \right\}, \\
        T^\epsilon_2(U_R) &= \left\{(U, W, \xi, \epsilon) : U = U_R,\, W = F(U_R) - \xi U_R,\, \max\{\lambda_a, \lambda_0\} < \xi \text{ at } U_R \right\}.
    \end{align*}
    
    \item The associated 3-dimensional (un)stable manifolds:
    \begin{align*}
        W^u(T^\epsilon_0(U_L)) &= \left\{(U, W, \xi, \epsilon) : U \in \Omega_\xi,\, W = W(U),\, \xi < \min\{\lambda_a, \lambda_0\} \text{ at } U_L \right\}, \\
        W^s(T^\epsilon_2(U_R)) &= \left\{(U, W, \xi, \epsilon) : U \in \Omega_\xi,\, W = W(U),\, \max\{\lambda_a, \lambda_0\} < \xi \  \text{at}\  U_R \right\},
    \end{align*}
    where \(\Omega_\xi\) is an open subset of \(U\)-space depending on \(\xi\), and \(W(U)\) is the solution of \eqref{eq:fast_time_sys_jump_accounted}.
\end{itemize}

Our goal is to show that the perturbed manifolds \(W^u(T^\e_0(U_L))\) and \(W^s(T^\e_2(U_R))\) intersect in the 5-dimensional extended phase space. Any such intersection corresponds to a viscous profile that connects \(U_L\) to \(U_R\), representing a valid traveling wave solution for fixed \(\e > 0\).

We now want to convert our system back to the variables $\rho \text{ and } u = \frac{m}{\rho}$. We have $\dot{u} = \frac{\dot{m}}{\rho} - \frac{m}{\rho^2}\dot{\rho}$ and 
$$\dot{U} = \begin{pmatrix} \dot{\rho} \\ \dot{m} \end{pmatrix} = \begin{pmatrix}
    m\left(1 - \fraction{}^a \right) - \xi \rho - w_1 \\
    \frac{m^2}{\rho} \left(1 - \fraction{}^a\right) - \xi m - w_2
\end{pmatrix}.$$
Thus,
\begin{align*}
    \dot{u} =& \, \frac{1}{\rho}\left\{\frac{m^2}{\rho} \left(1 - \fraction{}^a\right) - \xi m - w_2\right\} \\
    & \, -\frac{m^2}{\rho}\left\{m\left(1 - \fraction{}^a \right) - \xi \rho - w_1\right\} \\
    =& \, u^2 \left(1 - \fraction{}^a\right) - \xi u - \frac{w_2}{\rho} - u^2 \left(1 - \fraction{}^a\right) + u\xi + \frac{u}{\rho}w_1 \\
    =& \, \frac{u w_1 - w_2}{\rho}.
\end{align*}

Our system then becomes
\begin{align}
    \begin{cases}
        \dot{\rho} = u \rho \left(1 - \fraction{}^a\right) - \xi \rho - w_1, \\
        \dot{u} = \frac{u w_1 - w_2}{\rho}, \\
        \dot{w_1} = -\e \rho, \\
        \dot{w_2} = -\e \rho u, \\
        \dot{\xi} = \e, \\
        \dot{\e} = 0.
    \end{cases} \label{eq:before_a_value}
\end{align}
We also consider the reduced system from \eqref{eq:fast_time_sys_jump_accounted}
$$\dot{U} = F(U) - sU - W_L$$
When $\e = 0$, $\xi$ is equal to the speed of the delta shock $s=x'(t)$, $W_L = F(U_L) - sU_L$, and the equilibrium point $U_L$ acts as a source.

\begin{figure}
     \centering
     \begin{subfigure}[b]{0.48\textwidth}
         \centering
         \includegraphics[width=\textwidth]{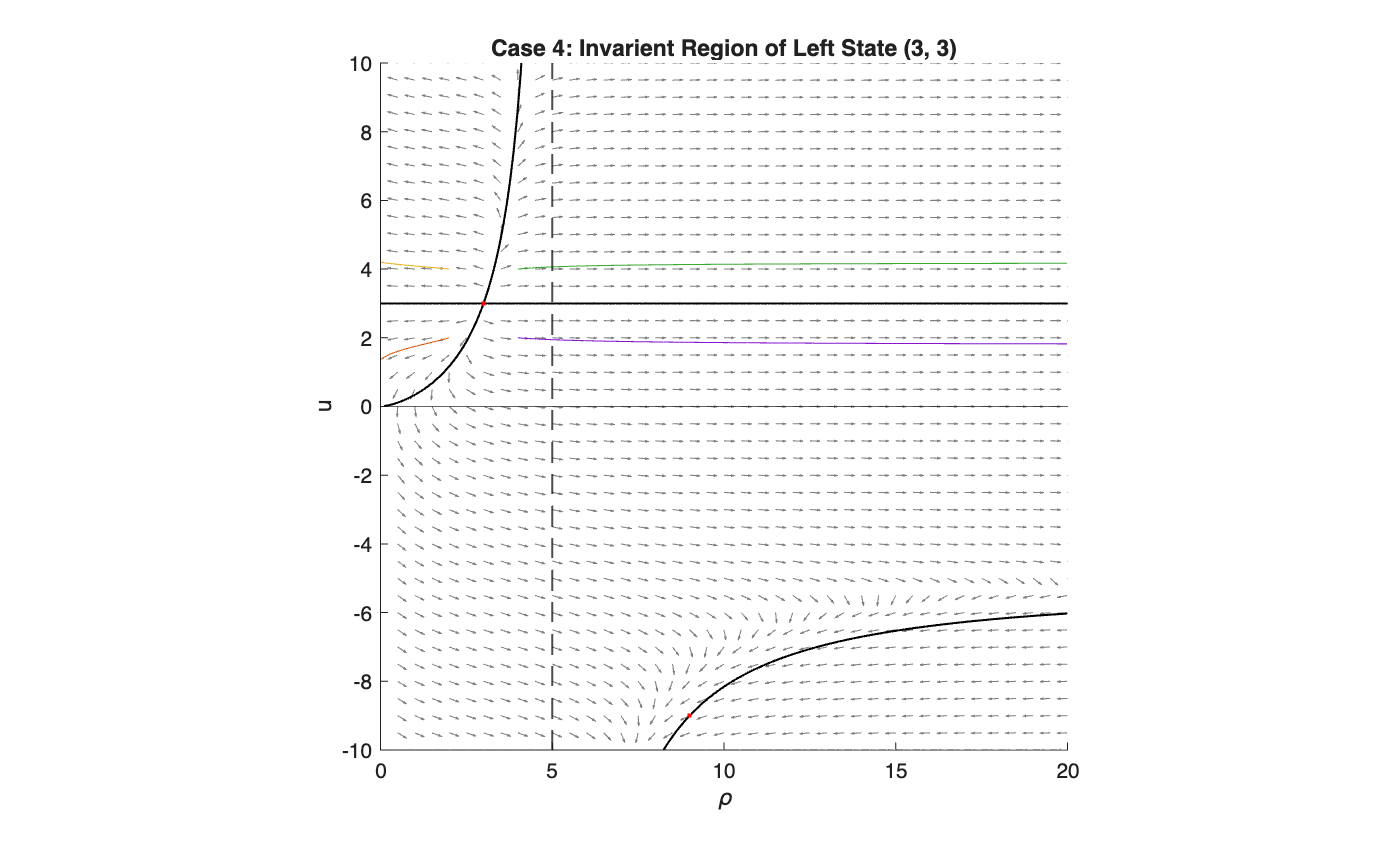}
         \caption{Invariant region of left state}
         \label{fig:11a}
     \end{subfigure}
     \begin{subfigure}[b]{0.48\textwidth}
         \centering
         \includegraphics[width=\textwidth]{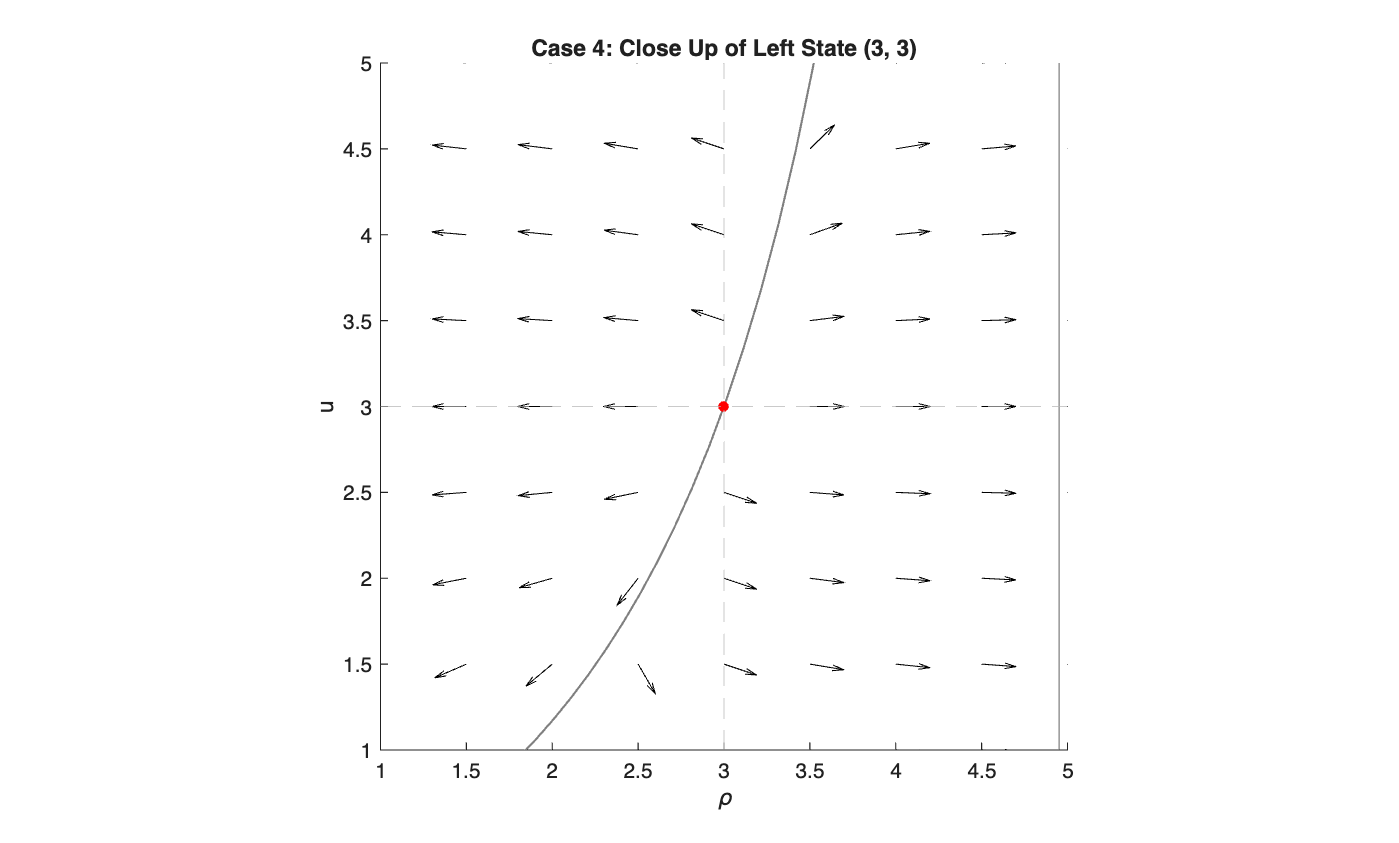}
         \caption{Close up of left state invariant region}
         \label{fig:11b}
     \end{subfigure}
     \caption{Case 4: Invariant Region of Left State. Left State: $(3,3)$. Parameters: $\bar{\rho}=5$, $a=-1.5$.}
     \label{fig:comboleft}
\end{figure}

\begin{proposition}
The planar system \( \dot{U} = F(U) - sU - W_L \) contains a negatively invariant region to the right of \( U_L \), bounded by the contact discontinuity curve of the $0$ family that passes through the left state and the straight line $u=u_L$ (see Figure \ref{fig:comboleft}).
\end{proposition}

\begin{proof}
A direct computation of \( \dot{U} \) along the curve and the horizontal line, using an argument analogous to Lemma 3.2 in \cite{Sc}, confirms that the vector field either points outward or is tangent to the boundary, following the curve or line precisely. Hence, the region is negatively invariant under the flow. 
\end{proof}

Similarly, if we consider
$$\dot{U} = F(U) - sU - W_R$$
with $\e = 0$, $\xi = s$, and \( W_R = F(U_R) - sU_R \), then the equilibrium \( U_R \) acts as a sink.

\begin{figure}
     \centering
     \begin{subfigure}[b]{0.48\textwidth}
         \centering
         \includegraphics[width=\textwidth]{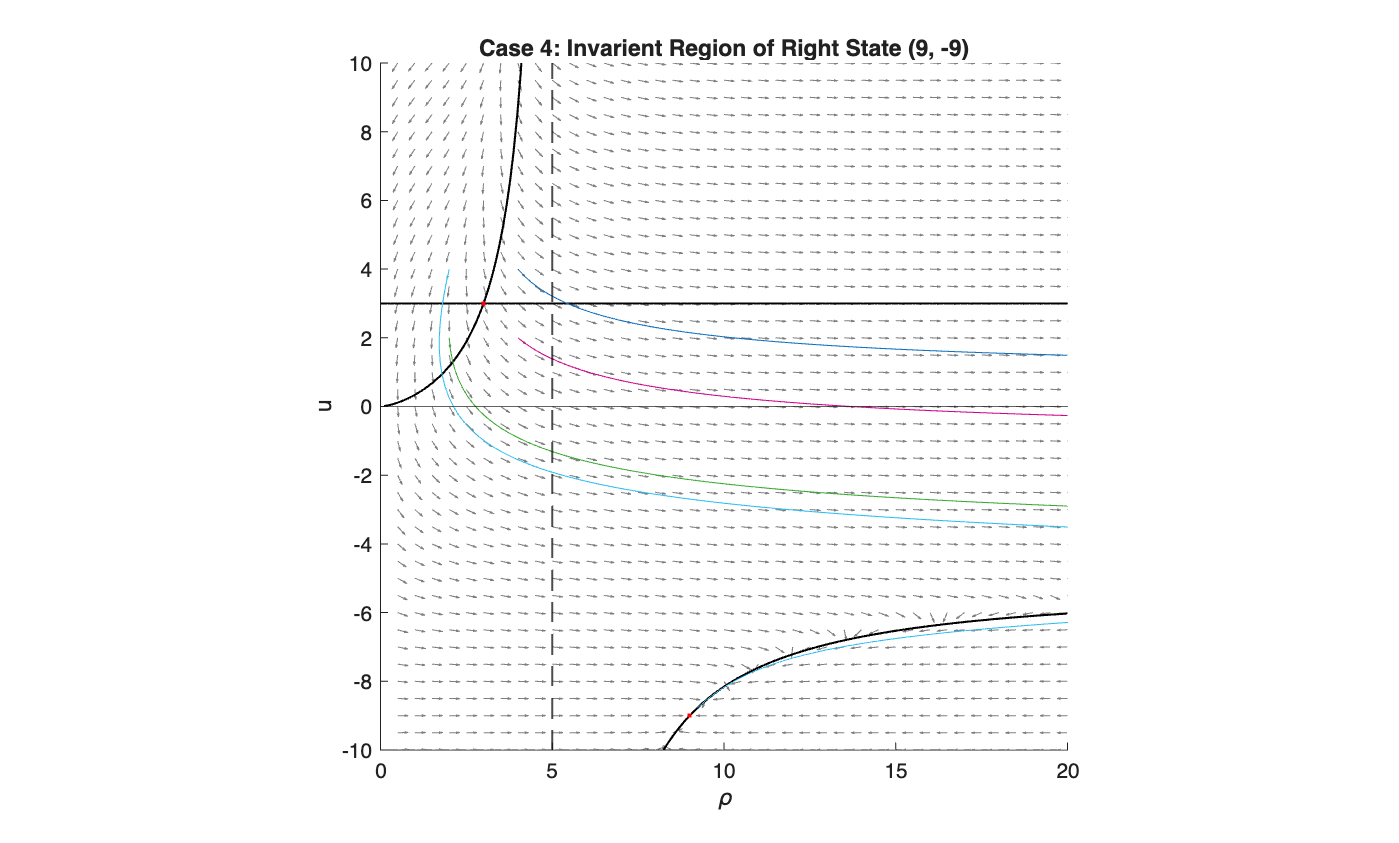}
         \caption{Invariant region of left state}
         \label{fig:12a}
     \end{subfigure}
     \begin{subfigure}[b]{0.48\textwidth}
         \centering
         \includegraphics[width=\textwidth]{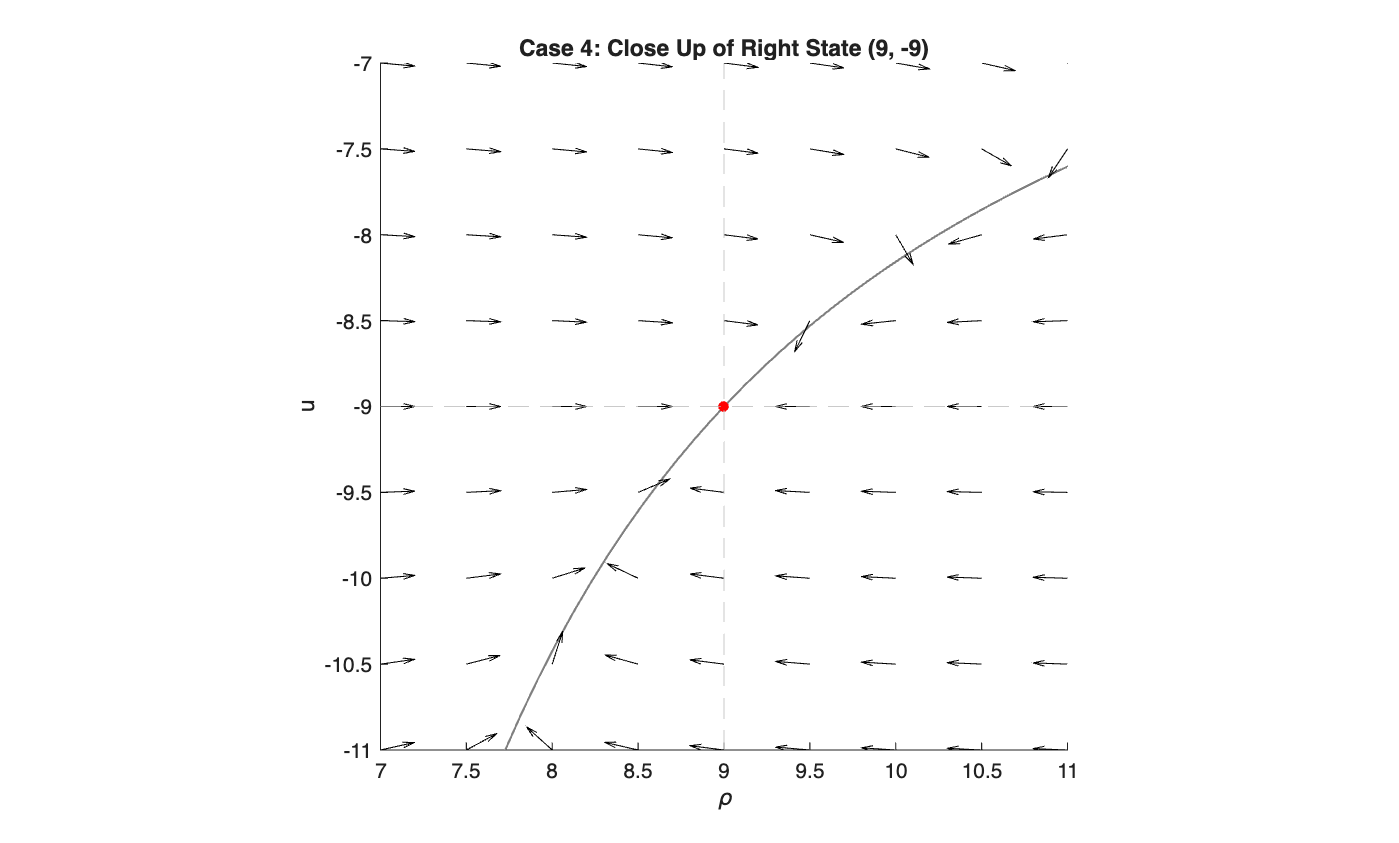}
         \caption{Close up of right state invariant region}
         \label{fig:12b}
     \end{subfigure}
     \caption{Case 4: Invariant Region of Right State. Right State: $(9,-9)$. Parameters: $\bar{\rho}=5$, $a=-1.5$.}
     \label{fig:comboright}
\end{figure}

\begin{proposition}
The planar system \( \dot{U} = F(U) - sU - W_R \) contains a positively invariant region to the right of \( U_R \), bounded by the contact discontinuity curve of the $0$ family that passes through the right state and the straight line $u=u_R$ (see Figure \ref{fig:comboright}).
\end{proposition}

\begin{proof}
A calculation of \( \cdot{U} \) along the curve and line, again similar to Lemma 3.2 in \cite{Sc}, verifies the result.
\end{proof}
In particular, the first proposition implies that all the trajectories contained within the curvilinear wedge have \( U_L \) as their \(\alpha\)-limit. Similarly, in the corresponding region near \( U_R \), the trajectories have \( U_R \) as their \(\omega\)-limit.

\begin{figure}
     \centering
     \begin{subfigure}[b]{0.48\textwidth}
         \centering
         \includegraphics[width=\textwidth]{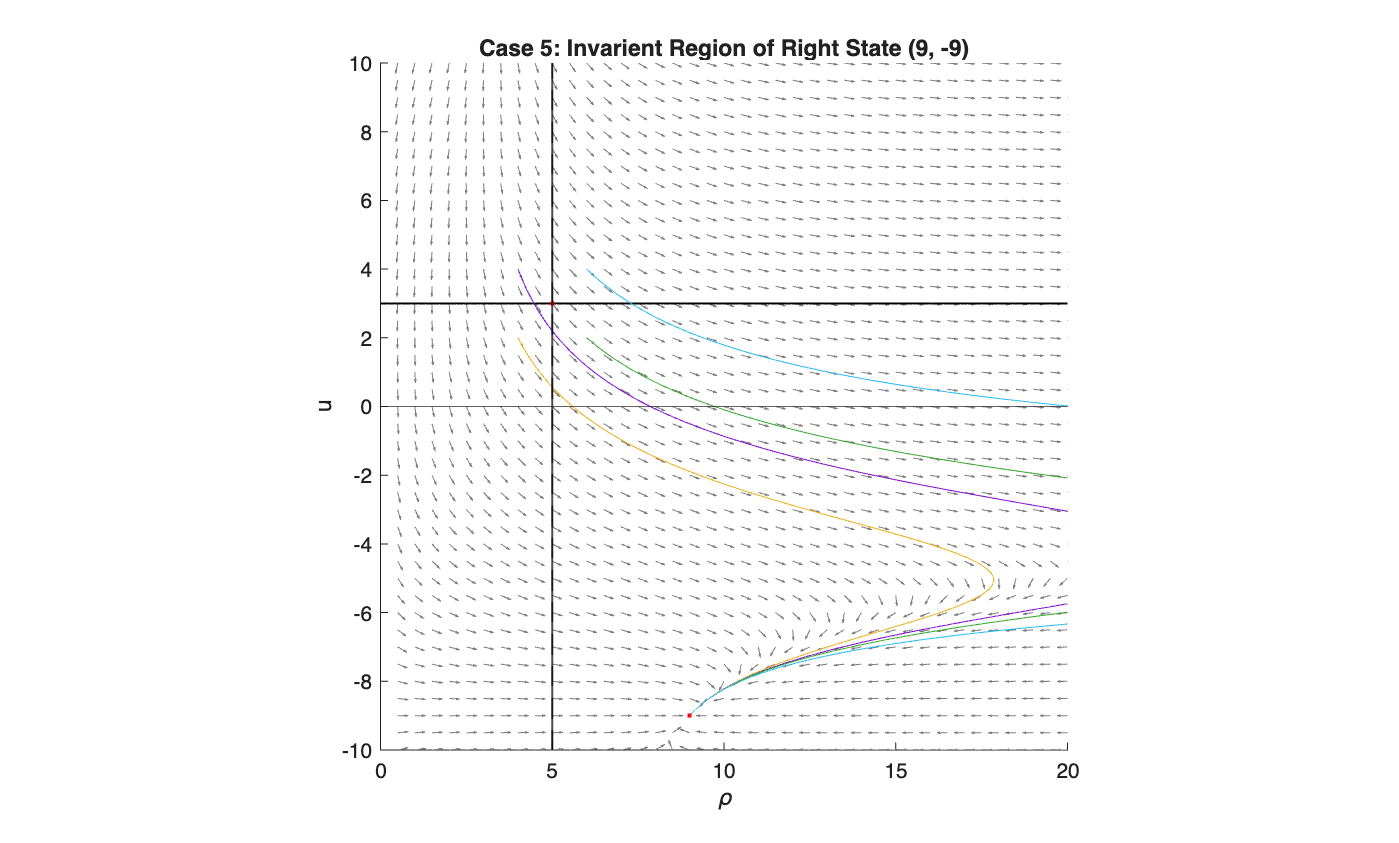}
         \caption{Invariant region of right state}
         \label{fig:13a}
     \end{subfigure}
     \begin{subfigure}[b]{0.48\textwidth}
         \centering
         \includegraphics[width=\textwidth]{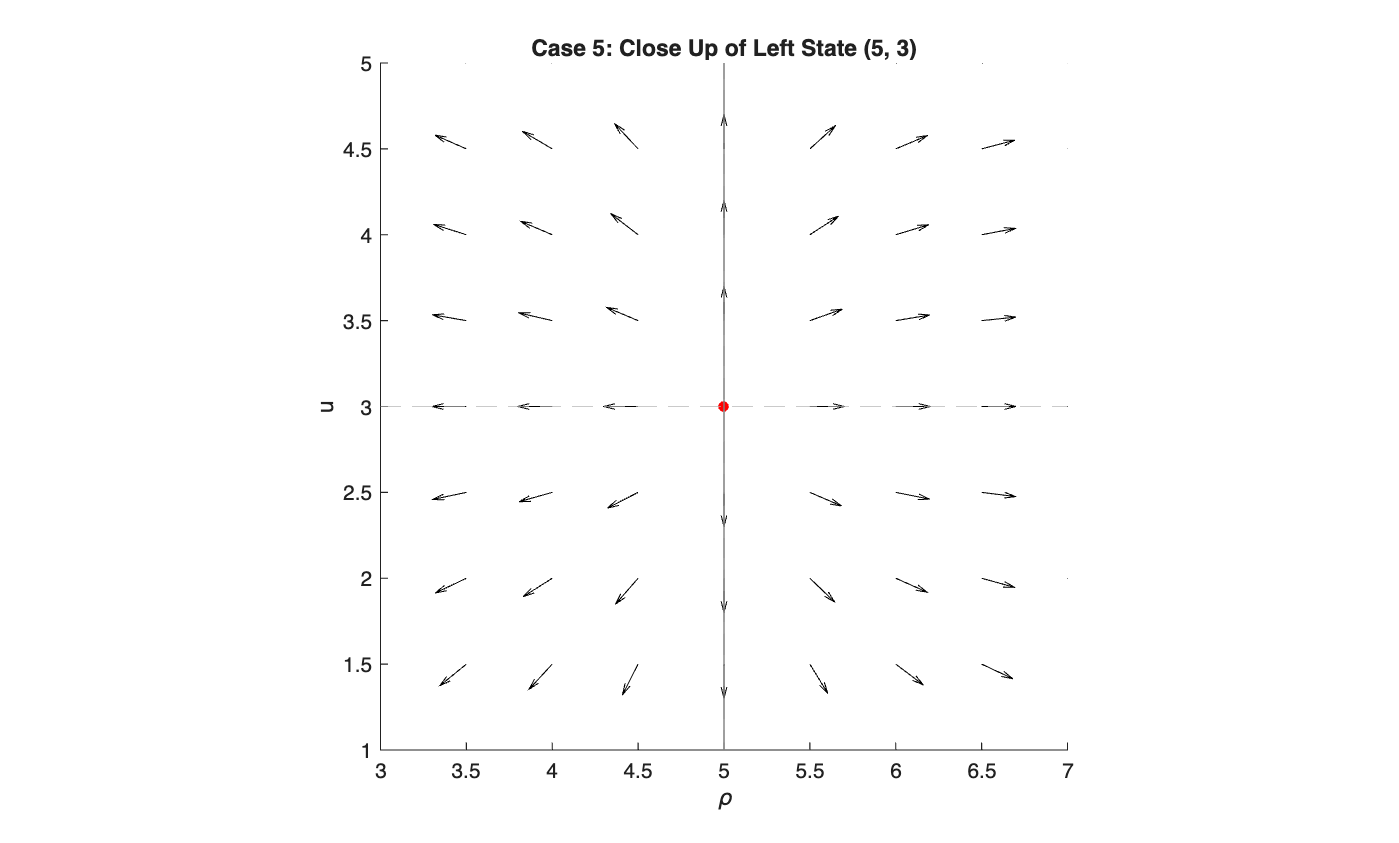}
         \caption{Close up of left state invariant region}
         \label{fig:13b}
     \end{subfigure}
     \caption{Case 5: Invariant Region of Right and Left States, where $\rho_L =\bar{\rho} $. Right State: $(9,-9)$. Left State: $(5,3)$.  Parameters: $\bar{\rho}=5$, $a=-1.5$.}
     \label{fig:case5}
\end{figure}
\begin{figure}
     \centering
     \begin{subfigure}[b]{0.48\textwidth}
         \centering
         \includegraphics[width=\textwidth]{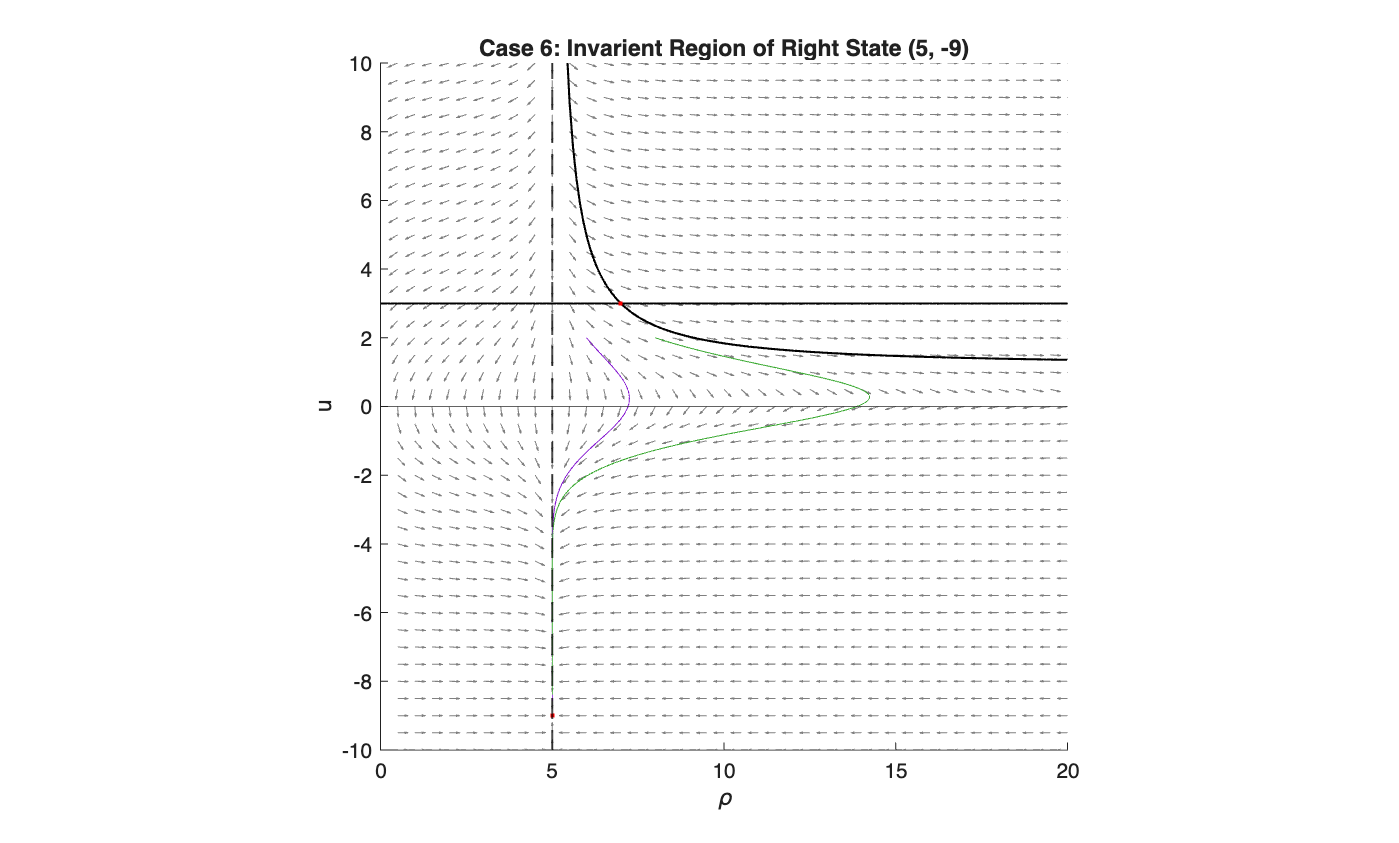}
         \caption{Invariant region of right state}
         \label{fig:14a}
     \end{subfigure}
     \begin{subfigure}[b]{0.48\textwidth}
         \centering
         \includegraphics[width=\textwidth]{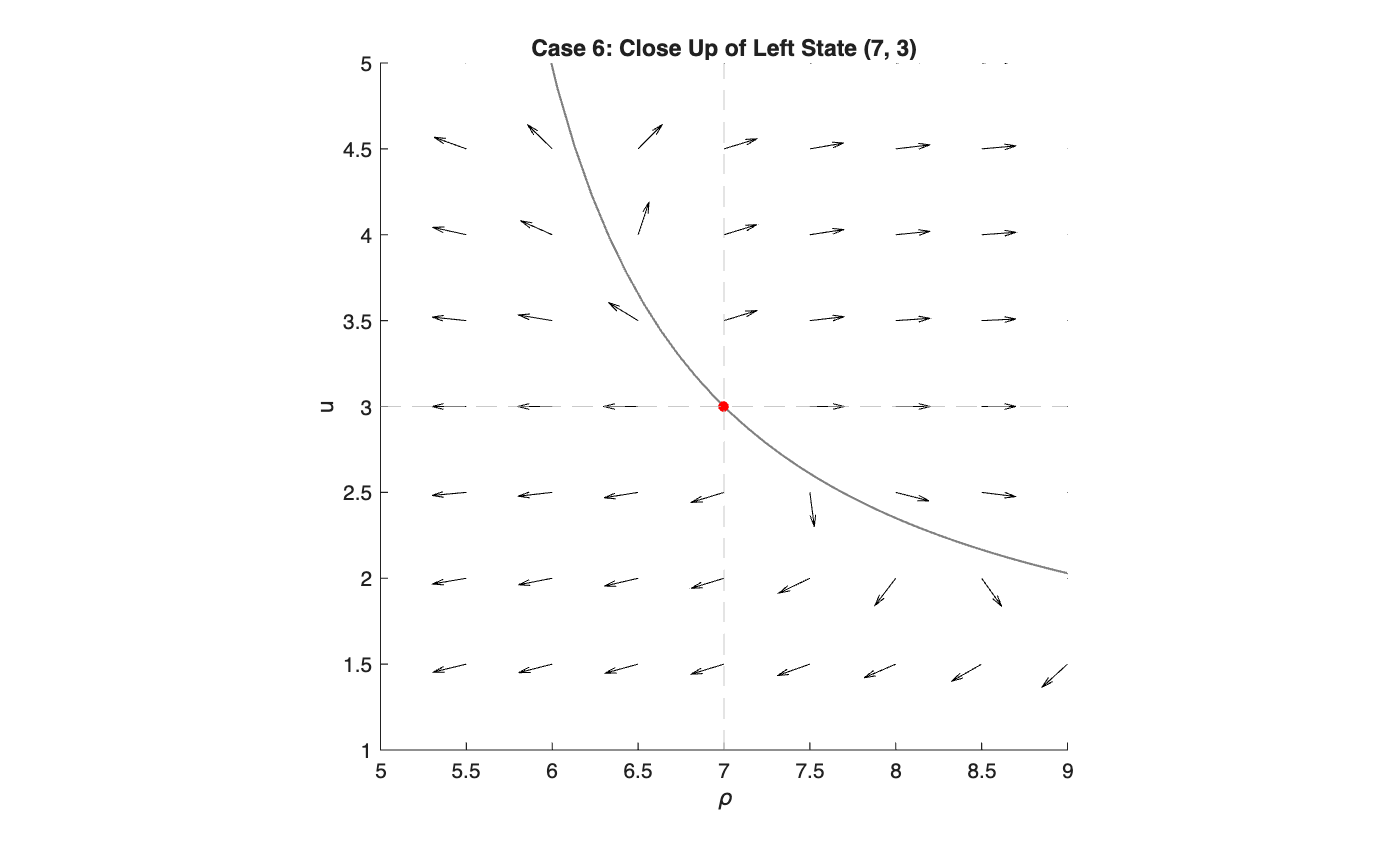}
         \caption{Close up of left state invariant region}
         \label{fig:14b}
     \end{subfigure}
     \caption{Case 6: Invariant Region of Right and Left States, where $\rho_R =\bar{\rho} $. Right State: $(5,-3)$. Left State: $(7,3)$.  Parameters: $\bar{\rho}=5$, $a=-1.5$.}
     \label{fig:case6}
\end{figure}

Note that these propositions still hold even when $\rho_L=\bar{\rho}$ or $\rho_R=\bar{\rho}$, as shown by Figures \ref{fig:case5} and \ref{fig:case6}.

\subsection{Delta Shocks for \texorpdfstring{$a < 0$}{a < 0}: The Blow-up Procedure} \label{GSPT_anegative}
Using the delta shock predicted by Section \ref{marko_anegative}, we now consider $\rho = \frac{\rho_0}{\e}$ (note that $\rho_0(\pm \infty)=0$). Since $\e$ is a constant with respect to $\tau$, our system then becomes
\begin{align}
    \begin{cases}
        \dot{\rho}_0 = u \rho_0 \left(1 - \left(\frac{\, \rho_0 \,}{\e \rhobar}\right)^a\right) - \xi \rho_0 - \e w_1, \\
        \dot{u} = \frac{(u w_1 - w_2)\e}{\rho_0}, \\
        \dot{w_1} = -\rho_0, \\
        \dot{w_2} = -u \rho_0, \\
        \dot{\xi} = \e, \\
        \dot{\e} = 0.
    \end{cases}
\end{align}
When $\e = 0$, the system has a fixed point at $\rho_0 = 0$ (also note that $\dot{\rho}_0=(u-\xi)\rho_0$ and $\dot{u} =0$). To analyze the dynamics near it, we first \emph{desingularize}, which involves a rescaling of time. Specifically, we multiply the right-hand side of the system by a factor of $\rho_0$, effectively changing the time variable. This is not merely an algebraic manipulation: it corresponds to a change of time scale that slows down the flow near the singularity and reveals a hidden structure in the dynamics. 

At the point $(\rho_0, u, \epsilon) = (0, \xi, 0)$, the fast system loses normal hyperbolicity. When we linearize the desingularized system at this equilibrium, we find that some eigenvalues of the Jacobian matrix vanish in directions transverse to the critical manifold. As a result, the classical theory of Fenichel does not apply. To analyze the dynamics near this degenerate point, we perform a blow-up transformation that replaces the singularity with a sphere. This transformation resolves the non-hyperbolic structure and reveals directions in different charts, allowing the use of tools such as the Exchange Lemma to rigorously track trajectories through this critical region. Thus, we introduce the \emph{blow-up transformation} by setting
\[\rho_0 = r \, \overline{\rho}_0, \quad u - \xi = r \, \overline{u}_0, \quad \e = r \, \overline{\e}, \quad \text{with} \quad \overline{\rho}_0^2 + \overline{u}_0^2 + \overline{\e}^2 = 1.\]
At this stage, we note that
\begin{align*}
    \dfrac{\rho_0}{\tau} &= \dfrac{r}{\tau} \rhozerobar + r \dfrac{\rhozerobar}{\tau} \\
    &= r \rhozerobar (r \uzerobar + \xi) - r \rhozerobar (r \uzerobar + \xi) \left(\frac{\, r \rhozerobar \,}{r \ebar \, \rhobar}\right)^a - \xi r \rhozerobar - r \ebar w_1, \\
    \implies \dfrac{\rhozerobar}{\tau} &= -\frac{\rhozerobar}{r} \dfrac{r}{\tau} + \rhozerobar (r \uzerobar + \xi) - \rhozerobar (r \uzerobar + \xi) \left(\frac{\rhozerobar}{ \ebar \, \rhobar}\right)^a - \xi  \rhozerobar -  \ebar w_1; \\
    \dfrac{(u - \xi)}{\tau} &= \dfrac{r}{\tau} \uzerobar + r \dfrac{\uzerobar}{\tau} = \frac{\bigl((r \uzerobar + \xi) w_1 - w_2\bigr) r \ebar}{r \rhozerobar} - r \ebar, \\
    \implies \dfrac{\uzerobar}{\tau} &= -\frac{\uzerobar}{r} \dfrac{r}{\tau} + \frac{\bigl((r \uzerobar + \xi) w_1 - w_2\bigr) \ebar}{r \rhozerobar} - \ebar; \\
    0 &= \dfrac{\e}{\tau} = \dfrac{r}{\tau} \ebar + r \dfrac{\ebar}{\tau} \implies \dfrac{\ebar}{\tau} = -\frac{\e}{r} \dfrac{r}{\tau}.
\end{align*}
Combining these with the relation $\rhozerobar \dot{\overline{\rho}}_0 + \overline{\u}_0 \dot{\overline{\u}}_0 + \ebar \, \dot{\ebar} = 0$, we obtain
$$\dfrac{r}{\tau} = \frac{\bigl((r \uzerobar + \xi) w_1 - w_2\bigr) \ebar \, \uzerobar}{\rhozerobar} - r \ebar \, \uzerobar + r \rhozerobar^2 (r \uzerobar + \xi) \left(1 - \left(\frac{\, \rhozerobar \,}{\ebar \, \rhobar}\right)^a\right) - \xi r \rhozerobar^2 - r \ebar \, \rhozerobar w_1.$$
Upon substituting the blow-up transformation into the previously derived system, we obtain a reformulated system where $r$ appears explicitly as an additional variable
\begin{align}\label{5.13}
    \begin{cases}
        \dfrac{r}{\tau} = \frac{\bigl((r \uzerobar + \xi) w_1 - w_2\bigr) \ebar \, \uzerobar}{\rhozerobar} - r \ebar \, \uzerobar + r \rhozerobar^2 (r \uzerobar + \xi) \left(1 - \left(\frac{\rhozerobar}{\, \ebar \, \rhobar \,}\right)^a\right) - \xi r \rhozerobar^2 - r \ebar \, \rhozerobar w_1, \\ \\
        \dfrac{\rhozerobar}{\tau} = \rhozerobar (r \uzerobar + \xi) \left(1 - \left(\frac{\rhozerobar}{\, \ebar \, \rhobar \,}\right)^a\right) - \xi \rhozerobar - \ebar w_1 - \frac{\bigl((r \uzerobar + \xi) w_1 - w_2\bigr) \ebar \, \uzerobar}{r} + \ebar \, \rhozerobar \, \uzerobar \\ \\ 
        \qquad\quad\!\! - \rhozerobar^3 (r \uzerobar + \xi) \left(1 - \left(\frac{\rhozerobar}{\, \ebar \, \rhobar \,}\right)^a\right) + \xi \, \rhozerobar^3 + \ebar \, \rhozerobar^2 w_1, \\ \\ 
        \dfrac{\uzerobar}{\tau} = \frac{\bigl((r \uzerobar + \xi) w_1 - w_2\bigr) \ebar}{r \rhozerobar} - \ebar - \frac{\bigl((r \uzerobar + \xi) w_1 - w_2\bigr) \ebar \, \uzerobar^2}{r \rhozerobar} + \ebar \, \uzerobar^2 +\xi \, \rhozerobar^2 \, \uzerobar \\ \\ 
        \qquad\quad\!\! - \rhozerobar^2 \, \uzerobar (r \uzerobar + \xi) \left(1 - \left(\frac{\rhozerobar}{\, \ebar \, \rhobar \,}\right)^a\right) - \ebar \, \rhozerobar \uzerobar w_1, \\ \\ 
        \dfrac{\ebar}{\tau} = -\frac{\bigl((r \uzerobar + \xi) w_1 - w_2\bigr) \ebar^2 \uzerobar}{r \rhozerobar} + \ebar^2 \uzerobar - \ebar \, \rhozerobar^2 (r \uzerobar + \xi) \left(1 - \left(\frac{\rhozerobar}{\, \ebar \, \rhobar \,}\right)^a\right) + \xi \, \ebar \, \rhozerobar^2 + \ebar^2 \rhozerobar w_1, \\ \\
        \dfrac{w_1}{\tau} = -r \rhozerobar, \\ \\ 
        \dfrac{w_2}{\tau} = -r \rhozerobar (r \uzerobar + \xi), \\ \\
        \dfrac{\xi}{\tau} = r \ebar.
    \end{cases}
\end{align}
We proceed with a second desingularization (multiplying by $r$ and $\rhozerobar \geq 0$ and dividing by $\ebar\geq 0$) and then linearize the resulting system and evaluate it at $r = 0$, which yields
\begin{align}\label{5.14}
    \begin{cases}
        \dfrac{r}{\tau} = 0 = \dfrac{w_1}{\tau} = \dfrac{w_2}{\tau} = \dfrac{\xi}{\tau}, \\ \\ 
        \dfrac{\rhozerobar}{\tau} = - \uzerobar \, \rhozerobar \left(\xi w_1 - w_2\right), \\ \\ 
        \dfrac{\uzerobar}{\tau} = (\xi w_1 -  w_2)\left(1 - \uzerobar^2\right), \\ \\
        \dfrac{\ebar}{\tau} = - \ebar \, \uzerobar \left(\xi w_1 - w_2\right).
    \end{cases}
\end{align}
Note that before setting $r=0$ we have the equation for the radial variable $\dfrac{r}{\tau} = \bigl(\xi w_1 - w_2\bigr) \uzerobar r+O(r^2),$ which governs the behavior of trajectories near the sphere.
Recall now that $\dfrac{w_1}{\tau} = - \rho_0$, so
$$w_1(\tau) = -\int_{-\infty}^\tau \rho_0(s)\d s + w_1(-\infty).$$
Likewise, since $\dfrac{w_2}{\tau} = - \rho_0 u$, we have
\begin{align*}
    w_2(\tau) &= -\int_{-\infty}^\tau \rho_0 u \d s + w_2(-\infty) \\
    &= w_2(-\infty) - \int_{-\infty}^\tau u \d\left(\int_{-\infty}^s\rho_0 \d\theta\right) \\
    &= w_2(-\infty) - \left(u \int_{-\infty}^s\rho_0 \d\theta\right) \bigg\vert_{-\infty}^\tau + \int_{-\infty}^\tau \left\{\dfrac{u(s)}{s} \int_{-\infty}^s \rho_0 \d\theta \right\}\d s.
\end{align*}
Using that
\begin{align*}
    w_2(-\infty) &= -V_{2,L} + F_2(U_L) - \xi U_{2,L} \\
    &= -0 + \rho_L \u_L^2 \left(1 - \fraction{L}^a\right) - \xi \rho_L \u_L \\
    &= u_L\left(\rho_L \u_L \left(1 - \fraction{L}^a\right) - \xi \rho_L \right) \\
    &= u_L w_1(-\infty),
\end{align*}
we obtain
$$\xi w_1 - w_2 = (\xi - \u_L) w_1(-\infty) - (\xi - u)\int_{\infty}^\tau \rho_0 \d s - \int_{-\infty}^\tau \left\{\dfrac{u(s)}{s} \int_{-\infty}^s \rho_0 \d\theta \right\}\d s.$$
When $r = 0$ (which corresponds to $\e = 0$) and $u = \xi$ (with $\xi$ constant and equal to the speed of the delta shock), the expression simplifies to
$$\xi w_1 - w_2 = (\xi - u_L) w_1(-\infty) = (\xi - u_L) \rho_L \left( \u_L \left(1 - \fraction{L}^a\right) - \xi \right),$$
which is negative in the inner layer. This holds because $\rho_L > 0$, $\xi < \lambda_0(U_L)$ (assumed in equation \eqref{eq:overcompressive}), and $\xi < u_L.$ As we discussed in Section \ref{numerics}, in each of the different cases, any region such that $u_R > u_L$ is accessible when $a < 0$ via classical waves. Thus, if a singular shock appeared (which we know does), it must be the case that $u_R < \xi < u_L$. We then return to \eqref{5.13}, desingularize once more by dividing by $-(\xi w_1 - w_2)$ (which corresponds to a rescaling of time), let $r=0$, and observe that the points $\uzerobar = \pm 1$ with $\rhozerobar = \ebar = 0$ are fixed points. Linearizing the system about the fixed points $\uzerobar = \pm 1$ and $\rhozerobar = \ebar = 0$ yields eigenvalues $\lambda = \uzerobar, 2\uzerobar, \uzerobar, \text{and } -\uzerobar$ in the $\rhozerobar, \uzerobar, \ebar, \text{and } r$ directions, respectively. We again note that $\dfrac{r}{\tau} = -\uzerobar r+O(r^2)$. This expression indicates that when \( \uzerobar = 1 \), we have \( \frac{dr}{d\tau} = -r + \mathcal{O}(r^2) < 0 \) for small \( r > 0 \), so trajectories enter the blown up sphere transversely near the point \( \uzerobar = 1 \); when \( \uzerobar = -1 \), we have \( \frac{dr}{d\tau} = r + \mathcal{O}(r^2) > 0 \), so trajectories exit the sphere transversely near \( \uzerobar = -1 \). 

The eigenvalues of the system on the blown-up sphere are nonzero, and the linearization yields a full set of linearly independent eigenvectors. Therefore, the fixed points are normally hyperbolic, and Fenichel's theory ensures the persistence of stable and unstable manifolds under perturbation.

\begin{figure}{h}
    \centering
    \begin{tikzpicture}[scale=2]
        \draw[->] (0,0)--(1.2,0) node[right] {$\ebar$};
        \draw[->] (0,0)--(0,1.2) node[above] {$\uzerobar$};
        \draw[->] (0,0)--(-.9,-.6) node[below left] {$\rhozerobar$};
        \draw[dashed] (0,0)--(-1,0);
        \draw[] (0,0)--(0,-1);
        \draw[dashed] (0,0)--(.285,.19);
    
        \draw plot[domain=-pi/2:pi/2] ({cos(\x r)},{sin(\x r)});
        \draw[dashed] plot[domain=pi/2:(pi+pi/2)] ({cos(\x r)},{sin(\x r)});
        \draw plot[domain=4.4:2*pi] ({cos(\x r)},{.2*sin(\x r)});
        \draw[dashed] plot[domain=0:4.4] ({cos(\x r)},{.2*sin(\x r)});
        \draw plot[domain=pi:2*pi] ({.293*sin(\x r)},{cos(\x r)});
        \draw[dashed] plot[domain=0:pi] ({.293*sin(\x r)},{cos(\x r)});
    
        \draw[->] (-.5,1.1)--(0,1);
        \draw[->] (0,-1)--(.5,-1.1);
        \draw[->] ({.293*sin((pi+3*pi/4) r)},{cos((pi+3*pi/4) r)})--({.293*sin((pi+3*pi/4) r)-0.293*cos((pi+3*pi/4) r)/3},{cos((pi+3*pi/4) r)+sin((pi+3*pi/4) r)/3});
        \draw[->] ({.293*sin((pi+2*pi/4) r)},{cos((pi+2*pi/4) r)})--({.293*sin((pi+2*pi/4) r)-0.293*cos((pi+2*pi/4) r)/3},{cos((pi+2*pi/4) r)+sin((pi+2*pi/4) r)/3});
        \draw[->] ({.293*sin((pi+1*pi/4) r)},{cos((pi+1*pi/4) r)})--({.293*sin((pi+1*pi/4) r)-0.293*cos((pi+1*pi/4) r)/3},{cos((pi+1*pi/4) r)+sin((pi+1*pi/4) r)/3});

        \draw[->] ({sin((3*pi/4) r)},{cos((3*pi/4) r)})--({sin((3*pi/4) r)+cos((3*pi/4) r)/3},{cos((3*pi/4) r)-sin((3*pi/4) r)/3});
        \draw[->] ({sin((2*pi/4) r)},{cos((2*pi/4) r)})--({sin((2*pi/4) r)+cos((2*pi/4) r)/3},{cos((2*pi/4) r)-sin((2*pi/4) r)/3});
        \draw[->] ({sin((1*pi/4) r)},{cos((1*pi/4) r)})--({sin((1*pi/4) r)+cos((1*pi/4) r)/3},{cos((1*pi/4) r)-sin((+1*pi/4) r)/3});

        \draw plot[domain=0:pi] ({0.5*sin(\x r)},{cos(\x r)});
        \draw[->] ({0.5*sin((3*pi/4) r)},{cos((3*pi/4) r)})--({0.5*sin((3*pi/4) r)+0.5*cos((3*pi/4) r)/3},{cos((3*pi/4) r)-sin((3*pi/4) r)/3});
        \draw[->] ({0.5*sin((2*pi/4) r)},{cos((2*pi/4) r)})--({0.5*sin((2*pi/4) r)+0.5*cos((2*pi/4) r)/3},{cos((2*pi/4) r)-sin((2*pi/4) r)/3});
        \draw[->] ({0.5*sin((1*pi/4) r)},{cos((1*pi/4) r)})--({0.5*sin((1*pi/4) r)+0.5*cos((1*pi/4) r)/3},{cos((1*pi/4) r)-sin((+1*pi/4) r)/3});
    \end{tikzpicture}
    \caption{Illustration of the blow-up transformation and the inner solutions}
    \label{fig:blown_up_sphere}
\end{figure}
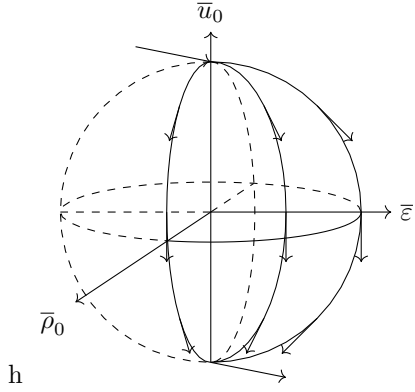

We seek a heteroclinic orbit on the blown-up sphere connecting the fixed points \( \uzerobar = 1 \) and \( \uzerobar = -1 \). Since the flow on the sphere is regular after desingularization, such a connecting orbit corresponds to the inner solution of a viscous profile. The flow for \( \uzerobar \in (-1, 1) \) satisfies \( \dfrac{\uzerobar}{\tau} = - (1 - \uzerobar^2) < 0 \) and \( \dfrac{\rhozerobar}{\tau} =\uzerobar \rhozerobar \), which implies that the system flows from \( \uzerobar = 1 \) to \( \uzerobar = -1 \) in the sphere. Therefore, there exists a heteroclinic orbit on the blown-up sphere connecting the north pole to the south pole. See Figure \ref{fig:blown_up_sphere}. This orbit corresponds to the inner layer of the singular shock profile. The Exchange Lemma ensures that this connection persists for sufficiently small \( \e > 0 \). Specifically, since the entry and exit points are normally hyperbolic and the dimensions of the corresponding stable and unstable manifolds match (both of dimension 3), a transverse intersection of these manifolds implies the existence of a trajectory in the perturbed system connecting the left and right states.

\subsection{Delta Shocks for \texorpdfstring{$a > 0$}{a > 0}: The Blow-up Procedure} \label{GSPT_apositive}
Again, following the suggested values from the shadow wave analysis in Section \ref{marko_apositive}, we consider $\rho = \frac{\rho_0}{\e}$ when $a > 0$. We substitute this into the system \eqref{eq:before_a_value} and obtain
\begin{align}
    \begin{cases}
        \dfrac{\rho_0}{\tau} = \rho_0 \u \left(1 - \left(\frac{\rho_0}{\e \rhobar}\right)^a\right) - \xi \rho_0 - \e w_1, \\
        \dot{u} = \frac{(\u w_1 - w_2)}{\rho_0} \e, \\
        \dot{w_1} = - \rho_0, \\
        \dot{w_2} = -\rho_0 \u, \\
        \dot{\xi} = \e, \\
        \dot{\e} = 0.
    \end{cases}
\end{align}
Note that if $u=\e^a u_0$, as discussed in Section \ref{marko_apositive}, then $\dfrac{\rho_0}{\tau} = -\rho_0 \u_0  \left(\frac{\rho_0}{\rhobar}\right)^a - \xi \rho_0$ when $\e=0.$
Again we proceed with our blow-up procedure. Noting that the value of $u$ at which the singularity occurred according to the shadow wave analysis was $\u = 0$, let $ \rho_0 = r^{k_1} \rhozerobar, \u = r^{k_2} \uzerobar, \e = r^{k_3} \ebar$ subject to the constraint $\uzerobar^2 + \rhozerobar^2 + \ebar^2 = 1$, where $k_i \in \mathbb{R}_+$ are scaling exponents to be determined based on the balancing of leading-order terms in the blown-up equations.

Observe that
\begin{align*}
    \dfrac{\rho_0}{\tau} =& \, k_1 r^{k_1-1} \dfrac{r}{\tau} \rhozerobar + r^{k_1} \dfrac{\rhozerobar}{\tau} \\
    =& \,  r^{k_1} \rhozerobar r^{k_2} \uzerobar \left(1 - \left(\frac{r^{k_1} \, \rhozerobar}{r^{k_3} \, \ebar \, \rhobar}\right)^a\right) - \xi r^{k_1} \, \rhozerobar - r^{k_3} \, \ebar w_1, \\
    \implies \dfrac{\rhozerobar}{\tau} =& \, -\frac{k_1 \rhozerobar}{r} \dfrac{r}{\tau} + \rhozerobar \uzerobar r^{k_2} \left(1 - r^{\left(ak_1 - ak_3\right)}\left(\frac{\rhozerobar}{\ebar \, \rhobar}\right)^a\right) -\xi \rhozerobar - r^{k_3 - k_1} \, \ebar w_1; \\
    \dfrac{\u}{\tau} =& \, k_2 r^{k_2-1} \dfrac{r}{\tau} \uzerobar + r^{k_2} \dfrac{\uzerobar}{\tau} = \frac{\bigl( r^{k_2} \uzerobar w_1 - w_2\bigr) r^{k_3} \, \ebar}{r^{k_1} \rhozerobar}, \\
    \implies \dfrac{\uzerobar}{\tau} =& \, -\frac{k_2 \uzerobar}{r} \dfrac{r}{\tau} + \frac{\bigl(r^{k_2} \uzerobar w_1 - w_2\bigr)}{\rhozerobar} r^{\left(k_3 - k_1 - k_2\right)} \ebar; \\
    0 =& \, \dfrac{\e}{\tau} = k_3 r^{\left(k_3 - 1\right)} \dfrac{r}{\tau} \ebar + r^{k_3} \dfrac{\ebar}{\tau} \implies \dfrac{\ebar}{\tau} = -\frac{k_3 \ebar}{r} \dfrac{r}{\tau}.
\end{align*}

Let $\kappa =  k_1 \rhozerobar^2 + k_2 \uzerobar^2 + k_3 \ebar^2$. Using our relation $\uzerobar \dot{\overline{u}}_0 + \rhozerobar \dot{\overline{\rho}}_0 + \ebar \, \dot{\overline{\e}}_0 = 0$, we obtain
\begin{align*}
    \kappa \dfrac{r}{\tau} =& \, \frac{\uzerobar \bigl(r^{k_2} \uzerobar w_1 - w_2\bigr) \ebar}{\rhozerobar} r^{\left(k_3 - k_1 - k_2 + 1\right)} + \rhozerobar^2 \uzerobar r^{\left(k_2 + 1\right)}\left(1 - r^{\left(ak_1 - ak_3\right)}\left(\frac{\rhozerobar}{\ebar \, \rhobar}\right)^a\right) \\
    & - \xi r \rhozerobar^2 - \rhozerobar r^{\left(k_3 - k_1 + 1\right)} \, \ebar w_1.
\end{align*}

The full system is
\begin{align}
    \begin{cases}
        \dfrac{r}{\tau} = \frac{1}{\kappa}\frac{\uzerobar \bigl(r^{k_2} \uzerobar w_1 - w_2\bigr) \ebar}{\rhozerobar} r^{\left(k_3 - k_1 - k_2 + 1\right)} + \frac{1}{\kappa} \rhozerobar^2 \uzerobar r^{\left(k_2 + 1\right)}\left(1 - r^{ak_1 - ak_3}\left(\frac{\rhozerobar}{\ebar \, \rhobar}\right)^a\right) \\ \\
        \qquad\enspace - \frac{1}{\kappa} \xi r \rhozerobar^2 - \frac{1}{\kappa} \rhozerobar r^{\left(k_3 - k_1 + 1\right)} \, \ebar w_1, \\ \\
        \dfrac{\rhozerobar}{\tau} = -\frac{k_1}{\kappa} \uzerobar \left(r^{k_2} \uzerobar w_1 - w_2\right) \ebar r^{\left(k_3 - k_1 - k_2\right)} - \frac{k_1}{\kappa}\rhozerobar^3 \uzerobar r^{k_2} \left(1 - r^{ak_1 - ak_3}\left(\frac{\rhozerobar}{\ebar \, \rhobar}\right)^a\right) \\ \\
        \qquad\quad\! + \frac{k_1}{\kappa} \xi \rhozerobar^3 + \frac{k_1}{\kappa} r^{\left(k_3 - k_1\right)}\rhozerobar^2 \ebar w_1 + \rhozerobar \uzerobar r^{k_2} \left(1 - r^{ak_1 - ak_3}\left(\frac{\rhozerobar}{\ebar \, \rhobar}\right)^a\right) - \xi \rhozerobar - r^{\left(k_3 - k_1\right)} \, \ebar w_1, \\ \\
        \dfrac{\uzerobar}{\tau} = -\frac{k_2}{\kappa}\frac{\uzerobar^2 \bigl(r^{k_2} \uzerobar w_1 - w_2\bigr) \ebar}{\rhozerobar} r^{\left(k_3 - k_1 - k_2\right)} - \frac{k_2}{\kappa} \rhozerobar^2 \uzerobar^2 r^{k_2}\left(1 - r^{ak_1 - ak_3}\left(\frac{\rhozerobar}{\ebar \, \rhobar}\right)^a\right) \\ \\
        \qquad\quad\! + \frac{k_2}{\kappa} \xi \rhozerobar^2 \uzerobar + \frac{k_2}{\kappa} \rhozerobar \uzerobar r^{\left(k_3 - k_1\right)} \, \ebar w_1 + \frac{\bigl(r^{k_2} \uzerobar w_1 - w_2\bigr)}{\rhozerobar} r^{\left(k_3 - k_1 - k_2\right)} \ebar, \\ \\
        \dfrac{\ebar}{\tau} = -\frac{k_3}{\kappa}\frac{\uzerobar \bigl(r^{k_2} \uzerobar w_1 - w_2\bigr) \ebar^2}{\rhozerobar} r^{\left(k_3 - k_1 - k_2\right)} - \frac{k_3}{\kappa} \rhozerobar^2 \uzerobar r^{k_2} \ebar \left(1 - r^{ak_1 - ak_3}\left(\frac{\rhozerobar}{\ebar \, \rhobar}\right)^a\right) \\ \\
        \qquad\quad\!\!\! + \frac{k_3}{\kappa} \xi \, \ebar \, \rhozerobar^2 + \frac{k_3}{\kappa} \rhozerobar r^{\left(k_3 - k_1\right)} \ebar^2 w_1, \\ \\
        \dfrac{w_1}{\tau} = - r^{k_1} \rhozerobar, \\ \\
        \dfrac{w_2}{\tau} = - r^{k_2 + k_1} \uzerobar \rhozerobar, \\ \\
        \dfrac{\xi}{\tau} = r^{k_3} \ebar.
    \end{cases}
\end{align}

Now, let $k_3 = k_2 = k_1 = 1 = \kappa$. Desingularize once again by multiplying by $r$ and $\rhozerobar$ and divide by $\ebar$ (which corresponds to a rescaling of time). The first equation, for the radial variable, is now $ \dfrac{r}{\tau} = -\uzerobar w_2  r + O(r^2)$. Similar to the case $a<0$, since $\dfrac{w_2}{\tau} = - \rho_0 u$, we once again have
\begin{align*}
    w_2(\tau) =u_L w_1(-\infty) - \left(u \int_{-\infty}^s\rho_0 \d\theta\right) \bigg\vert_{-\infty}^\tau + \int_{-\infty}^\tau \left\{\dfrac{u(s)}{s} \int_{-\infty}^s \rho_0 \d\theta \right\}\d s.
\end{align*}
When $r=0$ and $u=0$ (with $\xi$ equal to the speed of the delta shock, $s_+$), we get
$$w_2(\tau)=u_L \rho_L \left( \u_L \left(1 - \fraction{L}^a\right) - \xi \right),$$
which is negative in the inner layer, because $\xi<\lambda_0(U_L)$ and $u_L<0$, the latter of which is observed in Section \ref{numerics}. We then desingularize by dividing with $-w_2$ and substitute $r = 0$ to get
\begin{align}
    \begin{cases}
        \dfrac{\rhozerobar}{\tau} = -  \rhozerobar \uzerobar, \\ \\
        \dfrac{\uzerobar}{\tau} = -  \uzerobar^2  + 1, \\  \\
        \dfrac{\ebar}{\tau} = -  \uzerobar  \ebar, \\ \\
        \dfrac{w_1}{\tau} = 0 = \dfrac{w_2}{\tau} = \dfrac{\xi}{\tau}.
    \end{cases}
\end{align}
We observe that the system has fixed points at $\uzerobar = \pm 1$ with $\rhozerobar = \ebar = 0$. These correspond to the entry and exit points of the blown-up sphere during the matching procedure of the inner and outer solutions. To study the behavior near these fixed points, we linearize the system around them. At $\uzerobar = 1$, the eigenvalues are $-\uzerobar$, $-2\uzerobar$, and $-\uzerobar$ in the $\rhozerobar$, $\uzerobar$, and $\ebar$ directions, respectively. At $\uzerobar = -1$, the eigenvalues become $1$, $2$, and $1$. In analogy to the construction of Section \ref{GSPT_anegative}, the point \( \uzerobar = -1 \) serves as the entry point into the blown-up sphere along the unstable manifold, while the point \( \uzerobar = 1 \) serves as the exit point along the stable manifold. The nonzero eigenvalues and the full set of eigenvectors ensure normal hyperbolicity, confirming the transversality required to apply the Exchange Lemma. Hence, we conclude that, for sufficiently small \( \e > 0 \), the perturbed trajectory enters the sphere near \( \uzerobar = -1 \) and exits near \( \uzerobar = 1 \), producing a profile connecting the left and right states. Although the current case is not illustrated, it follows the construction of Section \ref{GSPT_anegative} in reverse, as shown in the case illustrated previously.

\begin{remark}
Regardless of the sign of $a$, we have the relation $\rho=\frac{\rho_0}{\e}=\frac{\rhozerobar}{\ebar}$. Therefore, when the inner solution passes through the point $\rhozerobar=1, \uzerobar=0, \ebar=0,$ we observe that the variable $\rho$ becomes unbounded due to division by $\ebar = 0$. This corresponds to a singular behavior in the original physical variables, where the density becomes infinite, signaling the presence of a delta shock and is consistent with the shadow wave approximation that predicts concentration phenomena in the singular limit $\e \to 0$. 
\end{remark}


\section{Conclusion} \label{conclusion}
In this work, we study the Riemann problem for a Keyfitz-Kranzer-type system without source term that models transport dynamics subject to density limitations. A loss of strict hyperbolicity gives rise to non-classical wave phenomena, including transitions through degenerate states.

We analyze self-similar Riemann solutions to this system, identifying both classical wave patterns, such as shocks, rarefactions, and contact discontinuities, and non-classical features, including overcompressive delta shocks. We address the open question of whether classical and non-classical solutions can coexist in this system for various values of the parameter $a \in \mathbb{R} \setminus \{0\}$.

We provide an affirmative answer by identifying the regions in state space where the Riemann problem is resolved non-classically (via delta shocks). 

We rigorously verify that delta shock solutions satisfy the system in the sense of distributions using two approaches. The first involves substituting an ansatz with Dirac delta distributions into the weak formulation and testing against smooth functions. The second approach, known as the shadow wave method, constructs smooth approximations with sharply localized internal layers that converge to singular limits.  Additionally, we explain the underlying singular structure using tools from geometric singular perturbation theory (GSPT), creating a framework for analyzing the formation of singular profiles in systems exhibiting degenerate or non-hyperbolic behavior. The main result of our work is the following theorem:

\begin{theorem}
Let $a\in\mathbb{R}\setminus \{0\}$, and suppose that the left and right states $(\rho_L, u_L)$ and $(\rho_R, u_R)$ satisfy the overcompressive condition \eqref{eq:overcompressive}. Then, for sufficiently small $\e> 0$, the Dafermos-regularized system admits a self-similar viscous profile $(\rho_\e, u_\e)(\xi)$ that connects the left state to the right. In the limit $\e \to 0$, this profile converges in the sense of distributions to a delta-shock solution of the original hyperbolic system, as defined in \eqref{eq:weak_sense_solutions1} and \eqref{eq:weak_sense_solutions2}.
\end{theorem}

The complexity of the wave structure is amplified by the lack of strict hyperbolicity and the presence of degenerate characteristic speeds, as well as by the loss of genuine non-linearity and the changing convexity of the flux. Lastly, we validate our analytical findings using numerical simulations using the Local Lax–Friedrichs scheme. Future work will investigate how these Riemann solutions can serve as building blocks for general initial data in the broader Cauchy problem. 

Additionally, while the LLF method produces reasonable results in areas of ample density, the dissipative nature of the scheme causes numerically ambiguous results as density approaches zero. Going forward, we hope to implement a high-order, non-dissipative numerical scheme to better analyze the behavior of solutions in the presence of vacuums. A promising option is a class of schemes called weighted essentially non-oscillatory (WENO) methods, which offer high-resolution results and can be of third-order or fifth-order accuracy, depending on the type. We are interested in exploring WENO-JS, a specific fifth-order WENO method developed by Guang-Shan Jiang and Chi-Wang Shu. \\

\noindent
{\bf Acknowledgments.} This work is supported by the National Science Foundation under Grant Number DMS-2349040 (PI: Tsikkou). Any opinions, findings, and conclusions or recommendations expressed in this
material are those of the authors and do not necessarily
reflect the views of the National Science Foundation. \\

\noindent
The authors gratefully thank Marko Nedeljkov for suggesting the problem and for his insightful discussions and guidance, which were essential to the development of this work. The authors also thank Griffin Paddock, Camden Toumbleston, and Sara Wilson for sharing their earlier MATLAB code, which served as the foundation for the numerical analysis presented in this paper, and for their valuable input, which facilitated the adaptation of the implementation to the present problem. \\

\noindent
{\bf Data Availability.} The data that support the findings of
this study are available from the corresponding author
upon reasonable request.

\printbibliography
\end{document}